\documentclass{article}

\usepackage{amsmath}
\usepackage{amssymb}
\usepackage[dvipsnames]{xcolor}
\usepackage{comment}
\usepackage{hyperref}
\usepackage{graphicx}
\usepackage[margin=2.5cm]{geometry}
\newtheorem{proposition}{Proposition}
\newtheorem{definition}{Definition}

\newtheorem{theorem}{Theorem}
\newtheorem{lemma}{Lemma}
\DeclareMathOperator{\Op}{Op}

\begin{document}

\title{Fredholm Theory of Non-Elliptic Operators in the Presence of Normally Hyperbolic Trapping}

\author{Selim Amar}

\maketitle
\noindent
Email: selama@stanford.edu \\
Department of Mathematics, Stanford University

\vspace{10pt}

\begin{abstract}

\noindent We present an improved Fredholm theory of non-elliptic operators for when the corresponding classical dynamical system exhibits normally hyperbolic trapping with smooth backward and forward trapped sets. It takes place on coisotropic Sobolev spaces with weak regularity at the backward trapped set $\Gamma_u$, which are roughly speaking made of distributions $v \in H^s$ satisfying $\Phi_u v \in H^{s+1}$ for some quantization $\Phi_u$ of the defining function $\phi_u$ of $\Gamma_u$. We then apply it to the case of wave operators on spacetimes.

\end{abstract}

\tableofcontents

\section{Introduction}

In recent years, non-elliptic Fredholm theory has had striking applications to dynamics and wave equations. For example, it was used to define dynamical resonances for Anosov flows \cite{dyatlov2016dynamical}, to show the meromorphic continuation of the resolvent on conformally compact spaces \cite{vasy2011microlocal} and to provide a framework for analyzing wave equations on asymptotically flat spacetimes \cite{hintz2023linearwavesasymptoticallyflat}. A crucial assumption commonly appeared in these works: the absence of trapping in the classical flow. More precisely, if $P$ is the operator in question, non-trapping roughly means that the corresponding Hamiltonian system $H_p$ of its principal symbol $p$ has global hyperbolic attractors and repellors \textit{on which the frequencies are growing/contracting exponentially}. This exponential behavior of the frequencies, which is often called the redshift/blueshift effect of event horizons in gravitational physics, is crucial to obtaining the required Fredholm estimates. In situations with trapping, the frequency behavior could be replaced by working on spaces of distributions which are growing in time, in which case estimates similar to the non-trapping case held. The purpose of this paper is show that if we allow for a normally hyperbolic trapped set, which roughly speaking corresponds to the case of an hyperbolic set on which frequencies \textit{are not} growing/contracting exponentially, then one can still obtain a well behaved Fredholm theory on \textit{slightly decaying} spaces, albeit on spaces with weaker regularity at the unstable manifold of this trapped set. These spaces are called \textit{coisotropic spaces} and are made up of distributions which are less regular \textit{exactly} at the unstable manifold. The applications we have in mind are to the solvability of wave equations. Before going to the general case, let us first give a simple application of the result to time-independent wave equations to introduce it in a simpler context. 
\newline

Consider a time-independent wave operator $\Box$, say on a compact manifold $N$. Upon taking the Fourier transform in time, we obtain a complex family of differential operators in the spatial variables $R(\sigma)$, $\sigma \in \mathbb{C}$. We wish to invert it so that we can write for $\Box v = f$

\begin{equation}
    v = \int_{Im(\sigma) = C} e^{-i\sigma t} R(\sigma)^{-1} \hat{f}(\sigma)
\end{equation}
If we were to show that $R(\sigma)^{-1}$, which is called the \textit{resolvent family}, is a meromorphic family of Fredholm operators, one could then use a complex contour to compute this integral. By the residue theorem, as the contour is deformed, the poles of $R(\sigma)^{-1}$, which are called \textit{resonances}, would then be picked up. Doing so leads to an expansion of solutions $v$ into \textit{quasinormal modes}

\begin{equation}
    v = \sum_{j} \sum_{k = 0}^{m(\sigma_j) - 1} t^{k} e^{-it\sigma_j} v_{jk}(x) + v_R(x, t)
\end{equation}
where $\sigma_j$ are the resonances with multiplicity $m(\sigma_j)$ of $R(\sigma)^{-1}$, $v_{jk}(x)$ are the modes themselves obtained as residues of $R(\sigma)^{-1}$ and $v_R$ is a reminder term whose decay in time is related to the domain in $\mathbb{C}$ for which we can have estimates on $R(\sigma)^{-1}$. A typical estimate on (high-energy) Sobolev spaces $H^s_{|\sigma|^{-1}}$ \footnote{These are the usual Sobolev spaces but with each derivative coming with a $| \sigma |^{-1}$ factor.} takes place in strips $Im(\sigma) \geq -\mu$ and is of the form

\begin{equation}
    \| v \|_{H^{s}_{|\sigma|^{-1}}} \lesssim |\sigma|^{d-2} \| R(\sigma) v \|_{H^{s-2+d}_{|\sigma|^{-1}}} + o(|\sigma|^{-1}) \| v \|_{_{H^{-N}_{|\sigma|^{-1}}}}
\end{equation}
for $s \in \mathbb{R}$ satisfying some constraints, $N \in \mathbb{R}$ arbitrary and $d > 0$ denoting the loss in derivatives. It implies that for large $Re(\sigma)$, $R(\sigma)^{-1}$ exists with bounds $\| R(\sigma)^{-1} \| = O(|\sigma|^{d-2})$. Furthermore, this \textit{resolvent estimate} leads to the following control on the reminder term $v_R$

\begin{equation}
    \| v_R \|_{s} \lesssim e^{-\mu t} \| f \|_{s-2+d}
\end{equation}
Notice that it controls both the decay and regularity of $v_{R}$. Furthermore, in practice, the modes are nice and smooth, therefore the regularity of $v$ is entirely dependent on the one of $v_R$. This brings us to our first interpretation of the result, which is a slightly improved resolvent estimate from the $O(|\sigma|^{2-2})$ resolvent bound Dyatlov derived in \cite{Dyatlov_2016}. In order not to get too pedantic with the details, the statement is heavily simplified and only the general idea is given, so one should refer to the theorems in sections \ref{EstimateNHT} and \ref{sec3} for precise results.

\begin{theorem}\label{TimeIndWave}(Imprecise; see Theorem \ref{SemiclassicalTrappingEstimate} for the corresponding semiclassical inequality) Let $\Box$ be a time-independent wave operator on a compact manifold N and $v$ satisfy $\Box v = f$. Assume the corresponding null geodesic flow exhibits global normally hyperbolic trapping. Then in a strip $Im(\sigma) > -\nu/2$, where $\nu$ is a dynamical quantity associated to the trapped set, we have the following uniform resolvent estimate for $R(\sigma)^{-1}$

\begin{equation}
    \| v \|_{H^{s}_{|\sigma|^{-1}}} + | \sigma | \| \Phi_u v \|_{H^{s+\alpha}_{|\sigma|^{-1}}} \lesssim |\sigma| \| R(\sigma) v \|_{H^{s-1}_{|\sigma|^{-1}}} + |\sigma|^2\| \Phi_u R(\sigma) v \|_{H^{s+\alpha-1}_{|\sigma|^{-1}}} + o(|\sigma|^{-1}) \| v \|_{H^{-N}_{|\sigma|^{-1}}}
\end{equation}
Where $\Phi_u$ is a quantization of the defining functions of the backward trapped set and $\alpha$ is a certain smooth function satisfying $0 \leq \alpha \leq 1$ and equal to $1$ near the trapped set. In particular, for any $\epsilon > 0$, there is a quasinormal mode expansion up to a $O(e^{-(\nu - \epsilon)/2t})$ reminder term which is coisotropic with respect to the backward trapped set, meaning that it satisfies an estimate of the form

\begin{equation}
    \| v_R \|_{s} + \| \Phi_u v_R \|_{s+\alpha}  \lesssim e^{-(\nu - \epsilon)t/2} (\| f \|_{s-1} + \| \Phi_u f \|_{s+\alpha-1})
\end{equation}
\end{theorem}
\textit{Remark:} The quasinormal mode expansion has already been known for some time now. For example, a polynomial estimate in the same strip gives it, say Dyatlov's $O(|\sigma|^{2-2})$ bound on the resolvent from \cite{Dyatlov_2016} or the polynomial bound of Wunsch and Zworski \cite{Wunsch_2011}. A concrete example where the required assumptions hold is the Kerr-de Sitter spacetime, where the quasinormal mode expansion was first derived in \cite{vasy2011microlocal} using the result of Wunsch and Zworski \cite{Wunsch_2011}. Here, the main point is that rather than having a loss of $d = 2$ due to the trapping, we can get away with an overall loss of $d=1$ by working with the $\Phi_u$ operators directly. We will have more to say about the loss of derivatives later on when discussing the Fredholm theory. Finally, there is an analog of Theorem \ref{TimeIndWave} for appropriate time-dependent wave equations, roughly speaking those where the coefficients of the wave operator also have a partial expansion into quasinormal modes. See for example \cite{AndrasSemilinear} where the method to obtain it is discussed.
\vspace{0.2in}

The above result is not exactly the same as having a robust solvability theory for wave operators. Indeed, one would like to have a theory which behaves well under perturbations. Furthermore, in the time-dependent setting, the result above can be interpreted as a regularity statement about solutions to wave equations, saying that growing solutions with decaying right hand side are decaying as well up to quasinormal modes. For a lot of purposes, it is preferable to have a solvability theory which directly relates the property of $f$ to the ones of $v$, which amounts to studying the invertibility of $\Box$ on different function spaces. Invertibility can be quite difficult to check in some contexts \footnote{For example, finding Feynman propagators on a given spacetime.}, so in the general case showing that $\Box$ is Fredholm is often the best we can do. It turns out that these issues are elegantly addressed by compactifying the problem, for example by compactifying the time coordinate using $\tau = e^{-t}$, and working on spaces of growing/decaying functions in $\tau$. This then leads to studying $\Box$ on a compact manifold $M$ with boundary and trying to show that it is Fredholm on such spaces. The time-dependent perturbations allowed will then be the ones that preserve the relevant dynamical and analytic features of the problem, which in our previous example roughly amounts to small perturbations smooth in $\tau = e^{-t}$. Thus, we move on to introduce the Fredholm theory presented in this paper for a general non-elliptic operator $P$ on a compact manifold with boundary $M$.
\newline

So let $\tau$ be a defining function for the boundary $\partial M$ and assume that $P$ is a \textit{b-differential operator}, which for now can be thought off as an operator smooth in $t = -\log(\tau)$. The corresponding Sobolev spaces on which $P$ acts on are b-Sobolev spaces $H_b ^{s, l}$, which are the usual $H^s$ distributions in $t = -\log(\tau)$ but decaying to order $\tau^l$ at the boundary. The trapped set $\Gamma$ is then a manifold in the characteristic set ${ \rm Char(P) } = p^{-1}(0)$ of $P$ at the boundary $\partial M$, with unstable manifold $\Gamma_u = \phi_u^{-1}(0)$. If $\Phi_u$ denotes a $0-$order quantization of $\phi_u$ and $\mu > 0$ denotes the decay order of $P$ to the corresponding stationary operator $P|_{\tau = 0}$ at the boundary, i.e. $P-P|_{\tau = 0} = O(\tau^{\mu})$ for some stationary operator $P|_{\tau = 0}$ independent of $\tau$, then the coisotropic spaces are defined by

\begin{equation}
    v \in H_b^{s, l} \quad \Phi_u v \in H_b ^{s+\alpha, l} \quad \tau^{\mu} v \in H_b ^{s+\alpha, l} \quad Pv \in H_b ^{s-m+\alpha, l}
\end{equation}
where $m$ is the order of $P$ and $\alpha$ is a certain smooth function satisfying $0 \leq \alpha \leq 1$ and equal to 1 near the trapped set. For now, one can think of $\alpha$ as equal to 1. The reason why a varying $\alpha$ is needed is because the global geometry of $\Gamma_u$ can interact poorly with other features of the problem and prevent the other necessary estimates from holding. Thus, one uses $\alpha$ to effectively 'cutoff' $\Gamma_u$, in the sense that $\alpha$ will be vanishing after a finite time propagation along the Hamiltonian flow of $H_p$ \footnote{It would be interesting to understand if $\alpha = 1$ can be made to work in some cases. For instance in a setting with radial points, that would mostly amount to understanding how $\Gamma_u$ intersects with radial points and how to incorporate $\Phi_u v$ terms in radial point estimates.}. Anyhow, one should thus think of these distributions as in $H_b^{s, l}$ in the zero set of $\phi_u, p, \tau$, which is exactly $\Gamma_u$, but in $H_b^{s+1, l}$ away from it. They are thus less regular on $\Gamma_u$. On these spaces, we will show the following statement (imprecise here).

\begin{theorem}
\label{ImpreciseFredholmTheory}
(Imprecise; see Theorems \ref{FredholmComplexPotential} and \ref{FredholmClosedManifold} for a more precise statement.) Let $P$ be a b-differential operator of order $m$ on a compact manifold with boundary, with real principal symbol $p$. Suppose the Hamiltonian flow $\varphi_t$ of $H_p$ in ${\rm Char}(p)|_{\partial M}$ has a smooth, symplectic \footnote{Symplecticity has to been interpreted in the b-sense when there is a boundary.} trapped set $\Gamma$ which is normally hyperbolic, with smooth coisotropic unstable/stable manifolds $\Gamma_{u/s}$. Then under some conditions on $s$ and $l$ which allows for $l$ to be slightly positive, $P$ is a Fredholm on the coisotropic spaces associated to $\Gamma_u$. Furthermore, the difference in coisotropic regularity between $v$ and $Pv$ is only $m-1$.
\end{theorem}
\textit{Remark:} In the setting of time-independent wave equations on compact manifolds, this result states that under the normally hyperbolic trapping assumption, the wave operator $\Box$ and its perturbations are Fredholm on spaces of distributions which decay in time like $\sim e^{-lt}$ for $l$ small. This will apply for example to the Kerr-de Sitter black holes and simple Minkowski-like spacetimes (see section \ref{SectionExamples}). 
\vspace{0.2in}

The typical conditions on $s, l$ are that $l$ lies outside a discrete set and that $s$ satisfies some constraints depending on $l$. They are thus easy to satisfy. This is not the first Fredholm theory in the presence of trapping that is presented. As discussed above, previous trapping estimates had an extra loss of regularity compared to standard non-elliptic estimates, which brought the total loss of regularity to two. The corresponding Fredholm theory would have to take place on spaces of the form 

\begin{equation}
    v \in H_b^{s, l} \quad Pv \in H_b ^{s-m+2, l}
\end{equation}
The issue with the condition $Pv \in H_b ^{s-m+2, l}$ is that in general it prevents the smooth functions from being dense in this space. Furthermore, arbitrary pseudodifferential operators $A$ do not act on them in the expected way, as you would need some conditions on its principal symbol of $a$. The coisotropic spaces presented here have the advantage of being a priori well-behaved spaces, in the sense that pseudodifferential operators do act on them in the expected way and they contain the smooth functions as a dense subset \footnote{See appendix \ref{AppendixA} for a proof of that statement, along with what makes it work compared to the condition $Pv \in H_b^{s-m+2, l}$.}. This is because the difference in differentiability order between $v$ and $Pv$ is only $m-1$, i.e. the only loss in regularity is due to the non-ellipticity of the problem which contributes a loss of $1$ order as usual. In other words, on these spaces, \textit{there is no regularity loss due to trapping}. Furthermore, the estimates were constrained to work on spaces of functions growing at the boundary. The precise constraint was that the decay order $l$ at the boundary should satisfy

\begin{equation}
    l < -\Tilde{p}_1
\end{equation}
where $\Tilde{p}_1$ is some rescaled version of the subprincipal term of $P$. In particular, unless $\Tilde{p}_1$ is negative, we were forced to work on growing spaces. The main constraint on $l$ for the setup of this paper is 

\begin{equation}
    l < -\Tilde{p}_1 + \frac{\nu}{2}
\end{equation}
where $\nu$ is a dynamical quantity associated to the Hamiltonian flow of $H_p$ at the trapped set, and is positive. Thus, if we can control $\Tilde{p}_1$ in some way, one can work on slightly decaying spaces, which is very useful for the linear stability theory of wave equations for example. The theory presented in this paper is made possible due to the recent estimates of Dyatlov \cite{Dyatlov_2016} and Hintz \cite{Hintz_2021} which play a key underlying role in the present paper. Dyatlov's original purpose was to prove a spectral gap while Hintz's application was to polynomial perturbations of the Kerr family. This and the papers of Dyatlov and Hintz are not the first works on normally-hyperbolic trapping nor on its application to Fredholm theory. This paper is very much in the spirit of \cite{Hintz_2014}, where trapping estimates were found on function spaces localized away from the trapping. These spaces, along with the trapping estimate of Dyatlov \cite{Dyatlov_2016}, were then used in \cite{AndrasSemilinear} to provide a first Fredholm theory for wave operators on the Kerr-de Sitter spacetimes, which was one of the ingredients in the proof of its stability \cite{KerrStab}. The loss of two order of regularity was not an issue for that purpose because smoothing operators were applied. Nonetheless, the Kerr-de Sitter metric was a key motivator for the development of the theory of NHT and it is also one for this paper, where our contribution is mostly to improving its linear theory as it is unclear whether the spaces presented can be made to be an algebra, which would be required for nonlinear purposes. Other relevant articles derived resolvent estimates in varying orders of generality, see for example \cite{Wunsch_2011} for the first proof of normally hyperbolic trapping in Kerr black holes with small angular momentum along with one of the first polynomial resolvent estimates on manifolds with Euclidean ends, \cite{DecayCorrelationTrapping} for a resolvent estimate with logarithmic loss when the unstable/stable manifolds have weak regularity and \cite{Dyatlov_2014} for a Weyl law for the number of resonances between resonance free strips. \cite{Dyatlov_2013} also has a very nice introduction of its role in the analysis of black holes. The reminder of this paper will proceed as follows:

\begin{itemize}
    \item In Section 2, we will provide a brief introduction to the main results and concepts in microlocal analysis that are relevant to the theory of this paper. In particular, we will describe the standard propagation of singularity theorem, radial point estimates and the definition of normally hyperbolic trapping that we will be using. Most of the content is taken from Hintz's notes \cite{PeterNote}.
    \item In section 3, we will state the main estimate of this paper theorems \ref{TheoremB} and \ref{TheoremClosed}, along with the (local) assumptions on $\Gamma$ that are required for it to hold.
    \item In section 4, we will present the proof of the estimate \ref{TheoremB}. It is a modification of Hintz's proof in the cusp setting \cite{Hintz_2014} hence we present it here mainly to highlight that the strengthening we need holds. 
    \item Finally in section 5, we will present the Fredholm theory, which now relies on global assumptions on the Hamiltonian flow of $p$, which we call global normally hyperbolic trapping \ref{GlobalNHT}. The conditions for the Fredholm theory, theorems \ref{TheoremB} and \ref{TheoremClosed}, will be stated there. This is also where we will discuss the standard method used to check them, which is proven in appendix \ref{AppendixB}. We will then apply them to two examples, which will be the wave operators on the Kerr-de Sitter black holes and on a simple Minkowski-like spacetime with trapping.
\end{itemize}

\section{Background Material}

\subsection{Short Review of Microlocal Analysis}
\label{IntroMicrolocalAnalysis}

We start by giving a short review of the important concepts of microlocal analysis for the the Frehdolm theory of non-elliptic operators. We won't delve into the proofs or the precise statements and instead refer to the notes of Hintz \cite{PeterNote} for a thorough introduction. First, by a \textit{pseudodifferential operator} $P \in \Psi^m(M)$ of order $m$ on a manifold $M$, we will mean an operator that in a chart locally looks like

\begin{equation}
    Pu(x) = \int p_{tot}(x, \xi) e^{i(x-x') \cdot \xi} u(y) d\xi dx'
\end{equation}
where $p_{tot}(x, \xi)$ is a symbol of order $m$, meaning it obeys the following bounds

\begin{equation}
\label{StandardSymbol}
    \partial_x^i \partial_{\xi}^j p_{tot}(x, \xi) \leq \frac{C_{ij}}{\sqrt{1+ \| \xi \|^2}^{(m-j)}}
\end{equation}
While $p_{tot}$ is coordinate dependent, its 'top order part', defined as $p_{tot}$ modulo the symbols of order $m-1$, is called the \textit{principal symbol} $p$ and is an invariant function on the cotangent bundle $T^{\ast}M$. It plays a key role as it captures the high frequency behavior of the operator $P$. For differential operators of order $m$, $p$ is homogeneous of degree $m$ in the frequency variable $\xi$. In general and in this paper, we tend to work with operators whose principal symbol is homogeneous since it then restrict, after rescaling, to a function on the cosphere bundle $S^{\ast}M$, which has compact fibers. We will denote by $\Op(\chi) \in \Psi^0$ a quantization of some $0^{th}$ order homogeneous symbol $\chi$, meaning any pseudodifferential operator such that its principal symbol is $\chi$. Furthermore, it can be shown that pseudo-differential operators of order $m$ act on Sobolev spaces continuously

\begin{equation}
    P: H^s \rightarrow H^{s-m}
\end{equation}
There are two main properties of an operator $P$ derived from its principal symbol $p$ that are relevent to this paper. The first is \textit{ellipticity} at some point $(x, \xi) \in T^{\ast}M$, which for homogeneous symbols amounts to non-vanishing of the principal symbol $p$ at $(x, \xi)$. In that case, one can obtain elliptic estimates in a neighborhood of that point, meaning that we can find a pseudodifferential operator $\Op(\chi)$ of order 0 whose principal symbol $\chi$ is elliptic at $(x, \xi)$ and such that we have the following estimate on Sobolev spaces $H^s$ and for an arbitrary error in $H^{-N}$:

\begin{equation}
\label{EllipticEstimates}
    \| \Op(\chi) u \|_s \lesssim \| Pu \|_{s-m} + \| u \|_{-N}
\end{equation}
This takes the form of a Fredholm estimate since the error $\| u \|_{-N}$ is compact relative to $H^s$ for $-N < s$ if we are working on a compact manifold. For manifolds with boundary, we will describe shortly what needs to be done to obtain a compact error. A similar estimate holds for the adjoint near $(x, \xi)$. Thus, we conclude that at points $(x, \xi)$ where our operator is elliptic, we automatically obtain (local) Fredholm estimates. It remains to consider what happens in the characteristic set $p = 0$ where this cannot be done. The main tool for obtaining Fredholm estimates there comes from the formula for the principal symbol of a commutator of two pseudodifferential operators $P$ and $A$ of order $m$ and $l$. It can be shown that $[P, A]$ is of order $m+l-1$, i.e. $P$ and $A$ commutes to first order, and that its principal symbol is 

\begin{equation}
    \sigma_{m+l-1}([P, A]) = \frac{1}{i} H_p(a)
\end{equation}
where $H_p$ denotes the \textit{Hamiltonian vector field} associated to the principal symbol $p$ of $P$

\begin{equation}
    H_p = \sum_i \frac{\partial p}{\partial \xi_i} \frac{\partial}{\partial x_i} - \frac{\partial p}{\partial x_i} \frac{\partial}{\partial \xi_i}
\end{equation}
and $a$ is the principal symbol of $A$. Note that for non-degenerate symbol $p$, meaning $dp \neq 0$ in $p=0$, $H_p(a)$ can be made to be nonzero hence can provide estimates within the characteristic set. The typical way this is used is through \textit{positive commutator estimates}, which involve using the formula 

\begin{equation}
    2 Im \langle P v, |A|^2 v \rangle = \langle i[P, |A|^2]v, v \rangle + \langle \frac{1}{i}(P - P^\dag) |A|^2 v, v \rangle
\end{equation}
with an appropriate $A$ so that $[P, |A|^2]$ is positive and 'large', while the remaining terms can be estimated by $Pv$. Making this work involves some assumptions on $H_p(a)$, which corresponds to the \textit{dynamics} of $H_p$. For the purposes of the Fredholm theory, such methods yields two important estimates. The first, is the \textit{propagation of singularity} theorem, which we recall here coloquially:

\begin{theorem}
\label{PropagationOfSingularities}
(Imprecise Propagation of Singularities; See [Theorem 9.10 \cite{petersen2023waveequationskerrdesitter}] Let $P - iQ$ be a pseudo-differential operator of order $m$ with principal symbol $p - iq$, $p$ and $q \geq 0$ real. Then Sobolev regularity propagates forward along the Hamiltonian vector field $H_p$ of $p$. More precisely, if $\chi$ and $\psi$ are two cutoffs on $S^{\ast}M$ for which the support of $\psi$ is transported to the one of $\chi$ by the forward flow of $H_p$, we have the following estimates on Sobolev spaces for the corresponding quantizations

\begin{equation}
    \| \Op(\chi) u \|_s \lesssim \| \Op(\psi) u \|_s + \| (P-iQ) u \|_{s-m+1} + \| u \|_{-N}
\end{equation}
If $q$ is $\leq 0$ instead, then we propagate backward along $H_p$ instead, meaning that the role of $\chi$ and $\psi$ are reversed.
 
\end{theorem}
\textit{Remark:} Note that relative to the elliptic estimates \ref{EllipticEstimates}, we lose one order of regularity for the norm on $Pu$. This comes from using a positive commutator arguments for the proof and is typical of non-elliptic estimates. This loss of regularity is the reason we cannot obtain a Fredholm theory of non-elliptic operators on standard Sobolev spaces and must modify them instead.
\vspace{0.2in}

The second are \textit{radial point estimates}. They allow us to obtain Fredholm estimates near invariant submanifolds $\mathcal{R} \subset {\rm Char}(p)$ of $H_p$. For that, one needs $\mathcal{R}$ to be a \textit{source}(+) or a \textit{sink}(-) for the flow of $H_p$ in the sense that the flow expands or contract exponentially fast near $\mathcal{R}$. If $p$ is homogeneous of degree $m$ and we let $\phi$ be a quadratic defining function for $\mathcal{R}$ in ${\rm Char}(p) \subset S^{\ast}M$ and some $\rho$ a non-vanishing homogeneous function of degree 1, this means that

\begin{equation}
    \Tilde{H}_p = \frac{1}{\rho^{m-1}} H_p \quad \Tilde{H}_p(\phi) \sim \pm w \phi \quad \Tilde{H}_p(\rho^{-1}) = \pm g \rho^{-1}
\end{equation}
where $w$ and $g$ are smooth functions positive near $\mathcal{R}$. If that is the case, one has the following estimates near the invariant submanifold $\mathcal{R}$. Note that $p_1$ denotes the symbol of the subprincipal term of $P$, i.e.

\begin{equation}
    p_1 = \frac{1}{2i} \frac{\sigma_{m-1}(P - P^{\dagger})}{\rho^{m-1}}
\end{equation}
for $\rho$ as above.

\begin{theorem}(Radial Point estimates; see [Theorem 11.7 \cite{PeterNote}] Suppose $\mathcal{R} \subset {\rm Char}(p)$ is an invariant submanifold of $H_p$ which is a source(+) or sink(-) for its flow. Suppose that the regularity order $s$ satisfies

\begin{equation}
    (s - \frac{m-1}{2})g \mp p_1 > 0
\end{equation}
Let $s_0$ be the infimum of $s$ satisfying this. Then we can estimate $u$ by $Pu$ near $\mathcal{R}$, in the sense that there exists a cutoff $\chi$ on $S^{\ast}M$ elliptic on a neighborhood of $\mathcal{R}$ such that for all $u \in H^{t}$, $t > s_0$, and $-N \in \mathbb{R}$

\begin{equation}
    \| \Op(\chi) u \|_{s} \lesssim \| Pu \|_{s-m+1} + \| u \|_{-N}
\end{equation}
On the other end, if 

\begin{equation}
    (s - \frac{m-1}{2})g \mp p_1 < 0
\end{equation}
Then we can propagate Sobolev regularity from a disjoint neighborhood of $\mathcal{R}$ into it, that is there exists a cutoff $\chi$ as above and a cutoff $\psi$ whose support does not intersect $\mathcal{R}$ and such that

\begin{equation}
    \| \Op(\chi) u \|_{s} \lesssim \| \Op(\psi) u \|_s + \| Pu \|_{s-m+1} + \| u \|_{-N}
\end{equation}
    
\end{theorem}
\textit{Remark:} Note that there is a restriction on the Sobolev order $s$ now. There is two things to note from that: the first is that there may be different radial sets that may come with different requirements on $s$ which we might not be able to satisfy if $s$ is constant. This often requires $s$ to vary on $S^{\ast}M$. The second is that when we work on the dual spaces to prove an estimate for the adjoint $P^{\dagger}$, we usually wish to reverse the direction of propagation along $H_p$. Doing so is compatible with the inequality above, since taking dual simply switches the sign of $s$ and $p_1$.
\vspace{0.2in}

As for the propagation of singularity theorem, radial point estimates come with a loss of one derivative which must be accounted for in the definition of our function spaces. The general idea for the Fredholm theory of non-elliptic operator will then be to obtain some local control 

\begin{equation}
    \| \chi u \|_s \lesssim \| Pu \|_{s-m+1} + \| u \|_{-N}
\end{equation}
somewhere in phase space $S^{\ast}M$ using either elliptic estimates or radial point estimates. Then we use the propagation of singularity theorem to propagates that control everywhere else to obtain a \textit{global Fredholm estimate}

\begin{equation}
    \| u \|_s \lesssim \| Pu \|_{s-m+1} + \| u \|_{-N}
\end{equation}
We finally do the same for the adjoint $P^{\dagger}$ on the dual space, which usually amounts to propagating backward along $H_p$. This works well for compact manifold without boundary, but requires modification when the manifold has a boundary because the error $\| u \|_{-N}$ will then no longer be compact. We discuss this in the next subsection.

\subsection{The b-pseudodifferential operators}
\label{bCalculusIntro}

The applications of our theory we have in mind are to compact manifolds with boundary. Hence, one must have a working theory of pseudodifferential operators on such geometries. There are multiple ways this can be achieved, corresponding to different classes of pseudodifferential operators. In this paper, we will work within the class of $b$-pseudodifferential operators. A formal introduction can be found in \cite{Vasy_2018_Grenoble} and [Chapter 2, 4, 5 \cite{Richard}]. Roughly speaking, if we think in term of compactifying 

\begin{equation}
    [0, \infty)_t \times Y
\end{equation}
where $Y$ is a compact manifold without boundary, then the $b$-pseudo-differential operators correspond to working with the compactification

\begin{equation}
    \tau = e^{-t} \quad [0, 1]_{\tau} \times Y
\end{equation}
In particular, the smooth structure on the original manifold changes to

\begin{equation}
    \partial_t, \partial_y \rightarrow \tau \partial_{\tau}, \partial_y
\end{equation}
hence the b-differential operator are defined as the ones generated by these vector fields, of the form

\begin{equation}
    P = \sum_{j+|k| \leq m} a_{jk}(\tau, y)(\tau \partial_{\tau})^j \partial_{y}^k 
\end{equation}
For the purposes of wave equations, usually $t$ is a timelike function and thus smooth functions on this compactification corresponds to waves which have an expansion into exponentially decaying or growing modes. The general case is defined on a compact manifold with boundary $M$. Let $\tau$ be a defining function for the boundary $\partial M$. Then a b-pseudo-differential operator is roughly speaking an integral operator of the form \footnote{To be precise, one must also impose conditions on the kernel of the operator for them to form an algebra, but colloquially at least this is the idea.}

\begin{align}
    Pu(x) &= \int p_{tot}(e^{-t}, y, \xi) e^{i(x - x') \cdot \xi} u(x') d\xi dx' \\
    &= \int p_{tot}(\tau, y, \sigma, \eta) \left( \frac{\tau}{\tau'} \right) ^{-i \cdot \sigma} e^{i(y - y') \cdot \eta} u(\tau', y') d\sigma d\eta \frac{d\tau'}{\tau'} dy'
\end{align}
where the variables $\sigma$ and $\eta$ are coordinates on covector defined by

\begin{equation}
    \xi = \sigma \frac{d \tau}{\tau} + \eta dy
\end{equation}
The vector bundle over $M$ whose sections are of this form is called the b-cotangent bundle $T^{\ast}_b M$ which plays the role of phase space for b-analysis. Similarly, its quotient by the dilation action $\xi \rightarrow \lambda \times \xi$ of $\mathbb{R}^{+}$ is the $b$-cosphere bundle $S^{\ast}_b M$. Note that it is defined at the boundary $\partial M$ even though $\frac{d\tau}{\tau}$ is not as a covector. Here, the total symbol $p_{tot}$ is \textit{conormal} to the b-cotangent bundle at the boundary $T^{\ast}_b M|_{\partial M}$, meaning that it satisfies

\begin{equation}
    (\tau \partial_{\tau})^i (\partial_y)^j \partial_{\xi}^k p_{tot}\leq \frac{C_{ijk}}{\sqrt{1+\| \xi \| ^2}^{(m-k)}}
\end{equation}
Notice that it agrees with the usual symbols \ref{StandardSymbol} when doing the change of variable $\tau = e^{-t}$. The dependence on $\tau$ of these symbols is typically of the form

\begin{equation}
    p_{tot} \sim \tau^{\beta} \times \log(\tau)^{\alpha} \times r(y, \xi, \sigma)
\end{equation}
for some complex $\beta, \alpha$. Notice that the quasinormal modes $t^k e^{-i\sigma t} = (-\log(\tau))^k \tau^{i \sigma}$ are of this form hence conormal functions come naturally in the study of expansions of waves into quasinormal modes. As before, we have a notion of principal symbol, which is $p_{tot}$ modulo symbols of order $m-1$. It is again invariantly defined, now on the b-cotangent bundle $T^{\ast}_b M$. The analog of the Sobolev spaces are then weighted Sobolev spaces

\begin{equation}
    H^{s, l} = \tau^{l} H^s(M) \quad \| u \|_{s, l} = \| \tau^{-l} u \|_{s}
\end{equation}
and the b-pseudo-differential operators act \ on them in the expected way

\begin{equation}
    P: H^{s, l} \rightarrow H^{s-m, l}
\end{equation}
Furthermore, ellipticity works in the same way as for the standard microlocal analysis, in the sense that if $p$ is elliptic near the support of $\chi$, we have estimates of the form

\begin{equation}
    \| \Op(\chi) u \|_{s, l} \lesssim \| Pu \|_{s-m, l} + \| u \|_{-N, l}
\end{equation}
Notice that the error $\| u \|_{-N, l}$ is not improved in the decay order sense $l$. Thus, it is not a compact error since to obtain a compact inclusion $H^{s, l} \rightarrow H^{s', l'}$, we need both $s' < s$ and $l' > l$. This is typical of estimates obtained from the principal symbol $p$ of $P$. We will describe how to improve the decay order $l$ at the end of the subsection. Just like for the standard pseudo-differential operators, we can also obtain estimates of $u$ in terms of $Pu$ through positive commutator estimates. This yields an analogue of the propagation of singularity theorem and radial point estimates, which we describe here. It should be said that the important dynamics in that case are the ones of $H_p$ at the b-cosphere bundle at boundary $S^{\ast}_b M|_{\partial M}$. Indeed, in typical settings, the integral curves of $H_p$ tend to the boundary hence by using the propagation of singularity theorem one can restrict to working in a neighborhood of the boundary. Near there, we can then apply radial point estimates.

\begin{theorem}(Imprecise Propagation of Singularities for b-pseudo-differential operators; See [Section 5.4.1 \cite{Vasy_2018_Grenoble}] Under the same assumptions on P-iQ, $\chi$ and $\psi$ as the propagation of singularity Theorem \ref{PropagationOfSingularities}, we have

\begin{equation}
    \| \Op(\chi) u \|_{s, l} \lesssim \| \Op(\psi) u \|_{s, l} + \| (P-iQ) u \|_{s-m+1, l} + \| u \|_{-N, l}
\end{equation}

\end{theorem}

\begin{theorem}(Radial Point estimates for b-pseudo-differential operators; see [Section 5.4.8 \cite{Vasy_2018_Grenoble}] \label{RadialPointsB} Suppose $\mathcal{R} \subset {\rm Char}(p) \cap \partial M$ is an invariant submanifold of $H_p$ which is a source(+) or sink(-) for its flow. Suppose that it is also a source(+) or sink(-) in the direction transversal to the boundary, that is

\begin{equation}
    \Tilde{H}_p(\tau) = \pm f \tau \quad f \geq 0
\end{equation}
Suppose that the regularity order $s$ and the decay order $l$ satisfy

\begin{equation}
    (s - \frac{m-1}{2})g + lf \mp p_1 > 0
\end{equation}
and let $s_0$ be the infimum $s$ achieving this. Then we can estimate $u$ by $Pu$ near $\mathcal{R}$, in the sense that there exists a cutoff $\chi$ on $S^{\ast}M$ elliptic on a neighborhood of $\mathcal{R}$ such that for all $u \in H^{t, l}$, $t > s_0$, and all $-N \in \mathbb{R}$

\begin{equation}
    \| \Op(\chi) u \|_{s, l} \lesssim \| Pu \|_{s-m+1, l} + \| u \|_{-N, l}
\end{equation}
On the other end, if 

\begin{equation}
    (s - \frac{m-1}{2})g + lf \mp p_1 < 0
\end{equation}
Then we can propagate Sobolev regularity from a disjoint neighborhood of $\mathcal{R}$ into it, that is there exists a cutoff $\chi$ as above and a cutoff $\psi$ whose support does not intersect $\mathcal{R}$ and such that

\begin{equation}
    \| \Op(\chi) u \|_{s, l} \lesssim \| \Op(\psi) u \|_{s, l} + \| Pu \|_{s-m+1, l} + \| u \|_{-N, l}
\end{equation}

\end{theorem}
\textit{Remarks:} Propagation of singularities and radial point estimates takes the same form as for the standard pseudodifferential operators. Note again that the error $\| u \|_{-N, l}$ is not weaker in the decay order $l$ sense. Furthermore, note that the decay order $l$ now enters in the inequality for radial point estimates, hence $s$ and $l$ are related to each other.

\vspace{0.2in}

There is also another version of radial point estimates worth mentioning for the b-pseudo-differential operators as it appears in the study of rotating black holes.  It concerns invariant submanifold $\mathcal{L}$ which are source or sink in $S^{\ast}_b(M)|_{\partial M}$, but in the direction transversal to the boundary have the opposite behavior. This means that $\Tilde{H}_p(\tau)$ is the opposite sign then for standard radial points

\begin{equation}
    \Tilde{H}_p(\tau) = \mp f \tau
\end{equation}
In that case, we have a similar theorem.
 
\begin{theorem} \label{HyperbolicRadialSets}
(Radial Point estimates for Saddle Radial Points; see [Section Proposition 2.1 \cite{AndrasSemilinear}] Suppose $\mathcal{L} \subset {\rm Char}(p)$ is an invariant submanifold of $H_p$ which is a source(+) or sink(-) for its flow in $S^{\ast}_b(M)|_{\partial M}$ but is a sink(+) or source(-) in the direction transversal to the boundary, as defined above. Suppose that the regularity order $s$ and the decay order $l$ satisfy

\begin{equation}
    (s - \frac{m-1}{2})g - lf \mp p_1 > 0
\end{equation}
with $s_0$ being the infimum $s$ achieving this. Then near $\mathcal{L}$, we can estimate $u$ by $Pu$ and $u$ away from the boundary, in the sense that there exists a cutoff $\chi$ on $S^{\ast}_b M$ elliptic on a neighborhood of $\mathcal{L}$ such that for all $u \in H^{t, l}$, $t > s_0$, and $-N \in \mathbb{R}$

\begin{equation}
    \| \Op(\chi) u \|_{s, l} \lesssim \| Pu \|_{s-m+1, l} + \| \chi_{\tau \neq 0} u \|_{s, l} + \| u \|_{-N, l}
\end{equation}
for some cutoff $\chi_{\tau \neq 0}$ whose wavefront set is disjoint from the boundary $\tau = 0$. On the other end, if 

\begin{equation}
    (s - \frac{m-1}{2})g - lf \mp p_1 < 0
\end{equation}
Then we can propagate Sobolev regularity from a disjoint neighborhood of $\mathcal{L} \cap S^{\ast}_b(M)|_{\partial M}$ into it, that is there exists a cutoff $\chi$ as above and a cutoff $\psi$ whose wavefront set does not intersect $\mathcal{L} \cap S^{\ast}_b(M)|_{\partial M}$ and such that

\begin{equation}
    \| \Op(\chi) u \|_{s, l} \lesssim \| \Op(\psi) u \|_{s, l} + \| Pu \|_{s-m+1, l} + \| u \|_{-N, l}
\end{equation}

\end{theorem}
\textit{Remarks:} Note that these estimates take a similar form as the ones for standard radial sets in Theorem \ref{RadialPointsB}, expect now we have different a priori assumptions on $u$ near the boundary $\partial M$.

\vspace{0.2in}

As the boundaryless case, the main way we obtain Fredholm estimates is by obtaining some a priori control somewhere in phase space $T^{\ast}_b M$, say using radial points, then propagating it everywhere else. This works well to obtain an estimate of the form

\begin{equation}
\label{AlmostFredholmEstimates}
    \| u \|_{s, l} \lesssim \| Pu \|_{s-m, l} + \| u \|_{-N, l}
\end{equation}
As alluded above, something else must be done to improve the error $\| u \|_{-N, l}$ to $\| u \|_{-N, l'}$ for $l' > l$. This is done through the \textit{normal operator} $P|_{\tau = 0}$. In order to define it, we restrict to $P$ being a $b$-differential operator whose coefficients are $O(\tau^{\mu})$ conormal perturbations for $\mu > 0$, meaning it is of the form \footnote{One can also work with a b-pseudodifferential operator whose principal symbol is of the form $p_0 + O(\tau^{\mu})$ for some stationary smooth function $p_0$ on $\partial M$.}

\begin{equation}
\label{DiffOpPlusPerturbation}
    P = \sum_{j+|k| \leq m} (a_{jk}(\tau, y) + b_{jk}(\tau, y))(\tau \partial_{\tau})^j \partial_{y}^k 
\end{equation}
where $a_{jk}$ is smooth (or some b-Sobolev space $H^s_b(M)$ of high enough regularity) and $b_{jk}$ is in

\begin{equation}
\label{DecayingConormal}
    b_{jk} \in \mathcal{A}^{\mu} = \tau^{\mu} \{ v: (\tau \partial_{\tau})^k \partial_{y_i}^l v \in L^2_b \forall k, l \}
\end{equation}
 In that case, its normal operator is obtained by freezing the coefficient at the boundary

\begin{equation}
    N(P) = \sum_{i+|j| \leq m} a_{jk}(0, y)(\tau \partial_{\tau})^j \partial_{y_i}^k 
\end{equation}
The reason to do so is that $N(P)$ is then stationary hence we can study the invertibility properties of $N(P)$ by taking its Mellin transform in $\tau$, or equivalently its Fourier transform in $t = -\log(\tau)$. We then obtain the \textit{indicial family} of operators by substituting $\tau \partial_{\tau}$ by $i\sigma$

\begin{equation}
\label{IndicialFamilyDefinition}
    \widehat{N(P)}(\sigma) = \sum_{j+|k| \leq m} a_{jk}(0, y)(i\sigma)^j \partial_{y}^k 
\end{equation}
Standard microlocal analysis will then typically show that these operators are Fredholm. Indeed, in most cases, these operators will be elliptic almost everywhere. Think for example of the corresponding family for the standard wave operator, which is the spectral family $\Delta - \sigma^2$ of the Laplacian. Unfortunately, there will usually be a characteristic set hence one must employ propagation estimates to show that $\widehat{N(P)}(\sigma)$ is Fredholm. This is discussed in more details in section \ref{ThreeDynamicalSystems}. The next step will be to obtain uniform estimates over $\sigma$ for the family $\widehat{N(P)}(\sigma)^{-1}$. To do so, we will convert this family of operator to a \textit{semiclassical} one, which ammounts to doing the following rescaling for $z$ in a compact subset of $\mathbb{C}$ with $Im(z) = O(h)$

\begin{equation}
    \widehat{N(P)}_h = h^m \widehat{N(P)}(\sigma) \quad \sigma = \frac{z}{h}
\end{equation}
Through the tools of semiclassical analysis, which ammounts to working with pseudodifferential operators depending on a small constant $h$, we will try to obtain \textit{high-energy estimates} for the family $\widehat{N(P)}(\sigma)$ on semiclassical Sobolev spaces as $Re(\sigma) \rightarrow \infty$. For an introduction to semiclassical analysis, see the book by Zworski \cite{SemiclassicalAnalysisZworski}. The main idea is to replace the standard Fourier transform by the semiclassical Fourier transform

\begin{equation}
    \mathcal{F}_h(f)(\xi) = \frac{1}{(2 \pi h)^{n/2}} \int e^{-i x \cdot \xi/h} f(x) dx
\end{equation}
and define quantization of symbols, principal symbols and semiclassical Sobolev spaces $H^s_h$ using $\mathcal{F}_h$ instead of the standard Fourier transform. By doing so, one can obtain analogs of the standard microlocal analysis estimates, but uniformly over the small constant $h$ and with typically $O(h^{\infty})$ reminder errors. This allows us to obtain estimates of the form:

\begin{equation}
\label{HighEnergyEstimates}
    \| u \|_{H^s_{h}} \lesssim h^{-\beta} \| \widehat{N(P)}_h u \|_{H^{s-m+1}_{h}} + o(1) \| u \|_{H^{-N}_h}
\end{equation}
By taking $h \rightarrow 0$, the error term $\| u \|_{H^{-N}_h}$ will then come with a small coefficient hence can be absorbed on the left hand side. This thus gives that for large $Re(\sigma)$, $\widehat{N(P)}(\sigma)$ is invertible. Thus, by the analytic Fredholm theorem, $\widehat{N(P)}(\sigma)$ is invertible everywhere besides a discrete set of poles $\sigma$, called the \textit{indicial roots}. In particular, 

\begin{equation}
    N(P): H^{s, l} \rightarrow H^{s-m+1, l}
\end{equation}
is invertible for all $l$ besides a discrete set determined by $Im(\sigma) = -l$, $\sigma$ a pole of $\widehat{N(P)}(\sigma)^{-1}$. Outside of this discrete set, we can then use $N(P)$ to improve the error term we obtain in our estimates for $u$ in terms of $Pu$. The philosophical reason is that

\begin{equation}
    Pu = N(P)u + O(\tau^{\mu})
\end{equation}
hence inverting $N(P)$ inverts $P$ as well up to a decaying error. The statement we have in the b-category is thus as follows:

\begin{equation}
    \text{Invertibility of $N(P)$} + \text{Estimates of the form \ref{AlmostFredholmEstimates}} = \text{Fredholmness of $P$}
\end{equation}
Hence for the b-pseudo-differential operators we also need an additional step for the Fredholm theory. Nonetheless, it remains very similar to the closed case in the sense that the global dynamics of $H_p$ are the key properties that go into obtaining the estimates. One last thing we'd like to talk about is \textit{trapping}. When the integral curves of $H_p$ all tend to a region of a priori control, say where elliptic estimates or radial points hold, then we can obtain the estimate \ref{AlmostFredholmEstimates}. Unfortunately, there are times where this is not the case: some integral curves stay indefinitely in a region of phase space. In that region, something else must be done to estimate $u$ by $Pu$. The main case it is possible to do so is when the trapped integral curves all lie on some invariant manifold $\Gamma$ of $H_p$ which is \textit{normally hyperbolic}. This paper is about deriving an estimate on that submanifold in order to have Fredholm estimates for $P$ on appropriate function spaces. In the next section, we will introduce normal hyperbolicity and then later on state the microlocal estimates that we can obtain at $\Gamma$.

\subsection{Definition of Normal Hyperbolicity}
\label{IntroTrapping}

The notion of normally hyperbolicity comes from dynamics and generalizes hyperbolic fixed point of vector fields. It concerns invariant submanifolds for which the nearby dynamics are contracting and expanding in different directions. For a thorough introduction to the concept, see \cite{InvariantManifolds}. One can define the notion for any continuous dynamical system $(N, \varphi_t)$, but we will focus on the case of $N = S^{\ast}M$ being the cosphere bundle of the cotangent space $T^{\ast}M$ (or $T^{\ast}_b M|_{\partial M}$) equipped with a (homogeneous of degree 0) Hamiltonian vector field $H$ generating $\varphi_t$. We say it is \textit{normally hyperbolic} at some invariant submanifold $\Gamma$ if the tangent bundle at $\Gamma$ breaks up as

\begin{equation}
\label{NHT}
    T_{\Gamma}N = E_u \oplus E_s \oplus T \Gamma
\end{equation}
for some $d\varphi_t$ invariant sub-bundle $E_u$ and $E_s$ satisfying

\begin{equation}
\begin{split}
    \| d\varphi_t |_{E_s} \| \lesssim e^{-\Tilde{\nu}_s t} \quad \forall t > 0 \\
    \| d\varphi_{-t} |_{E_u} \| \lesssim e^{-\Tilde{\nu}_u t} \quad \forall t > 0
\end{split}
\end{equation}
with respect to some Riemannian metric on $N$ and for constants $\Tilde{\nu}_{u/s} > 0$ measuring the expansion/decay rates respectively. We will denote by $\nu_{u/s}$ the infimum over all such $\Tilde{\nu}_{u/s}$. Thus, the flow expand and contract normally to $\Gamma$ with exponential rate $\nu_{u/s}$. As it stands, this definition is not perturbation invariant hence has limited applicability. A stronger requirement is \textit{r-normally hyperbolic trapping}, which demands that the dynamics on $T\Gamma$ are small in comparison to what happens normally to it:

\begin{equation}
    r \beta < \nu 
\end{equation}
for some $r > 0$ and where $\nu = \min(\nu_u, \nu_s)$ and $\beta$ is the maximal expansion rate on $\Gamma$, i.e. $\| d\varphi_t|_{T\Gamma} \| \lesssim e^{\Tilde{\beta} |t|}$, $\beta = inf(\Tilde{\beta})$. A generalization of the unstable/stable manifold theorem [Theorem 4.1 \cite{InvariantManifolds}] then states that if $\Gamma$ is $r$-normally hyperbolic for every $r$ there exists near $\Gamma$ locally invariant submanifolds $\Gamma_{u/s}$ such that

\begin{equation}
    \Gamma = \Gamma_u \cap \Gamma_s
\end{equation}
and whose tangent spaces at $\Gamma$ are 

\begin{equation}
    T_{\Gamma} \Gamma_u = E_u \oplus T\Gamma \quad T_{\Gamma} \Gamma_s = E_s \oplus T\Gamma
\end{equation}
See theorem \ref{UnstableStableManifoldTheorem} for a more precise statement. Furthermore, these properties are perturbation stable in this case. The unstable manifold will be playing the key role in the Fredholm theory. To make use of these properties for microlocal purposes, one needs to relate it to derivatives $H(a)$ of function on $T^{\ast}(M)$. Here we will assume that $\Gamma_u$ and $\Gamma_s$ are locally level sets

\begin{equation}
    \Gamma_u = \phi_u^{-1}(0) \quad \Gamma_s = \phi_s^{-1}(0)
\end{equation}
for some smooth functions $\phi_{u/s}$ defined in a neighborhood of $\Gamma$ and such that the expansion and contraction behavior is quantified by

\begin{equation}
    H \phi_u = -w_u \phi_u \quad H \phi_s = w_s \phi_s
\end{equation}
for weights $w_u$, $w_s$ positive near $\Gamma$. For higher codimension $\Gamma_{u/s}$, these should be taken as component wise conditions for vectors, in the sense that $w_{u/s}$ are now pointwise positive-definite matrices at $\Gamma$ with respect to the standard inner-product. Finally, we make the further assumption that 

\begin{itemize}
    \item $\Gamma$ is symplectic 
    \item $\Gamma_{u/s}$ are coisotropic and the unstable/stable spaces span the symplectic complements of $T \Gamma$: $E_{u} \oplus E_{s} = T\Gamma^{\perp}$
\end{itemize}
This also implies that $E_{u/s}$ have the same dimension with $E_{u/s}^{\ast} \approx E_{s/u}$ through the symplectic form. Thus, if we let $\rho$ be a non-vanishing function homogeneous of degree 1, $\rho H_{\phi_s^1|_{\Gamma}} \in E_s$ has a symplectic dual in $E_u$ hence we can pick functions $f_i$ defined on $\Gamma$ such that

\begin{equation}
    \sum_i \rho f_i H_{\phi_u^i}(\phi_s^1)|_{\Gamma} = 1
\end{equation}
We then extend them so that $H_{\phi_s^1}(f_j) = 0$ near $\Gamma$ for all $j$. This is possible since $H_{\phi_s^1}$ is transverse to $\Gamma$. These assumptions can be hard to verify. An argument by Dyatlov in [Lemma 5.1 \cite{Dyatlov_2014}] states that in the codimension 1 case with $\Gamma$ orientable, such defining functions $\phi_{u/s}$ exist, and they can be chosen so that the weights $w_{u/s}$ are bounded below by the rates $\nu_{u/s}$

\begin{equation}
\label{LowerBoundWeight}
    w_{u/s} > \nu_{u/s} - \epsilon
\end{equation}
For the higher codimension case, we must first check that defining functions for $\Gamma_{u/s}$ exists, i.e. that they have trivial normal bundles. Then if the resulting weight $w_{u/s}$ is symmetric, Dyatlov's argument goes through and we can average $w_{u/s}$ over the flow to obtain a lower bound as in equation \ref{LowerBoundWeight}. For non symmetric $w_{u/s}$, it fails and we need to assume that $w_{u/s}$ are positive definite. We still denote by $\nu_{u/s}$ their (pointwise now) lower bound. These assumptions on $\Gamma$ are not so hard to check in practice. We refer to Appendix \ref{AppendixB} for a common method that can be applied to do so. In the rest of this paper, the assumptions above will just be referred as \textbf{normally hyperbolic trapping}, referred by NHT. 

\subsection{Comparison between NHT and Radial Sets}
\label{NHTvsRadialSets}

In order to justify the necessity of our estimate for NHT submanifolds, it is relevent to  make a quick comparison between the definition of a radial set $\mathcal{R}$ and a normally hyperbolic trapped set $\Gamma$. In both cases, the behavior of the flow of $H_p$ transversal to them is exponential contraction and expansion in different directions. That is, for defining functions $\psi$ of $\mathcal{R}$ and $\phi_{u/s}$ for $\Gamma_{u/s}$, we have in both cases for the rescaled Hamiltonian vector field

\begin{equation}
    \Tilde{H}_p \phi_{u/s} = \mp w_{u/s} \phi_{u/s} \quad \Tilde{H}_p(\psi) = \pm w_r \psi
\end{equation}
for some positive weights $w_{u/s}$, $w_r$ and a rescaling $\Tilde{H}_p$ of $H_p$. Similarly, for b-pseudo-differential operators, the behavior transversal to the boundary is similar in both cases 

\begin{equation}
    \Tilde{H}_p(\tau) = \pm f \tau
\end{equation}
The crucial difference between the two is the \textit{behavior of frequencies} at the submanifold. For radial sets, they grow or contract exponentially fast, while for NHT trapped set, they have sub-exponential behavior. More precisely, there is an homogeneous of degree 1 non-vanishing function $\rho$ for which

\begin{equation}
    \Tilde{H}_p(\rho) = \mp g \rho \quad g > 0 \text{ on $\mathcal{R}$ vs. } g \text{ possibly $0$ on $\Gamma$} 
\end{equation}
In the case of radial set, it is this weight $g$ that allows us to obtain estimates at $\mathcal{R}$ hence applying radial point estimates to a NHT invariant submanifold $\Gamma$ would lead to degenerate estimates since $g \sim 0$ there. For example, in the b-setting, the inequality between $s$ and $l$ that must be satisfied for saddle radial point estimates is

\begin{equation}
    (s - \frac{m-1}{2})g - lf \pm p_1 > 0
\end{equation}
which for $g = 0$ reduces to

\begin{equation}
   \pm p_1 > lf
\end{equation}
Note that for $p_1 \sim 0$, this means that applying a radial point estimates to a NHT submanifold $\Gamma$ would only yield an estimate for $l < 0$, that is on \textit{growing} spaces. This is an issue since for many applications, say to the stability of wave equations, one wants to work on decaying spaces with $l < \epsilon$ for some small number $\epsilon$. One of the purposes of the microlocal estimate for NHT is to have an estimate on $\Gamma$ which allows for the decay order $l$ to be slightly positive. This is the content of Theorem \ref{TheoremB} and was one of the motivations for this project. Another difference between the two worth pointing out is their geometry: usually, radial sets $\mathcal{R}$ or saddle points $\mathcal{L}$ are Lagrangian manifolds (often conormal bundles), while trapped sets will typically be symplectic manifolds (often cotangent bundles). Thus, they should really be treated as distinct objects even though the dynamics nearby are similar.

\subsection{The Three Different Dynamical Systems}

\begin{figure}[ht!]
\label{TrappingCotangentBundle}
    \centering
    \includegraphics[width=0.7\linewidth]{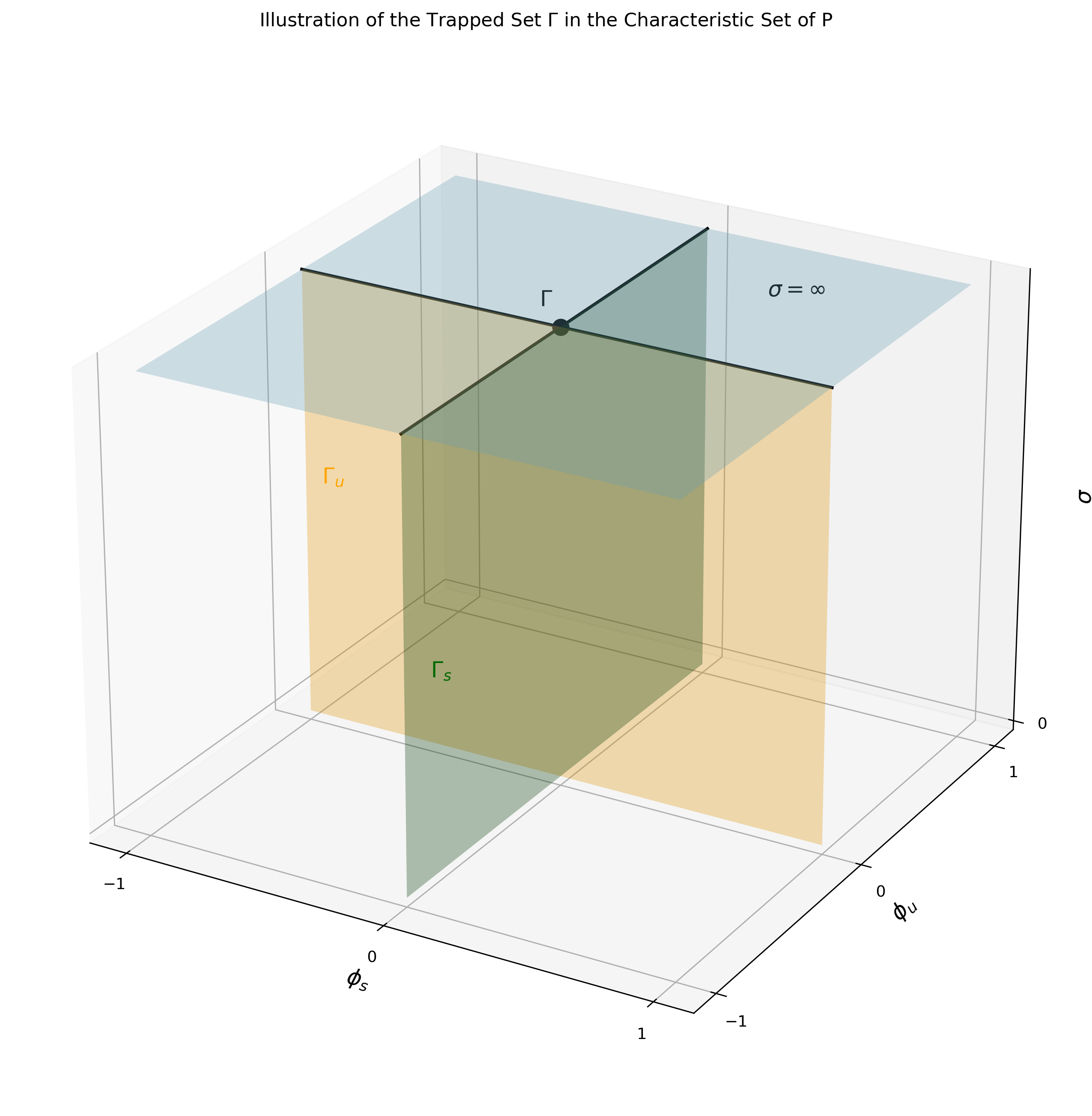}
    \caption{An illustration of the trapped set in the characteristic set $p^{-1}(0) \subset T^{\ast}_b M$. It's intersection with fiber-infinitity $\sigma = \infty$ is $\Gamma$ while it's intersection with finite frequencies $\sigma = const.$ gives the semiclassical trapping. It is assumed to not intersect $\sigma = 0$ hence the 'classical' problem is non-trapping.}
\end{figure}

Here we discuss the three different dynamical systems one encounters when trying to prove that a b-operator $P$ is Fredholm. The first is the Hamiltonian flow of the principal symbol $p$ of $P$. Assuming that all its integral curves tend to the boundary, which is usually the case, propagation and radial point estimates estimates allows us to reduce to studying the dynamics at the boundary. Thus, given a b-differential operator $P$ as in equation \ref{DiffOpPlusPerturbation} with $a_{jk}$ smooth and $b_{jk}$ a $O(\tau^{\mu})$ perturbation for $\mu > 0$

\begin{equation}
    P = \sum_{j+|k| \leq m} (a_{jk}(\tau, y) + b_{jk}(\tau, y))(\tau \partial_{\tau})^j \partial_{y}^k 
\end{equation}
we need to consider the Hamiltonian of 

\begin{equation}
    p|_{\tau = 0} = i^m \sum_{j+|k| = m} a_{jk}(0, y) \sigma^j \eta^k
\end{equation}
In the \textit{non-trapping} case, where all its integral curves at the boundary tend to regions where we can estimate $u$ by $Pu$, say radial sets, we then have a global estimate of the form \ref{AlmostFredholmEstimates}. In this paper, we will generalize that and assume that there is normally hyperbolic trapping, that is there exists some invariant submanifold $\Gamma$ on which the integral curves stay indefinitely. Note that the Hamiltonian satisfies

\begin{equation}
    H_p(\sigma)|_{\tau = 0} = 0
\end{equation}
where $\sigma$ is the dual variable to $\tau \partial_{\tau}$, hence it is tangent to the $\sigma = cont.$ hyperplanes. It will thus be natural to rescale the flow by $\sigma$ where it does not vanish in order to make it homogeneous of degree 0, i.e. consider the vector field

\begin{equation}
    \Tilde{H}_p = \frac{1}{\sigma^{m-1}} H_p
\end{equation}
which is then a vector field on the b-cosphere bundle $S^{\ast}_b M$. Assuming $\sigma \neq 0$ on $\Gamma$, we can thus consider $\Gamma$ as a subset of $S^{\ast}_b M$ and we assume it is NHT with respect to $\Tilde{H}_p$ as defined in section \ref{IntroTrapping}. We think of it as living in $\sigma = \infty$. We will show that one can still obtain an estimate analogous to \ref{AlmostFredholmEstimates}, except on more degenerate function spaces. Afterwards, the next step is to work with the indicial family $\widehat{N(P)}(\sigma)$ \ref{IndicialFamilyDefinition} of the normal operator $N(P)$ in order to invert it. The standard approach is to show that the assumptions of the analytic Fredholm theorem are satisfied which would give that $\widehat{N(P)}(\sigma)^{-1}$ is a meromorphic family of operators. In order to show that $\widehat{N(P)}(\sigma)$ is Fredholm, one typically has to use non-elliptic estimates since $\widehat{N(P)}(\sigma)$ will not be elliptic everywhere. Its principal symbol is the restriction of $p$ to the $\{ \sigma = 0 \}$ hyperplane:

\begin{equation}
    p|_{\tau = 0, \sigma = 0} = i^m \sum_{|k| = m} a_{jk}(0, y) \partial_{y}^k
\end{equation}
hence the corresponding dynamical system is the Hamiltonian flow of $\Tilde{H}_p$ in $p^{-1}(0) \cap S^{\ast}_b M|_{\partial M} \cap \{ \sigma = 0 \}$. It is often called the \textit{classical} flow. In practice, this system will be non-trapping for a judicious choice of $\tau$, meaning that the trapped set $\Gamma$ does not intersect $\{ \sigma = 0 \}$ \footnote{The theory also works in the case where there is trapping in $\sigma = 0$, but since our examples don't have it we don't cover it.}. Thus, using the non-trapping theory, one can show $\widehat{N(P)}(\sigma)$ is Fredholm for any $\sigma$. Finally, to show that it is invertible for some $\sigma$, we must derive high-energy estimates of the form \ref{HighEnergyEstimates}. This ammounts to working with the semiclassical principal symbol $p_{h}$ of $\widehat{N(P)}_h$ for $\sigma = \frac{z}{h}$:

\begin{equation}
    p_h = i^m \sum_{j+|k| = m} a_{jk}(0, y) z^j \eta^k
\end{equation}
Note that this is just 

\begin{equation}
    p_h = p(y, \eta, \sigma = Re(z))
\end{equation}
hence the corresponding dynamical system is the flow of $\Tilde{H}_p$ on the finite frequencies $\sigma = Re(z)$ hyperplane. This will also exhibits NHT trapping and necessitate the trapping estimate. These different dynamics are summarized in table \ref{ThreeDynamicalSystems} while the trapped set is illustrated in figure \ref{TrappingCotangentBundle}.

\begin{table}[ht!]
\label{ThreeDynamicalSystems}
    \centering
    \begin{tabular}{c|c|c|c}
        & $P$ & $\widehat{N(P)}(\sigma)$ & $\widehat{N(P)}_h$  \\
         Dynamical System & $H_p|_{S^{\ast}_b M|_{\partial M}}$ & $H_p|_{S^{\ast}_b M|_{\partial M}} \cap \{ \sigma = 0 \}$ & $H_{p}|_{T^{\ast }_b M|_{\partial M} \cap \{ \sigma = Re(z) \}}$ \\
         Trapping & Yes & No & Yes
    \end{tabular}
    \caption{A summary of the different dynamical systems and trapping dynamics we need to deal with when proving that a b-operator $P$ is Fredholm.}
\end{table}

\subsection{Relationship with the Decay of Waves}

\begin{figure}[ht!]
    \centering
    \includegraphics[width=0.7\linewidth]{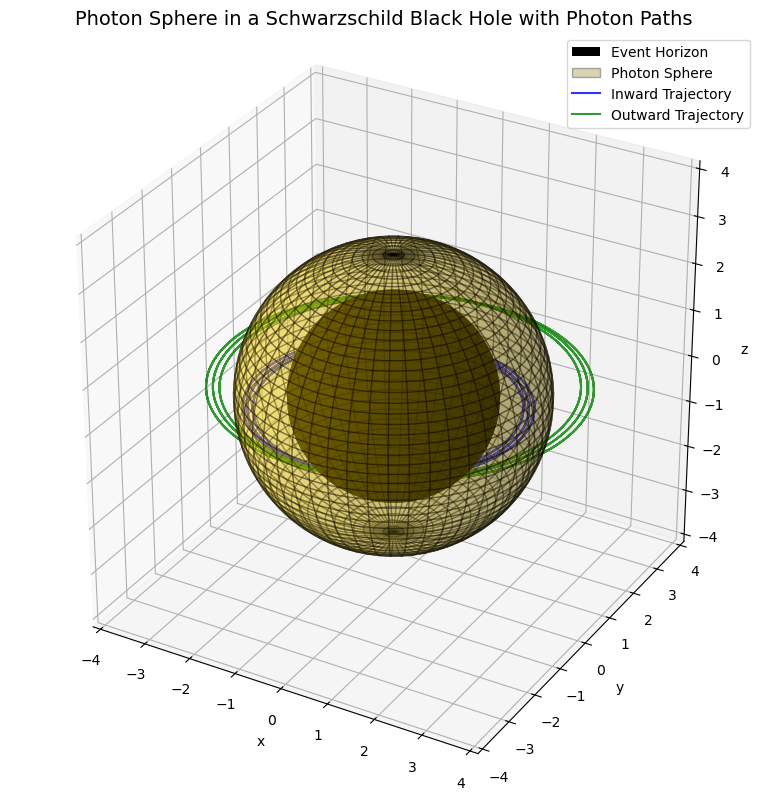}
    \caption{An illustration of the photon sphere (in yellow) around the event horizon of a black hole (in black). Photons starting their path tangentially to the photon sphere will orbit on it forever, while the other ones either escape to infinity or spiral towards the black hole. The behavior of photons normal to the photon sphere is normally hyperbolic. The trajectories shown in this figure are not physical but are shown there for illustrative purposes.}
    \label{fig:PhotonSphere}
\end{figure}

We now give a quick explanation of why NHT plays such a central notion in the linear theory of waves on black holes spacetime. Non-rotating black holes are described by the Schwarzschild metric:

\begin{equation}
    ds^2 = -(1-\frac{r_s}{r})dt^2 + (1 - \frac{r_s}{r})^{-1}dr^2 + r^2 d\Omega^2_{\mathbb{S}^2}
\end{equation}
In other to understand the decay of waves on such a background, a first step is to check whether the geodesic flow is dissipative, that is whether geodesics escape the exterior black hole region. Unfortunately, it is not in this case; there exists a sphere, called the \textit{photon sphere} fig. \ref{fig:PhotonSphere}, from which null geodesics never leave. More precisely, the cotangent space of the submanifold $r = \frac{3}{2} r_s$ is invariant under the null geodesic flow:

\begin{equation}
    \Gamma = T^{\ast}S \cap \{ g^{-1}(\xi, \xi) = 0 \} \quad S = \{ r = \frac{3}{2} r_s\}
\end{equation}
It is thus an obstacle to obtaining decay estimates. The saving grace is that one can check that the geodesic flow is $r$-normally hyperbolic for every $r$ at $S$ hence one can hope to use that to obtain estimate there. This is indeed possible and is a key ingredient in its linear stability theory. A similar behavior is observed for rotating black holes, but the analogue of the photon sphere does not project to physical space as nicely. See section \ref{SectionExamples} for a more technical description of how our estimate for NHT is applied to the study of wave operators.

\section{Microlocal Estimate for NHT}
\label{EstimateNHT}
\subsection{The Standard Microlocal Estimate}

We now present how to obtain a microlocal estimate at an invariant manifold $\Gamma$ satisfying NHT. For the necessary microlocal techniques we will use, we refer to our reminder in section \ref{IntroMicrolocalAnalysis} and the excellent notes of Hintz \cite{PeterNote}. We will first be considering a pseudodifferential operator $P \in \Psi^m$ on a compact manifold $M$ equipped with some volume density. We assume that the principal symbol $p$ is real, homogeneous of degree $m$ and non-degenerate (i.e. $dp \neq 0$ on $p^{-1}(0) \subset T^{\ast}M \backslash 0$). Furthermore, we assume that $\Gamma \subset {\rm Char}(P) = p^{-1}(0) \subset S^{\ast}M$ is an invariant manifold of the Hamiltonian flow on the cosphere bundle $S^{\ast}M$ satisfying NHT, meaning that the induced vector field on $S^{\ast}M$

\begin{equation}
    \Tilde{H}_p = \frac{1}{\rho^{m-1}} H_p
\end{equation}
has a normally hyperbolic invariant submanifold $\Gamma$ there as defined in section \ref{IntroTrapping}. Here $\rho$ is a smooth homogeneous function on $T^{\ast}M$ of degree 1, nonvanishing in a neighborhood of $\Gamma$ and such that

\begin{equation}
    \Tilde{H}_p \rho = g \rho \quad \text{in a neighborhood of $\Gamma$} 
\end{equation}
for some weight function $g$. We extend the defining functions $\phi_u$, $\phi_s$ and the weights $w_u, w_s$ to be functions on $T^{\ast}M\backslash0$ homogeneous of degree 0. This means that these extensions satisfy

\begin{align}
    \Tilde{H}_p \phi_{u/s} &= \mp w_{u/s} \phi_{u/s} + r_{u/s} \frac{p}{\rho^m} \\
    \Gamma_{u/s} &= \left\{ p = \phi_{u/s} = 0 \right\}
\end{align}
for some error $r_{u/s}$ homogeneous of order $0$. We recall that as a consequence of these assumptions, $\Phi_u v$ satisfies a better system of equations:

\begin{equation}
    (P+\frac{1}{i} W_u)\Phi_u v = \Phi_u P v + R_1 P v + R_2v
    \label{DampingEquation}
\end{equation}
where $W_u \in \Psi^{m-1}$ is a matrix quantization of the weight $w_u \rho^{m-1}$, $\Phi_u$ is a self-adjoint one of order $0$ for the defining function of the unstable manifold $\Gamma_u$ (which is vectorial $\Phi_u = \Phi_{u,j}$ if $\Gamma_u$ is of higher codimension than 1) and $R_1 \in \Psi^{-1}$, $R_2 \in \Psi^{m-2}$ are reminder errors. The term $W_u$ plays the role of damping hence allows for stronger estimates outside the unstable manifold. The symbol of the subprincipal term of $P$ will also play a role, which we define by

\begin{equation}
    p_1 = \frac{1}{2i} \frac{\sigma_{m-1}(P - P^{\dagger})}{\rho^{m-1}}
\end{equation}
Under these assumptions, a slight strengthening of Dyatlov's \cite[Theorem 1]{Dyatlov_2016} and Hintz's \cite[Theorem 3.1]{Hintz_2021} gives the following estimate near $\Gamma$.

\begin{theorem}
\label{TheoremClosed}
Suppose $P \in \Psi^m(M)$ has an invariant submanifold $\Gamma \subset p^{-1}(0)$ which is NHT in $S^{\ast}M$. Assume that $p_1$ and the regularity order $s$ satisfy

\begin{equation}
\begin{split}
    &\left(2 \nu_u - 2p_1  - (2s - m + 3)g \right)|_{\Gamma} > 0 \\
    &\left(\nu_s - 2p_1  - (2s - m + 1)g \right)|_{\Gamma} > 0 
\end{split}
\end{equation}
Let $s_0$ be the infimum over all $s \in \mathbb{R}$ satisfying the two conditions above. Then for any neighborhood $U$ of $\Gamma$ and number $N$, there are cutoffs $\Op(\chi) \in \Psi^0(M)$ elliptic at $\Gamma$, $\Op(\chi_{\phi_u \neq 0}) \in \Psi^0(M)$ with wave front set disjoint from a neighborhood of $\Gamma_u$ and $G, E \in \Psi^0$ elliptic on $WF(\Op(\chi))$, with all their wave front sets contained in U,  s.t the following inequality holds for all $v \in H^{-N}$ such that $WF^{t}(v) \cap \Gamma = \emptyset$, $t > s_0$:

\begin{equation}
\begin{split}
      &\| \Op(\chi) v \|_{s} + \sum_{j=1}^{codim(\Gamma_u)} \| \Op(\chi) \Phi_{u,j} v \|_{s+1} \\ 
      & \qquad \lesssim \| G P v\|_{s-m+1} + \sum_{j=1}^{codim(\Gamma_u)} \| G \Phi_{u,j} P v \|_{s-m+2} + \| \Op(\chi_{\phi_u \neq 0}) v \|_{s+1} + \| E v \|_{-N}
\end{split}
\end{equation}
The same result holds for the adjoint $P^{\dag}$ but with $\phi_s$ replacing $\phi_u$ and $\nu_s$ switched with $\nu_u$.
\label{theorem1}
\end{theorem}
\textit{Remark:} Under the same assumptions, Hintz's estimate (Dyatlov's estimate is in the semiclassical setting but is of similar form) takes the following form

\begin{equation}
    \| \Op(\chi) v \|_s \lesssim \| GPv \|_{s-m+2} + \| \Op(\chi_{\phi_u \neq 0}) \|_{s+1} + \| E v \|_{-N}
\end{equation}
The problem with this estimate for the purpose of Fredholm theory is that the loss of regularity compared to elliptic estimate is 2, i.e. we have a $H^s$ norm on the left handside and a $H^{s-m+2}$ norm on the right. This is problematic because we would have to work on spaces of the form $\{ v \in H^s: Pv \in H^{s-m+2} \}$, which are not well behaved as explained in the introduction. So the novelty in the one presented above is the inclusion of the $\| \Phi_u v \|_{s+1}$ and $\| \Phi_u P v \|_{s-m+2}$ terms, which allows the difference in regularity order between $v, \Phi_u v$ and $Pv, \Phi_u Pv$ to be just 1. The corresponding Fredholm theory then takes place on coisotropic spaces, which are well behaved.
\vspace{0.2in}

A simple example of such an operator on a closed manifold is given by Hintz in \cite[Example 3.5]{Hintz_2021} as 
\begin{equation}
    P = -i(\partial_y + \sin(x) \partial_x)
\end{equation}
on the 3-torus $\mathbb{S}^1_x \times \mathbb{S}^1_y \times \mathbb{S}^1_z$, with unstable and stable manifold
\begin{equation}
    \phi_u = \frac{\xi_x}{|\rho|} \quad \phi_s = x  \quad \rho = \xi_z \quad \quad \xi = \xi_x dx + \xi_y dy + \xi_z dz
\end{equation}
Thus, for that case, the estimate says that the $\partial_x$ derivative of $v$ has one more order of regularity than expected.
\newline

Now for the application to waves, one must work on manifolds with boundaries. If we want to work with expansions in terms of exponentially decaying quasinormal modes, this ammounts to working in the $b$-calculus. For a very short introduction to the $b$-calculus we refer to section \ref{bCalculusIntro}, while for a more in-depth treatment to \cite{Vasy_2018_Grenoble} and [Chapter 2, 4, 5 \cite{Richard}]. The crucial distinction is that the important dynamics are now the ones of $\Tilde{H}_p$ on the b-cosphere bundle at the boundary $S_{b}^{\ast}M|_{\partial M}$. More precisely, we will let $P$ now be a b-pseudodifferential operator with symbol decaying to a 'stationary' state at $\partial M$, that is there exists some homogeneous $p_0 \in C^{\infty}(T^{\ast}_{b}M|_{\partial M} \backslash 0)$ such that in a tubular neighborhood of the boundary

\begin{equation}
    p - p_0 \in \mathcal{A}^{\mu}(M)
\end{equation}
where $\mathcal{A}^{\mu}$ is the space of conormal distributions at $\partial M$ decaying with rate $\mu > 0$ as described in section \ref{bCalculusIntro} equation \ref{DecayingConormal}. We assume that $\Tilde{H}_{p_0}$ has a NHT invariant submanifold $\Gamma$ in $p_0^{-1}(0) \subset S_{b}^{\ast}M|_{\partial M}$ as described in section \ref{IntroTrapping}. In order to carry out the analysis, one must then use functions defined on the whole of $T^{\ast}_{b}M$. We thus need extensions $\Bar{\phi}_u$, $\Bar{\phi}_s$ and $\Bar{\rho}$ away from the boundary. In the case of $r$-normally hyperbolic trapping for every $r$, this is provided by the following dynamical result of Hintz [Theorem 2.6 \cite{Hintz_2021}]. In plain words, it says that the unstable and stable manifold $\Gamma_u$ and $\Gamma_s$ at the boundary $\tau = 0$ have an extensions to a tubular neighborhood of the boundary $\tau < \epsilon$.

\begin{theorem} (Hintz \cite[Theorem 2.6]{Hintz_2021}) Assume $\Tilde{H}_{p_0}$ has an invariant submanifold $\Gamma \subset p_0^{-1}(0)$ at the b-cosphere bundle at the boundary $S^{\ast}_b M|_{\partial M}$ which is $r-$normally hyperbolic for every $r$. Let $\Gamma_{u/s}$ be its unstable/stable manifold with defining functions $\phi_{u/s}$. Then there exist functions $\Bar{\phi}_{u/s}$ defined on a tubular neighborhood of the boundary such that $\Tilde{H}_p \Bar{\phi}_{u/s} = \Bar{w}_{u/s} \Bar{\phi}_{u/s}$ in ${\rm Char}(p)$ and

\begin{equation}
    \Bar{\phi}_{u/s} - \phi_{u/s} = O(\tau^{\mu})
\end{equation}
    
\end{theorem}
\textit{Remark:} The extension of the stable manifold $\Bar{\phi}_{s}$ away from the boundary is unique, but the one of the unstable manifold $\Bar{\phi}_u$ is not. Intuitively, one can see that from how $\Tilde{H}_{p_0}$ acts on the boundary defining function $\tau$. We shall discuss later that $\Tilde{H}_{p_0}(\tau) = -f \tau$ for some positive function $f$, so $\tau$ acts like a defining function for the unstable manifold $\Gamma_u$ as well. This means that $\Bar{\phi}_u$ is only defined up to $O(\tau^{\mu})$ perturbations, hence is not unique. In particular, for our Fredholm theory, $\Bar{\phi}_u$ and $\tau^{\mu}$ will need to be treated equivalently.
 
\vspace{0.2in}

In the general case, one must take for granted those extensions. The assumptions in the $b$-setting are thus as follows:

\begin{itemize}

    \item Pick a homogeneous smooth function $\rho$ of degree 1 on $T_{b}^{\ast}M|_{\partial M}$ nonvanishing near $\Gamma$. Define $g$ by $\Tilde{H}_{p_0} \rho = g \rho$. In the case of $\tau \partial_{\tau}$ being elliptic on $\Gamma$, which it will be for applications, one can assume $g = 0$ by choosing its dual coordinate.
    
    \item There exists smooth defining functions $\phi_u$, $\phi_s$ defined in a neighborhood of $\Gamma \subset {\rm Char}(p_0)$ satisfying $\Tilde{H}_{p_0} \phi_u = -w_u \phi_u$, $\Tilde{H}_{p_0} \phi_s = w_s \phi_s$ and such that $\Gamma = \phi_u^{-1}(0) \cap \phi_s^{-1}(0)$. We assume they extend on a tubular neighborhood of the boundary \footnote{It is not strictly necessary to assume an extension for $\phi_u$ since we can estimate $\tau^{\mu}$ errors, but the extension for $\phi_s$ is necessary for the argument hence why we include this assumption. In practice, $r$-normal hyperbolicity is required to check that, and it gives the extension for both $\phi_u$ and $\phi_s$ anyway.} to $\Bar{\phi}_{u/s}$ with $\Bar{\phi}_{u/s} - \phi_{u/s} = O(\tau^{\mu})$ and $\Tilde{H}_p \Bar{\phi}_{u/s} = \mp \Bar{w}_{u/s} \Bar{\phi}_{u/s} + r_{u/s}p$ with $\Bar{w}_{u/s} - w_{u/s} = o(1)$. Fix such extensions.
    
    \item In a tubular neighborhood of the boundary, the Hamiltonian flow converges to the boundary $\partial M$ near the trapped set, that is $\Tilde{H}_{p_0} \tau = - f \tau$ with  $\inf_{\Gamma} f > 0$.

\end{itemize}
Under these assumptions, we obtain the analog of Theorem \ref{TheoremClosed} in the b-setting. Furthermore, to ease the notation, we will drop the overhead bar from extensions of $\phi_{u/s}$, $\rho$ and $w_{u/s}$ and instead denote them by the same symbol.

\begin{theorem}
\label{TheoremB}

Let $P \in \Psi_b^m(M)$ be a b-pseudo-differential operator whose principal symbol at the boundary $p_0$ has an invariant submanifold $\Gamma \subset p_0^{-1}(0)$ satisfying NHT. Suppose that its subprincipal term $p_1$, the regularity order $s$ and the decay order $l$ satisfy at the boundary $\tau = 0$

\begin{equation}
\begin{split}
    &\left( 2 \min(\nu_u, f\mu) - 2p_1 - 2lf - (2s - m + 3)g \right)|_{\Gamma} > 0 \\
    &\left( \nu_s - 2p_1 - 2lf - (2s - m + 1)g \right)|_{\Gamma} > 0
\end{split}
\end{equation}
where $f, g, \mu$ are as described above. Let $s_0$ be the infimum over all $s \in \mathbb{R}$ satisfying the two conditions above. Then for any neighborhood $U$ of $\Gamma$ and number $N$, there are cutoffs $\Op(\chi) \in \Psi^0_b(M)$ elliptic at $\Gamma$, $\Op(\chi_{\phi_u \neq 0}) \in \Psi^0_b(M)$ with wave front set disjoint from a neighborhood of $\Gamma_u$, $\Op(\chi_{\tau \neq 0})$ with wave front set disjoint from a neighborhood of the boundary $\partial M$ and $G, E \in \Psi^0_b(M)$ elliptic on $WF(\Op(\chi))$, with all their wavefront set contained in $U$, s.t the following inequality holds for all $v \in H^{-N, l}$ with $WF^{t, l}(v) \cap \Gamma = \emptyset$, $t > s_0$:

\begin{equation}
\begin{split}
      &\| \Op(\chi) v \|_{s, l} + \sum_{j=1}^{codim(\Gamma_u)} \| \Op(\chi) \Phi_{u, j} v \|_{s+1,l} + \| \Op(\chi) \tau^{\mu} v \|_{s+1,l} \\
      &\qquad \lesssim \| G P v\|_{s-m+1, l} + \sum_{j=1}^{codim(\Gamma_u)} \| G \Phi_{u, j} P v\|_{s-m+2, l} + \| G \tau^{\mu} P v\|_{s-m+2, l} \\
      &\qquad \hspace{4em}+ \| \Op(\chi_{\phi_u \neq 0}) v \|_{s+1, l} + \| \Op(\chi_{\tau \neq 0}) v \|_{s+1, l} + \| E v \|_{-N, l}
\end{split}
\end{equation}
The same results holds for the adjoint $P^{\dag}$ with $\phi_s$ replacing $\phi_u$, $\nu_s$ switched with $\nu_u$ and with the term $\chi_{\tau \neq 0}$ dropped. Finally, the same conclusion holds but with the $\chi \tau^{\mu} v$ and $G \tau^{\mu} P v$ terms dropped on both side assuming only
 
 \begin{equation}
 \begin{split}
     &\left( \min(\nu_s, 2 \nu_u) - 2p_1 - 2lf - (2s - m + 3)g \right)|_{\Gamma} > 0 \\
         &\left( \nu_s - 2p_1 - 2lf - (2s - m + 1)g \right)|_{\Gamma} > 0
\end{split}
 \end{equation}
i.e. with no restriction on $l$ coming from $\mu$.
\end{theorem}
 \textit{Remark:} Note that this takes place near the boundary and that this statement is independent of the extensions outside $\partial M$ of $\phi_{u}$ we picked if we include the $\tau^{\mu}$ term. On the other end, if we do not, the estimate depends on the extension beyond the boundary of $\phi_u$ we chose. These are non-unique in general and depends on the particular perturbation $p - p_0$ we are dealing with, which isn't ideal. Thus, the estimate without the $\tau^{\mu}$ term is mainly useful in the stationary setting, where $p = p_0$.
 \vspace{0.2in}
 
 Examples of operators satisfying these assumptions will be given in section \ref{SectionExamples}. Combining the constraints from $P$ and $P^{\dag}$ and setting $g = 0$, the resulting condition on the decay rate $l$ is

 \begin{equation*}
     l < \min(\frac{\nu}{2f}, \mu) - p_1
 \end{equation*}
 where $\nu$ is the minimum of the contraction rate on $\Gamma_s$ and expansion rate on $\Gamma_u$ \footnote{This condition is exactly what Dyatlov found in \cite{Dyatlov_2016} for his spectral gap, so it should come as no surprise that it shows up here.}. Thus, the stronger the normal hyperbolic behavior is, the larger we can allow the decay to be and still obtain Fredholm estimates. Furthermore, this solves the problem mentionned in section \ref{NHTvsRadialSets} with radial point estimates applied to $\Gamma$ since now we have an estimate on slightly decaying spaces. The proof will be given in the next section \ref{ProofNHTEstimate} and is a small variation of Hintz's proof \cite{Hintz_2021} in the cusp setting \footnote{The cusp algebra is another pseudodifferential calculus one can use on a manifold with boundary. It ammounts to working with polynomial decay rather than exponential. }. The main difference is that we obtain restriction on $l$ from the terms $\Tilde{H}_p(\tau)$ and we can obtain control on $\tau^{\mu} v$ also rather than just $\phi_u v$. As usual in microlocal analysis, this is a commutator estimate but it uses three of them. We will present the $b$-case; for a closed manifold simply drop the $\tau$ terms. To simplify the notation, we will use $|y| \sim x$ to denote 'of order $x$', i.e. $y \in [-ax, bx]$ for some constants $a$, $b$. Similarly, $\lesssim$ and $\simeq$ denotes $\leq$ and $=$ up to constants. We typically use the symbol $\chi$ and $\psi$ to denote microlocal cutoffs.

 \subsection{The Semiclassical Estimate}

 As explained in section \ref{bCalculusIntro}, for our application to the Fredholm theory of $b$-operators, it is important to have a semiclassical version of the estimate, i.e. an estimate that applies to operators depending on a small constant $h$ which is uniform over $h$. For an introduction to semiclassical analysis, we refer to the book by Zworski \cite{SemiclassicalAnalysisZworski} as a reference. These estimates will be used in order to obtain 'high-energy estimates' for the normal operators by setting $h \sim 1/Re(\sigma)$. This was already obtained in the original paper of Dyatlov \cite{Dyatlov_2016} with a loss of $h^{-2}$ relative to elliptic estimates, but it will follow from the same proof as given below that a slightly stronger estimate holds if we include the operators $\Phi_u$ in the estimate. Roughly speaking, to convert to the semiclassical estimate, each loss of derivatives now comes with a power of $h^{-1}$ and we switch the Sobolev spaces used to the semiclassical ones $H^s_h$, i.e. based on the semiclassical Fourier transform. More precisely, we now assume that $P_h$ is a semiclassical pseudo-differental operator of order $m$ whose principal symbol $p_h$ is such that $H_{p_h}$ has a NHT invariant submanifold $\Gamma$. In that case, one has the following estimate.

 \begin{theorem} \label{SemiclassicalTrappingEstimate} Let $P_h$ be a semiclassical operator of order $m$ whose real principal symbol $p_h$ has a NHT invariant submanifold $\Gamma$. Then for

\begin{equation}
\begin{split}
    &\left(2 \nu_u - 2p_{h, 1}  - (2s - m + 3)g \right)|_{\Gamma} > 0 \\
    &\left(\nu_s - 2p_{h, 1}  - (2s - m + 1)g \right)|_{\Gamma} > 0 
\end{split}
\end{equation}
the following estimate near $\Gamma$ holds for all $v \in H^{-N, l}$ with $WF^{t, l}(v) \cap \Gamma = \emptyset$, $t > s_0$ and $s_0$ being the infimum over all $s \in \mathbb{R}$ satisfying the conditions above:

 \begin{equation}
 \label{SemiclassicalNHT_Estimate}
\begin{split}
     &\| \Op(\chi) v \|_{H^s_h} + \sum_{j=1}^{codim(\Gamma_u)} h^{-1}\| \Op(\chi) \Phi_{u, j} v \|_{H^{s+1}_h} \\ 
      & \qquad \lesssim h^{-1} \| G P v\|_{H^{s-m+1}_h} + \sum_{j=1}^{codim(\Gamma_u)} h^{-2}\| G \Phi_{u, j} P v \|_{H^{s-m+2}_h} \\
      & \qquad \hspace{4em} + h^{-1} \| \Op(\chi_{\phi_u \neq 0}) v \|_{H^{s+1}_h} + h^{N}\| E v \|_{H^{-N}_h}
\end{split}
\end{equation}
where the semiclassical pseudodifferential operators $\Op(\chi), G, \chi_{\phi_u \neq 0}, E$ are as described in theorem \ref{TheoremClosed} and g, $p_{h, 1}$ is defined analogously as $g = \frac{1}{\rho^{m}}H_{p_h}(\rho)$ and $p_{h, 1} = \frac{1}{2ih} \frac{\sigma_{m-1}(P_h - P^{\dagger}_h)}{\rho^{m-1}}$ for some homogeneous of degree 1 function $\rho$ on $T^{\ast}M$.

 \end{theorem}
\vspace{0.2in}
 In particular we have a high-energy estimate for the normal operator $P|_{\tau = 0} = N(P)$, which will end up showing it is invertible for $Re(\sigma) >> 1$. In practice, we actually do not really need this optimal estimate as the invertibility of the normal operator is only necessary to improve the error term $\| v \|_{-N, l}$ in our estimates to $\| v \|_{-N, l'}$, $l' > l$. But to do so, one may take an arbitrary loss in $h$ and regularity $s$, hence Dyatlov's result \cite[Theorem 1]{Dyatlov_2016} 

 \begin{equation}
     \| v \|_{H^s_h} \lesssim h^{-2} \| Pv \|_{H^{s-m+2}_h} + h^N \| v \|_{H^{-N}_h}
 \end{equation}
 is sufficient for that purpose. Nonetheless, it is worth mentioning.

\section{Proof of Theorem \ref{TheoremB}}
\label{ProofNHTEstimate}

\begin{figure}[ht!]
    \centering
    \includegraphics[width=0.7\linewidth]{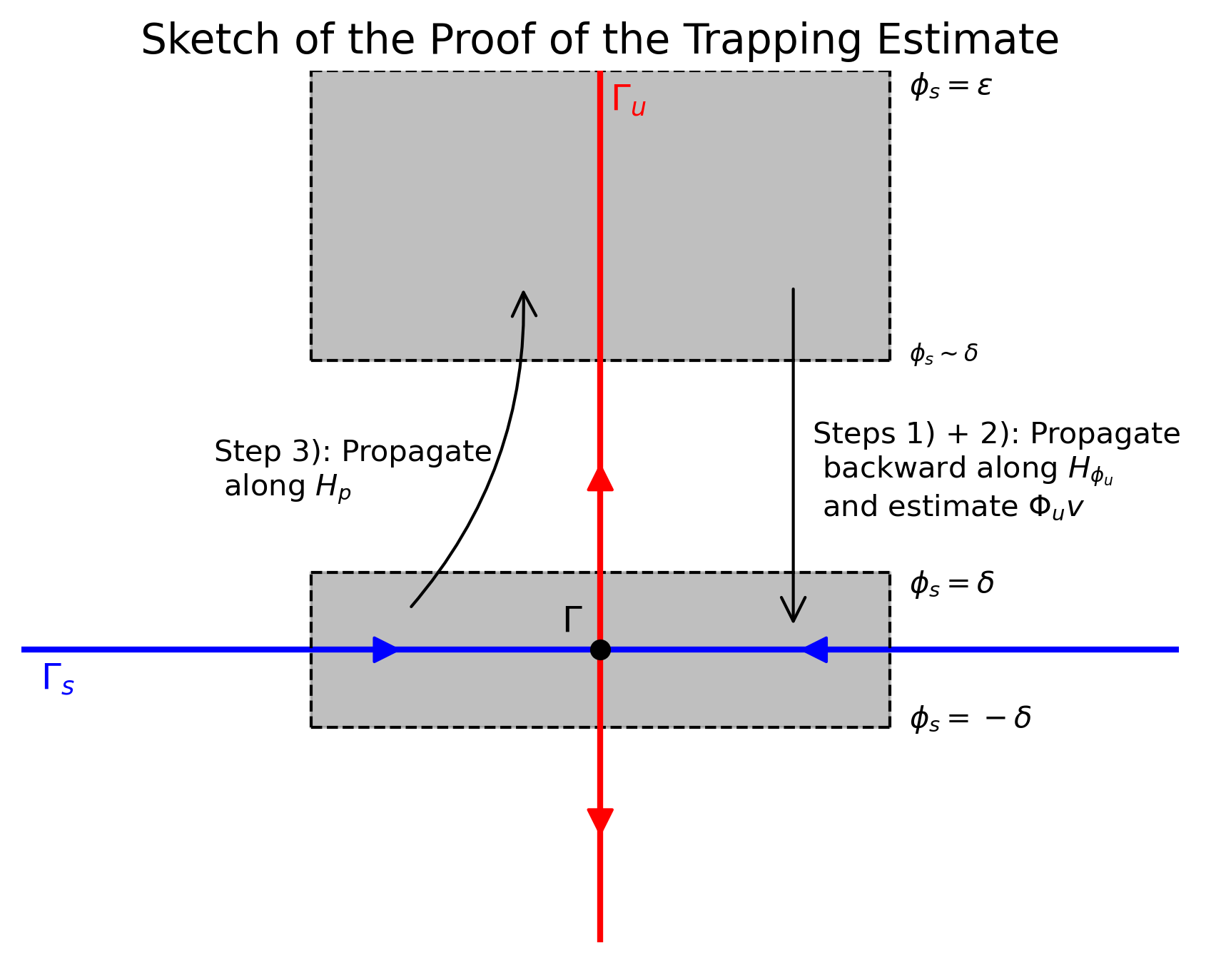}
    \caption{A sketch of the proof, illustrating the propagation along $H_{\phi_u}$ and $H_p$ we are doing. The arrows are not shown to scale; one should think of $\phi_s$ as growing exponentially fast along $H_p$ and linearly along $H_{\phi_u}$.}
    \label{fig:ProofNHT}
\end{figure}

We first give a short description of how Dyatlov's and Hintz's proof goes. To estimate $v$ on a neighborhood of the trapped set, we will use propagation estimates to control it by $v$ somewhere else. We can't move off $\Gamma$ using $\Tilde{H}_p$, so the only remaining thing to try is to propagate parallel to the unstable manifold along $\Tilde{H}_{\phi_u}$. Since our intersection is transversal, that has the effect of leaving the stable manifold in linear time $\sim \delta$. Thus, it now remains to control what we pick up by transporting which is $\sim \Phi_u v$ and $v$ off the stable manifold in $\phi_s \gtrsim \delta$. The first can be done using that $\Phi_u v$ is more regular than $v$ due to the damping $W_u$ provides for it. As for estimating $v$ in $\phi_s \gtrsim \delta$, we can control it by propagating along $\Tilde{H}_p$. Since $\phi_s$ grows exponentially fast near the trapped set, it only takes a logarithmic time $\sim \log(\delta^{-1})$ for it to reach the region. Thus, we can control $v$ in $\phi_s \gtrsim \delta$ by a relatively small amount (in the sense of the dependence on $\delta$) of $v$ in our original neighborhood. By comparing the two scaling constants in these propagation estimates, one then obtains a negligible overall factor. This then yields an estimate assuming only control off the unstable manifold, which is what we want. These steps can be graphically represented in figure \ref{fig:ProofNHT}.
Where our proof differs from theirs is really simple: we preserve the $\| \Phi_u v \|_{s+1, l}$ and $\| \Phi_u Pv \|_{s-m+2, l}$ terms rather than bounding them by the corresponding norms $\| v \|_{s+1, l}$ and $\| Pv \|_{s-m+2, l}$. The rest of the proof remains the same. This gives the local Fredholm estimate on coisotropic space, which we will then globalize in section \ref{sec5}. This will prove to be a bit trickier and is addressed in appendix \ref{AppendixA} and \ref{AppendixB}. Furthermore, as we are working in the b-calculus, we must take into account the constraints that $l$ puts on the estimate by $\Tilde{H}_p(\tau)/\tau$. Finally, one last distinction is that we can estimate $\| \tau^{\mu} v \|_{s+1, l}$ using the weight $\Tilde{H}_p(\tau)/\tau = -f$. This would fail in the cusp setting because $\Tilde{H}_p(\tau)/\tau$ would vanish at the boundary hence offer no control for $\| \tau^{\mu} v \|_{s+1, l}$, which is important for our Fredholm theory to be independent of the extensions of $\phi_u$ outside the boundary. It should also be noted that all of the remaining results assume that $v$ is sufficiently regular so that the pairings and integrations by parts we do are defined. We will discuss in the end how to remove that assumption through a regularization procedure. Furthermore, as the proof is very similar to Hintz's proof in the cusp setting [\cite{Hintz_2021} pages 38-44], we will be rather concise with the explanations and refer to his paper for a more detailed exposition.

\subsection{Control away from the Unstable Manifold}

We start the proof by estimating $v$ away from the unstable manifold. To do so, we make use of the fact that $\Phi_u v$ satisfies a better equation \ref{DampingEquation}, as the weight $W_u$ provides damping away from the unstable manifold.

\begin{proposition}
Let $\epsilon > 0$ be small enough and assume the following 

\begin{equation}
    \left( 2 \min(\nu_u, f\mu) - 2p_1 - 2lf - (2s - m + 3)g \right)|_{\Gamma} > 0
\end{equation}
Then there is a a cutoff $\Op(\chi_{\epsilon}) \in \Psi^0_b(M)$ microlocally supported in a $\epsilon$-neighborhood of $\Gamma$, i.e. $|\tau|, |p|, \|\phi_s\|, \|\phi_u\| < \epsilon$, elliptic at $\Gamma$ s.t. the following inequality holds:

\begin{equation}
\label{ControlAwayUnstableEstimate}
\begin{split}
      & \sum_{j=1}^{codim(\Gamma_u)} \| \Op(\chi_{\epsilon}) \Phi_{u, j} v \|_{s+1, l} + \| \Op(\chi_{\epsilon}) \tau^{\mu} v \|_{s+1, l} \\
      & \qquad \lesssim \sum_{j=1}^{codim(\Gamma_u)} \| G \Phi_{u, j} P v\|_{s-m+2, l} + \| G \tau^{\mu} P v\|_{s-m+2, l} \\ 
      & \qquad \hspace{4em} + \| G P v\|_{s-m+1, l} + \| \Op(\chi_{\phi_u \neq 0}) v \|_{s+1, l} + \lVert \Op(\chi_{\tau \neq 0}) v \rVert_{s+1, l} + \| E v \|_{s, l}
\end{split}
\end{equation}
For some $\Op(\chi_{\phi_u \neq 0}) \in \Psi^0_b(M)$ whose wavefront set is disjoint from the unstable manifold, $\Op(\chi_{\tau \neq 0}) \in \Psi^0_b(M)$ whose wavefront set is disjoint from the boundary and some $G, E \in \Psi^0_b(M)$ elliptic on $WF(\Op(\chi_{\epsilon})$. All their wavefront sets can be taken to be contained in the $\epsilon$-neighborhood of $\Gamma$.
\label{proposition}
\end{proposition}
\textbf{Proof:}
Pick cutoff functions $\chi_{\tau}$, $\chi_s$, $\chi_u$ and $\chi_p$ that are 1 in $(-\epsilon^2, \epsilon^2)$ and such that $\sqrt{- \chi_{\tau}' \chi_{\tau}}$ is smooth. We define the commutant

\begin{equation}
    a^2 = \tau^{-2l} \rho^{2(s+1)-m+1} \chi_s^2(\| \phi_s \|^2) \chi_u^2(\| \phi_u \|^2) \chi_p^2(p/\rho^{m}) \chi^2_{\tau}(\tau^2)
\end{equation}
 Let $A$ be a quantization of $a$. Furthermore, we take $A = A_{\rho} \Op(\chi)$ where $A_{\rho}$ quantize the weights $\rho^{(2(s+1)-m+1)/2} \tau^{-l}$ and $\Op(\chi)$ the cutoffs, both self-adjoint. Then we have for $P_u = P - iW_u$, $|A|^2 = A^{\dagger} A$ and the tuple $v_u = \Phi_u v$

\begin{equation}
    2 Im \langle P_u v_u, |A|^2 v_u \rangle = \langle i[P_u, |A|^2]v_u, v_u \rangle + \langle \frac{1}{i}(P_u - P_u^\dag) |A|^2 v_u, v_u \rangle
\end{equation}
with the pairing being

\begin{equation*}
    \langle u, v \rangle = \int \sum_j u_j v_j dy d\rho / \rho
\end{equation*}
The principal symbol of the commutator is then

\begin{equation}
\begin{split}
    H_p(a^2) &= b_u - b_s + b_{\tau}^2 - b_e + hp \\
    b_s &= -4 \chi_s' \chi_s \langle w_s \phi_s, \phi_s \rangle  \tau^{-2l} \rho^{2(s+1)} \chi^2_u \chi^2_p \chi_{\tau}^2 \\
    b_u &= -4 \chi_u' \chi_u \langle w_u \phi_u, \phi_u \rangle  \tau^{-2l} \rho^{2(s+1)} \chi^2_s  \chi^2_p \chi_{\tau}^2 \\
    b_{\tau} &= \sqrt{-2 \chi_{\tau}' \chi_{\tau} f} \tau^{-l + 1/2} \rho^{2(s+1)} \chi^2_s  \chi^2_p \chi^2_{\tau} \\
    b_{e} &= (- 2lf - (2(s+1) - m + 1) g) \rho^{m-1} a^2 \\
    h &= (4 r_s \chi_s' \chi_s \phi_s \chi_u^2 \chi^2_p - 4 r_u \chi_u' \chi_u \phi_u \chi_s^2 \chi^2_p - 2 \chi_p' \chi_p m \frac{1}{\rho^{m}} H_p(\rho) \chi^2_u \chi^2_s) \\
    & \times \tau^{-2l} \rho^{2(s+1)-m} \chi^2_{\tau}
\end{split}
\end{equation}
Note that besides the term factoring through $p,$ this is a sum of terms with definite sign. One should think of minus as the 'good' sign, meaning  it will provide the estimate or can be ignored, and plus as the 'bad' sign, meaning we must assume a priori control of these terms. The subprincipal term has principal symbol

\begin{equation}
    \sigma_{2(s+1)}(\frac{1}{i}(P_u - P_u^\dag) |A|^2) = (2p_1 - w_u - w_u^{T}) \rho^{m-1} a^2
 \end{equation}
Note that $w_u + w_u^T - 2 p_1 - 2lf - (2s-m+3)g > 0$ is a positive definite matrix at the trapped set by assumption. Denote by $W$ a quantization of this term; it will provide the estimate for $v_u$. Quantizing everything and rearranging gives the following

\begin{align*}
    &\langle W |A|^2 v_u, v_u \rangle + \langle B_s v_u, v_u \rangle = \\
    &\langle B_u v_u, v_u \rangle  + \langle B_{\tau} v_u, B_{\tau} v_u \rangle + 2 Im \langle |A|^2 v_u, P_u v_u \rangle + \langle HP v_u, v_u \rangle + \langle R_3 v_u, v_u\rangle
\end{align*}
With $R_3$ an error term of order $2s + 1$ microlocally supported near the trapped set. For any $\delta > 0$, the third term $Im(...)$ on the right is bounded by

\begin{align}
    2 Im \langle |A|^2 v_u, P_u v_u \rangle & \leq 
     2 \lVert |A_{\rho}|^2 \chi v_u \rVert_{-s + m - 2, -l} \lVert \Op(\chi) P_u v_u \rVert _{s - m + 2, l} \\
     & \lesssim \delta \lVert \Op(\chi) v_u \rVert_{s+1, l}^2 + \frac{1}{\delta} \lVert \Op(\chi) P_u v_u \rVert _{s - m + 2, l}
\end{align}
Using the expression for $P_u v_u$ above we can then bound $\| \Op(\chi) P_u v_u \rVert _{s - m + 2, l}$ as well to obtain

\begin{align*}
    &\langle W |A_\rho| \Op(\chi) v_u, |A_{\rho}| \Op(\chi) v_u \rangle - \delta \lVert \Op(\chi) v_u \rVert_{s+1, l}^2 + \langle B_s v_u, v_u \rangle \leq  \\
    & \qquad \frac{1}{\delta} (\lVert  \Op(\chi) \Phi_u P v \rVert _{s - m + 2, l}^2 + \lVert \Op(\chi) R_1 P v \rVert _{s - m + 2, l}^2 + \lVert \Op(\chi) R_2 v \rVert _{s - m + 2, l}^2) & \\
    & \qquad \hspace{4em} + \langle B_u v_u, v_u \rangle + \langle B_{\tau} v_u, B_{\tau} v_u \rangle + \langle HP v_u, v_u \rangle + \langle R_3 v_u, v_u\rangle &
\end{align*}
We can then use G\r{a}rding inequality to obtain for $N$ arbitrary

\begin{equation}
     \| \Op(\chi) v_u \|_{s+1, l}^2 \lesssim \langle W \Op(\chi) |A_\rho| \Op(\chi) v_u, |A_{\rho}| \Op(\chi) v_u \rangle + \lVert E v_u \rVert^2_{-N, l}
\end{equation}
Where $E$ is some cutoff elliptic on $WF(\chi)$. Then, taking $\delta$ small enough to absorb the $\delta \lVert \Op(\chi) v_u \rVert_{s+1, l}^2$ term, discarding the $B_s$ contribution (up to $\| E v_u \|_{-N, l}$ error using G\r{a}rding's inequality) and bounding the pairings on the right by the corresponding norms we get 

\begin{equation}
\begin{split}
    & \sum_j \lVert \Op(\chi) \Phi_{u, j} v \rVert_{s+1, l}^2 \lesssim \sum_j \lVert \Op(\chi) \Phi_{u, j} P v \rVert _{s - m + 2, l}^2 + \lVert \Op(\chi) P v \rVert _{s - m+1, l}^2 + \lVert E v \rVert _{s, l}^2 \\
    & \qquad + \lVert \Op(\chi_{\phi_u \neq 0}) v \rVert_{s+1, l}^2 + \lVert \Op(\chi_{\tau \neq 0}) v \rVert_{s+1, l}^2 + \lVert E \Phi_u v \rVert _{s + 1/2, l}^2 + \lVert E v \rVert_{-N, l}
\end{split}
\end{equation}
 In order to get the similar estimate for $\tau^{\mu} v$, we apply the exact same method to

\begin{equation}
    (P+\frac{\mu}{i} F)\tau^{\mu} v = \tau^{\mu} P v + R_3 P v + R_4v
\end{equation}
where $F$ quantizes $f = \frac{1}{\tau} \Tilde{H}_p(\tau)$. We now use 

\begin{equation}
    \left( 2 \mu f - 2lf - 2 p_1 - (2s-m+3)g \right)|_{\Gamma} > 0
\end{equation}
in order to obtain the estimate. Combining the two inequality, taking the square root and repeatedly applying it to improve the error of order $s + 1/2$ gives the result.
$\Box$

\subsection{Propagation Parallel to the Unstable Manifold}\label{subsec2}

We now describe the \textit{quantitative} propagation estimate along $H_{\phi_u}$ we will use. It consists of estimating the function in a $\sim \delta$ sized region around the trapped set $\Gamma$ by its norm in a bigger region times a small constant. As the 'distance' travelled from the trapped set (measured using $\phi_s$) grows linearly in time when propagating along $H_{\phi_u}$, the scaling constant in the inequality will be $\sqrt{\delta} \sim \frac{\| \chi_{\delta} v \|_s}{\| \chi_{\| \phi_s \| > \delta} v \|_s}$. This is encapsulated in the following proposition.

\begin{figure}[ht!]
    \centering
    \includegraphics[width=0.7\linewidth]{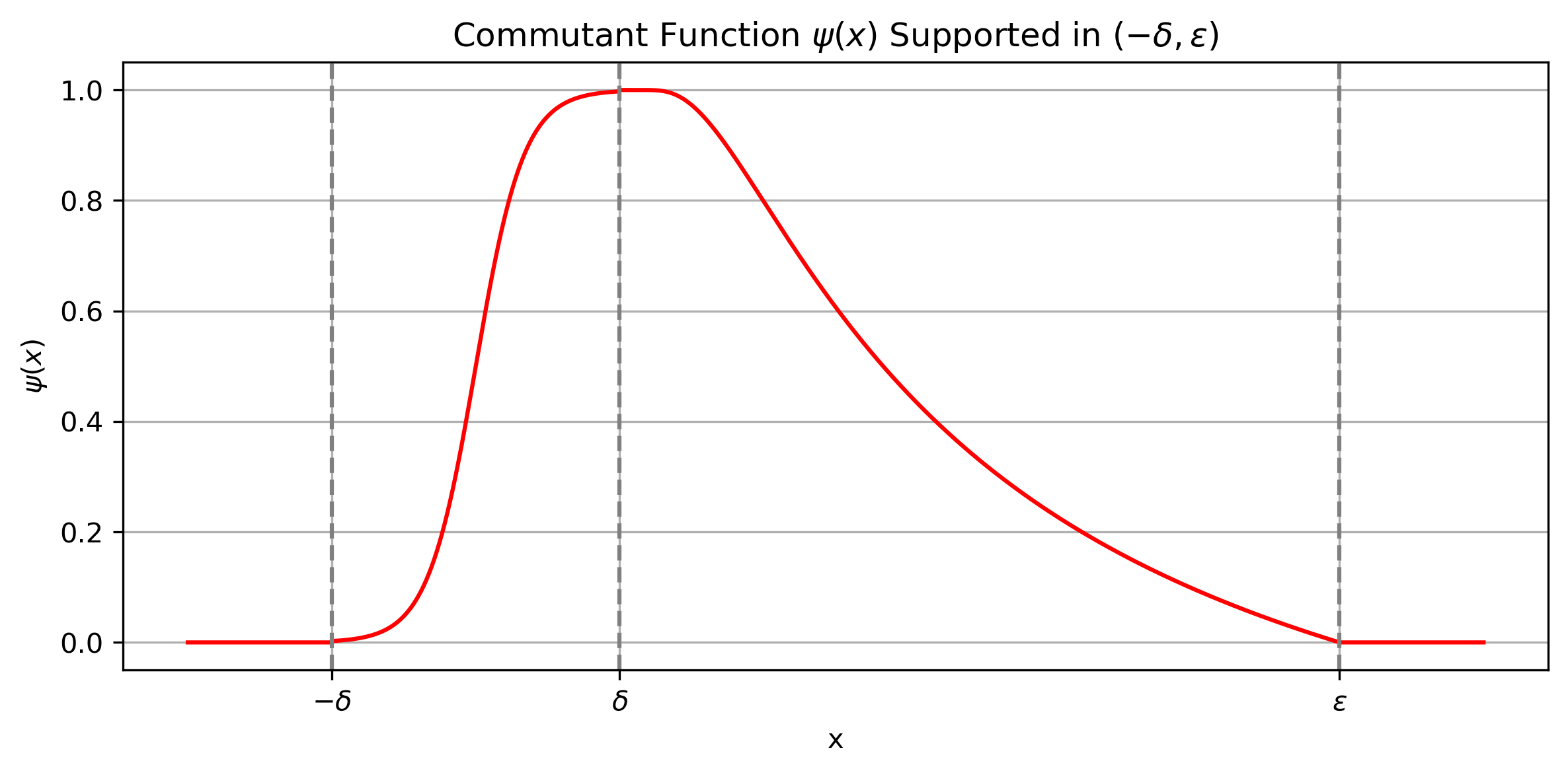}
    \caption{An illustration of the cutoff $\psi$ we use for the quantitative propagation estimate along $\sim H_{\phi_u}$. The rapid increase in $(-\delta, \delta)$ provides the large positive commutator we need for the estimate, while the controlled decay in $(\delta, \epsilon)$ corresponds to the assumed control in $\delta \lesssim \phi_s \lesssim \epsilon$.}
    \label{fig:CommutantDelta}
\end{figure}

\begin{proposition}
Let $\sum_j f_j \rho H_{\Phi_{u, j}}(\phi_s^1)|_{\Gamma} = 1$ be the symplectic dual of $H_{\phi_s^1}$ on $\Gamma$. Then for a small enough $\epsilon$, the following holds for any propagation time $0 < \delta < \epsilon$ along the Hamiltonian flow of $\sum_i \rho f_i \phi_u^i$ 

\begin{equation}
\label{PropagationAlongHphiEstimate}
\begin{split}
    &\| \Op(\chi_{\delta}) v \|_{s, l} \lesssim \sqrt{\delta} ( \| \Op(\chi_{\phi_u \neq 0}) v \|_{s, l} + \| \Op(\chi_{\tau \neq 0}) v \|_{s, l} + \| G P v \|_{s-m, l} + \| \Op(\chi_{\| \phi_s \| > \delta}) v \|_{s, l} \\
     & \qquad + \sum_{j=1}^{codim(\Gamma_u)} \| \Op(\chi_{\epsilon}) \Phi_{u, j} v \|_{s+1, l} ) + C_E \| E v \|_{s-1/2, l} 
\end{split}
\end{equation}
For some constant $C_E$ depending on both $\delta$ and $\epsilon$ and some microlocal cutoffs supported in the following neighborhoods of $\Gamma$ (when not specified, $\tau \lesssim \epsilon$)

\begin{equation}
    \begin{split}
        &WF(\Op(\chi_{\delta})), \chi_{\delta} \geq 1 \quad in \quad \{ \| \phi_u \| \lesssim \epsilon, |p| \lesssim \epsilon, \| \phi_s \| \lesssim \delta \} \\
        &WF(\Op(\chi_{\phi_u \neq 0})) \quad \subset \quad  \{ \frac{\epsilon}{2} \lesssim \| \phi_u \| \lesssim \epsilon, |p| \lesssim \epsilon, \| \phi_s \| \lesssim \epsilon \} \\
         &WF(\Op(\chi_{\tau \neq 0})) \quad \subset \quad  \{ \| \phi_u \| \lesssim \epsilon, |p| \lesssim \epsilon, \| \phi_s \| \lesssim \epsilon, \frac{\epsilon}{2} \lesssim \tau \lesssim \epsilon \} \\
        &WF(\Op(\chi_{\| \phi_s \| > \delta})), \chi_{\| \phi_s \| > \delta} \leq 1 \quad in \quad \{ \| \phi_u \| \lesssim \epsilon, |p| \lesssim \epsilon, \delta \lesssim \| \phi_s \| \lesssim \epsilon \} \\
        &WF(\Op(\chi_{\epsilon})), WF(G), WF(E) \quad in \quad \{ \| \phi_u \| \lesssim \epsilon, |p| \lesssim \epsilon, \| \phi_s \| \lesssim \epsilon \}
    \end{split}
\end{equation}

\label{proposition2}
\end{proposition}
\textbf{Proof:} We define the principal symbol of the commutant exactly as in Hintz's paper 

\begin{equation}
    \begin{split}
    a_{\pm}^2 &= \tau^{-2l} \rho^{2s+1} \psi_{\pm}^2(\phi_s^1) \chi^2(\sum_{j \neq 1} (\phi_s^j)^2 + \| \phi_u \|^2) \chi_p^2(p/\rho^{m}) \chi^2_{\tau}(\tau) \\
    a^2 &= a_{-}^2 + a_{+}^2
    \end{split}
\end{equation}
The $\chi$ cutoffs are as before i.e. 1 in $\sim (-\epsilon^2, \epsilon^2)$ but this time, we have to be quantitative about the behavior of the cutoffs $\psi_{\pm}$ in $\phi_s^1$ since these terms will be the one providing the estimate. We arrange so that for $\psi^2 = \psi_{-}^2+\psi_{+}^2$ we have $2 \psi \psi' \simeq \frac{1}{\delta} \psi_{-}^2 + \frac{1}{\delta} \Tilde{b}_{\delta}^2 - \Tilde{e}$ with $\Tilde{b}_{\delta} \gtrsim 1$ in $\sim [-\delta, \delta]$ being positive and giving the control we want while $|\Tilde{e}| \lesssim \frac{1}{\epsilon}$ is the error where we require a priori regularity. The functions $\psi_{-}$ and $\psi_{+}$ are supported respectively in $supp (\psi_{-}) \sim [-\delta, \delta]$ and $supp (\psi_{+}) \sim [\delta, \epsilon]$, while the error $\Tilde{e}$ is supported in $\sim[\delta, \epsilon]$ too. See fig. \ref{fig:CommutantDelta} for an illustraton of $\psi$.
The other commutant we will use is 

\begin{equation*}
    \sum_j F_j \Phi_{u, j}
\end{equation*}
with $F_j$ being self adjoint quantization of $f_j$. Its Hamiltonian vector field is then

\begin{equation}
    H_{\sum_j f_j \phi_{u, j}} = \sum_j f_j H_{\phi_{u, j}} + \sum_j \phi_{u, j} H_{f_j}
\end{equation}
hence is exactly the symplectic dual of $H_{\phi_s^1}$ at $\Gamma$. We then have

\begin{equation}
    2 Im \langle \sum_j F_j \Phi_{u, j} v, |A|^2 v \rangle = \langle i [\sum_j F_j \Phi_{u, j}, |A|^2] v, v \rangle + \langle \frac{1}{i}(\sum_j [F_j, \Phi_{u,j}]) |A|^2 v_u, v_u \rangle
\end{equation}
If we let $w$ be defined by $\rho H_{\sum_j f_j \phi_{u, j}}(\phi_s^1) = w$, with $w|_{\Gamma} = 1$, the principal symbol of the commutator is then given by

\begin{equation}
\begin{split}
    &H_{\sum_j f_j \phi_{u, j}} (a^2) \simeq \frac{1}{2 \rho \delta} a_{-}^2 + b_{-}^2 + b_{\delta}^2 - e - b_p - b_{\tau} - b_{\chi'} + \sum_j \phi_{u, j} H_{f_j}(a^2) \\
    &b_{-} =  \sqrt{ \frac{1}{\delta} (w - 1/2) - \frac{2l}{\tau} \sum_j \rho f_j H_{\phi_{u, j}}(\tau) + (2s+1) \sum_j f_j H_{\phi_{u, j}}(\rho) - \sum_j H_{f_j}(\phi_{u, j})} \\
    &\hspace{4em} \times \tau^{-2l} \rho^{s} \psi_{-} \chi_u \chi_p \chi_{\tau} \\
    &b_{\delta} =  \sqrt{ \frac{w}{\delta}} \Tilde{b}_{\delta} \tau^{-2l} \rho^{s} \chi_u \chi_p \chi_{\tau} \\
    &b_{\chi'} = - 2 \chi \chi' \psi^2 \chi_p^2 \chi_{\tau}^2 \tau^{-2l} \rho^{2s+1} \sum_j f_j H_{\phi_{u, j}} (\sum_{i \neq 1} (\phi_s^i)^2) \\
    &b_p =  - 2 \chi^2 \psi^2 \chi_p \chi_p' \chi^2_{\tau} \tau^{-2l} \rho^{2s+1} \sum_j f_j H_{\phi_{u, j}} (p/\rho^{m}) \\
    &b_{\tau} =  - 2 \chi^2 \psi^2 \chi^2_{p} \chi_{\tau} \chi_{\tau}' \tau^{-2l} \rho^{2s+1} \sum_j f_j H_{\phi_{u, j}} (\tau) \\
    &e = (w \Tilde{e} + \psi_+^2 ((2s+1) \sum_j f_j H_{\phi_{u, j}}(\rho) - \frac{2l}{\tau} \sum_j \rho f_j H_{\phi_{u, j}}(\tau) - \sum_j H_{f_j}(\phi_{u, j}))) \chi^2 \chi_p^2 \chi_{\tau}^2 \tau^{-2l} \rho^{2s} 
\end{split}
\end{equation}
The $a_{-}$ and $b_{\delta}$ terms are what we will use to estimate $v$, $e$ is where we require the appriori control and $b_p$, $b_{\tau}$ are errors supported away from the characteristic set and the boundary respectively. Note that the terms in the square root are positive for $\delta$ small enough. As for $b_{\chi'}$, it is supported in an annulus

\begin{equation}
    \left\{ \epsilon^2 \lesssim \sum_{j \neq 1} (\phi_{s,j})^2 + \| \phi_u \|^2 \lesssim 2 \epsilon^2 \right\}
\end{equation}
which is contained in the regions
\begin{align}
    &\left\{ \sum_{j \neq 1} (\phi_{s, j})^2 \lesssim \delta^2 \quad  \epsilon^2 - \delta^2 \lesssim \| \phi_u^2 \| \lesssim 2 \epsilon^2 \right\} \\
    \bigcup &\left\{ \delta^2 \lesssim \sum_{j \neq 1} (\phi_{s, j})^2 \lesssim 2 \epsilon^2 \quad \| \phi_u^2 \| \lesssim 2 \epsilon^2 - \delta^2 \right\}
\end{align}
The first of these is located away from the unstable manifold $\Gamma_u$ hence is in the control region while on the second region the function needs to be estimated. Quantizing all these symbols and rearranging gives the following equation:

\begin{equation}
    \begin{split}
    & \langle \frac{1}{\rho \delta} |A_{-}|^2 v, v \rangle + \langle B_{\delta} v, B_{\delta} v \rangle + \langle B_{-} v, B_{-} v \rangle \simeq \\
    & \qquad \langle B_p v, v \rangle + \langle B_{\tau} v, v \rangle + \langle B_{\chi'} v, v \rangle + \langle \Op(\chi_e) v, v \rangle - \sum_j \langle \Phi_{u, j} [F_j, |A|^2] v, v \rangle \\
    & \qquad \hspace{4em}+ 2 Im \langle\sum_j F_j \Phi_{u, j} v, |A|^2 v \rangle + \langle R v, v \rangle
    \end{split}
\end{equation}
with $R$ of order $2s-1$. As in the previous section, we estimate the imaginary term using Cauchy-Schwarz and estimate the resulting product as sums of squares with a big and a small constant in front. Furthermore, we bound the reminder term $R$ and the a priori control terms $E, B_p, B_{\tau}, B_{\chi'}$ by microlocalized $s-1/2$ and $s$ norms respectively. After doing these estimates, we obtain the following.

\begin{align*}
    &\frac{1}{\rho \delta} \langle |A_{-}|^2 v, v \rangle + \langle B_{\delta} v, B_{\delta} v \rangle + \langle B_{-} v, B_{-} v \rangle \\
    &\qquad \lesssim \langle B_p v, v \rangle + \langle B_{\tau} v, v \rangle + \langle B_{\chi'}v, v \rangle + \langle \chi_e v, v \rangle + \sum_j |\langle \Phi_{u, j} [F_j, |A|^2] v, v \rangle| \\ 
    &\qquad \hspace{2em} +( \frac{1}{\delta} \| A_{-} v \|^2_{-1/2} + \delta \| A_{-} \sum_j F_j \Phi_{u, j} v \|^2_{1/2} + \frac{1}{\epsilon} \| A_{+} v \|^2_{-1/2} + \epsilon \| A_{+}  \sum_j F_j \Phi_{u, j} v \|^2_{1/2} ) + \langle R v, v \rangle \\
    & \qquad \lesssim C_p \| \Op(\chi_{p \neq 0}) \|^2_{s, l} + C_{\tau} \| \Op(\chi_{\tau \neq 0}) v \|_{s, l}^2 + C_{\chi'} (\| \Op(\chi_{\phi_u \neq 0}) v\|^2_{s, l} + \| \Op(\chi_{\|\phi_s\| > \delta}) v \|^2_{s, l}) \\
    & \qquad \hspace{2em} + \frac{C_{e}}{ \epsilon } \| \Op(\chi_{\| \phi_s \| > \delta}) v \|_{s, l} + \sum_j |\langle \Phi_{u, j} [F_j, |A|^2] v, v \rangle| \\ 
    & \qquad \hspace{4em} + ( \frac{1}{\delta} \| A_{-} v \|^2_{-1/2} + \delta \| A_{-} \sum_j F_j \Phi_{u, j} v \|^2_{1/2} + \epsilon \| A_{+} \sum_j F_j \Phi_{u, j} v \|^2 _{1/2} ) + C_R \| \chi_{R} v \|^2_{s-1/2, l}
 \end{align*}
with $\chi_{\phi_u \neq 0}$ elliptic in the region $\| \phi_u \|^2 \gtrsim \epsilon^2 - \delta^2$, $\chi_{\| \phi_s \| > \delta}$ elliptic in the region $\| \phi_s \| \gtrsim \delta$ and $\chi_R$ elliptic in the whole $\epsilon$ neighborhood. Now, the first term in parentheses cancels the one on the left handside while the $B_{-}$ term on the left can be neglected. Furthermore, we multiply by $\delta$ to factor out the $\frac{1}{\delta}$ dependence of $b_{\delta}^2$. Thus, we obtain the following inequality (where the constant $C_E$ depends on $\delta$)

\begin{equation}
\begin{split}
    &\| \Op(\chi_{\delta}) v \|^2_{s, l} \lesssim \delta (\| \Op(\chi_{\phi_u \neq 0}) v \|^2_{s, l} + \| \Op(\chi_{\tau \neq 0}) v \|^2_{s, l} + \| \Op(\chi_{p \neq 0}) v \|^2_{s, l} \\
    & \qquad + \| \Op(\chi_{\| \phi_s \| > \delta}) v \|_{s, l}^2 + \sum_j \| \Op(\chi_{\epsilon}) \Phi_{u, j} v \|^2_{s+1, l}) + C_E \| Ev \|^2_{s-1/2, l}
\end{split}
\end{equation}
where the cutoffs are micro locally supported as described in proposition \ref{proposition2}. Taking the square-root of this inequality and applying an elliptic estimate for the $\chi_{p \neq 0} v$ term gives the desired result. 
$\Box$

\subsection{Propagation along the Hamiltonian Flow of \texorpdfstring{$H_p$}{Hp}}
\label{subsec3}

We now propagate along the $H_p$ flow for a logarithmic time in order to be able to estimate the term $\| \chi_{\| \phi_s \| > \delta} v \|_s$ by $v$ in the smaller region $\| \chi_{\delta} v \|_s$. Due to this logarithmic behavior, we only need a $\delta^{-\beta} \sim \frac{\| \chi_{\| \phi_s \| > \delta} v \|_s}{\| \chi_{\delta} v \|_s}$ scaling constant. Once we combine everything, taking $\beta$ small enough will result in a small overall factor $\delta^{1/2 - \beta}$ multiplying $\| \chi_{\delta} v \|_s$, which allows us to neglect that term. 

\begin{proposition}
For a propagation time $\sim \log(\delta^{-1})$ along the flow of $H_p$, any $\beta > 0$ and 

\begin{equation*}
    \frac{2 \beta \langle w_s \phi_s, \phi_s \rangle}{\| \phi_s \|^2} - 2p_1 - 2lf - (2s - m + 1)g > 0
\end{equation*}

the following holds

\begin{equation}
\label{PropagationAlongHpEstimate}
\begin{split}
    & \| \Op(\chi_{\| \phi_s \| > \delta}) v \|_{s, l} \lesssim \delta^{-\beta} (\| \Op(\chi_{\delta}) v \|_{s, l} + \| \Op(\chi_{\phi_u \neq 0}) v \|_{s, l} + \| \Op(\chi_{\tau \neq 0}) v \|_{s, l} + \| G Pv \|_{s-m+1, l}) \\
    & \qquad + C_E \| E v \|_{s-1/2, l}
\end{split}
\end{equation}
For some constant $C_E$ depending on $\delta$ and some microlocal cutoffs supported in the following neighborhoods of $\Gamma$ (when not specified, $\tau \lesssim \epsilon$)

\begin{equation}
    \begin{split}
    & WF(\Op(\chi_{\| \phi_s \| > \delta})), \chi_{\| \phi_s \| > \delta} \geq 1 \quad in \quad \{ \| \phi_u \| \lesssim \epsilon, |p| \lesssim \epsilon, \delta \lesssim \| \phi_s \| \lesssim \epsilon  \} \\
    & WF(\Op(\chi_{\phi_u \neq 0})) \quad \subset \quad  \{ \| \phi_u \| \gtrsim \epsilon, |p| \lesssim \epsilon, \| \phi_s \| \lesssim \epsilon \} \\
    & WF(\Op(\chi_{\tau \neq 0})) \quad \subset \quad  \{ \| \phi_u \| \lesssim \epsilon, |p| \lesssim \epsilon, \| \phi_s \| \lesssim \epsilon, \tau \gtrsim \epsilon \} \\
    & WF(\Op(\chi_{\delta})) \quad \subset \quad \{ \| \phi_u \| \lesssim \epsilon, |p| \lesssim \epsilon, \| \phi_s \| \lesssim \delta  \} \\
    & WF(G), WF(E) \quad \subset \quad \{ \| \phi_u \|, |p|, \| \phi_s \| \lesssim \epsilon \}
    \end{split}
\end{equation}

\label{proposition3}
\end{proposition}

\begin{figure}[ht!]
    \centering
    \includegraphics[width=0.7\linewidth]{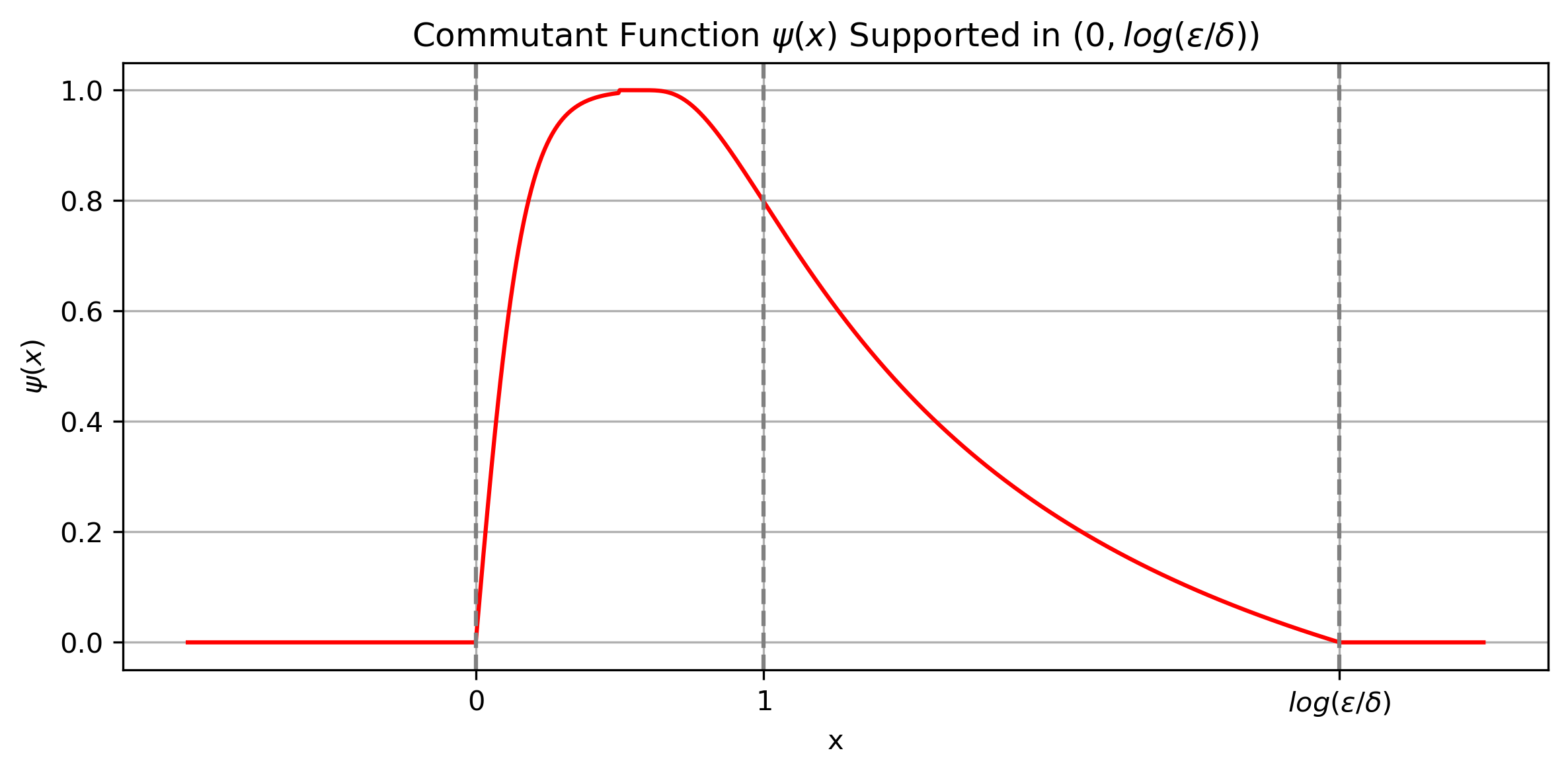}
    \caption{An illustration of the cutoff $\psi$ we use for the quantitative propagation estimate along $H_{p}$. The initial increase in $(0,1)$ corresponds to the region of a priori control $\| \phi_s \| \sim \delta$ while the slow decrease in $(1, \log(\epsilon / \delta))$ is how we will obtain our positive commutator estimate in $\delta \lesssim \| \phi_s \| \lesssim \epsilon$.}
    \label{fig:CommutantLog}
\end{figure}

\textbf{Proof:}
Once again, we follow Hintz's treatment. The commutant we will use is

\begin{equation}
    \begin{split}
    a_{\pm}^2 &= \tau^{-2l} \rho^{2s- m + 1} \| \phi_s \|^{-2\beta} \psi_{\pm}^2 (\log (\| \phi_s \|^2 / \delta^2)) \chi_u^2(\| \phi_u\|^2) \chi_p^2(p/\rho^{m}) \chi^2_{\tau}(\tau) \\
    a^2 &= a_{-}^2 + a_{+}^2
    \end{split}
\end{equation}
and we must compute the principal symbol of

\begin{equation}
    2 Im \langle P v, |A|^2 v \rangle = \langle i[P, |A|^2]v, v \rangle + \langle \frac{1}{i}(P - P^\dag) |A|^2 v, v \rangle
\end{equation}
Once again, we arrange so that $\psi^2 = \psi_{+}^2 + \psi_{-}^2$ satisfies some quantitative properties. We take it so that $\psi \psi' = -\psi^2 \Tilde{b}^2 + \Tilde{e}$, with $\Tilde{b} \geq 0$ and $\Tilde{e} \leq 1$, with supports $supp(\psi_{+}) \sim [1, \log(\epsilon/\delta)]$, $supp(\psi_{-}) \sim [0,1]$, $supp(\Tilde{b}) \sim [1, \log(\epsilon/\delta)]$ and  $supp(\Tilde{e}) \sim [0, 1]$ as illustrated in figure \ref{fig:CommutantLog}.
Given this, the principal symbol is (where $c$ is some small constant s.t. $0 < c < \frac{2 \beta \langle w_s \phi_s, \phi_s \rangle}{\| \phi_s \|^2} - 2p_1 - 2lf - (2s - m + 1)g $):

\begin{equation}
\begin{split}
    &H_p(a^2) + 2 p_1 a^2 = -c \rho^{m-1} a^2 - b_{+}^2 - b_{-}^2 + e_1 + e_2 + e_3 + hp \\
    & b_{\pm} =   \sqrt{ \frac{2(\beta + \Tilde{b}_{\pm}^2) \langle w_s \phi_s, \phi_s \rangle}{\| \phi_s \|^2} - 2p_1 - 2lf - (2s - m + 1)g - c } \\
    & \qquad \times \tau^{-2l} \rho^{s} \| \phi_s \|^{-\beta} \psi_{\pm} \chi_u \chi_p \chi_{\tau} \\
    &e_1 = 4 \chi_u \chi_u' \langle H_p(\phi_u), \phi_u \rangle \tau^{-2l} \rho^{2s-(m-1)} \| \phi_s \|^{-2\beta} \psi^2 \chi_p^2 \chi_{\tau}^2 \\
    &e_2 = 2 \Tilde{e} \frac{\langle w_s \phi_s, \phi_s \rangle}{\| \phi_s \|^2} \tau^{-2l} \rho^{2s} \| \phi_s \|^{-2\beta} \chi_u^2 \chi_p^2 \chi_{\tau}^2 \\
    &e_3 = 4 \chi_{\tau} \chi_{\tau}' H_p(\tau) \tau^{-2l} \rho^{2s-(m-1)} \| \phi_s \|^{-2\beta} \psi^2 \chi_u^2 \chi_p^2\\
    &h = (2 (\psi \psi' - \beta) \langle r_s, \phi_s \rangle \chi_p^2 - 2m \chi_p \chi_p' \frac{1}{\rho^m} H_p(\rho) \psi^2)\tau^{-2l} \rho^{2s-m} \| \phi_s \|^{-2\beta} \chi_u^2 \chi_{\tau}^2
\end{split}
\end{equation}
Note that the expression within the squared root is positive by our assumptions. Thus, quantizing everything, rearranging and using the usual bounds we obtain

\begin{equation}
    \begin{split}
    & c \langle |A|^2 v, v \rangle_{(m-1)/2} + \langle B_{+} v, B_{+} v \rangle + \langle B_{-} v, B_{-} v \rangle = \\
    & \qquad \langle E_1 v, v \rangle + \langle E_2 v, v \rangle + \langle E_3, v \rangle + \langle H P v, v \rangle + 2 Im \langle |A|^2 v, P v \rangle + \langle R v, v \rangle \\
    & \qquad \hspace{2em} \leq  \langle E_1 v, v \rangle + \langle E_2 v, v \rangle + \langle E_3, v \rangle + c \| Av \|_{(m-1)/2}^2 + \frac{2}{c} \| AP v \|_{-(m-1)/2}^2 \\
     & \qquad \hspace{4em} + C_H \| G P v \|_{s-m}^2 + C_R \| \Op(\chi_R) v \|_{s-1/2}^2
    \end{split}
\end{equation}
where the $C's$ are constant controlling the norms of the operators $H$ and $R$ respectively and $\Op(\chi_R), G \in \Psi_b^0$ are elliptic near the support of $a$. We can cancel the $c \| Av \|_{(m-1)/2}$ norms on each side, neglect the $B_{-}$ term and use that $e_2$ is $O(\delta^{-2\beta})$ to obtain the following bound
\begin{equation}
\begin{split}
    &\| \Op(\chi_{\| \phi_s \| > \delta}) v \|_s^2 \lesssim \delta^{-2\beta} (\| \Op(\chi_{\delta}) v \|_s^2 + \| \Op(\chi_{\phi_u \neq 0}) v \|_s^2 + \| \Op(\chi_{\tau \neq 0}) v \|_s^2 + \| G Pv \|_{s-m+1}^2) \\
    & \qquad \hspace{6em} + C_E \| E v \|_{s-1/2}^2
\end{split}
\end{equation}
where the cutoffs are supported as in proposition \ref{proposition3}. Taking the square root gives the result.
$\Box$

\subsection{Combining Everything to Prove Theorem \ref{TheoremB}}\label{subsec4}

Bringing everything together, we combine the three inequalities above \eqref{ControlAwayUnstableEstimate} \eqref{PropagationAlongHphiEstimate} \eqref{PropagationAlongHpEstimate}l to estimate $v$ near the trapped set by $v$ away from it. As discussed above, we will first propagate along the flow lines of $\phi_u$ to control $v|_{\Gamma}$ by $v$ slightly off the stable manifold. Then using our estimate for $\Phi_u v$ and propagating by $H_p$, we can control the terms that involve $v$ on the unstable manifold. All and all, it leads to the following chain of inequalities
\begin{align*}
         &\| \Op(\chi_{\delta}) v \|_{s, l} + \sqrt{\delta} \sum_{j} \| \Op(\chi_{\delta}) \Phi_{u, j} v \|_{s+1, l} + \sqrt{\delta} \| \Op(\chi_{\delta}) \tau^{\mu} v \|_{s+1, l} \\ 
         & \qquad \lesssim \sqrt{\delta} (\| \Op(\chi_{\phi_u \neq 0}) v \|_{s, l} + \| \Op(\chi_{\tau \neq 0}) v \|_{s, l} + \| G P v \|_{s-m, l} + \| \Op(\chi_{\| \phi_s \| > \delta}) v \|_{s, l} \\
         & \qquad \hspace{3em} + \sum_j \| \Op(\chi_{\epsilon}) \Phi_{u, j} v \|_{s+1, l} + \| \Op(\chi_{\epsilon}) \tau^{\mu} v \|_{s+1, l}) + C_E \| E v \|_{s-1/2, l}   \\      
        & \qquad \lesssim \sqrt{\delta} (\| \Op(\chi_{\phi_u \neq 0}) v \|_{s+1, l} + \| \Op(\chi_{\tau \neq 0}) v \|_{s+1, l} + \sum_j \| G \Phi_{u, j} P v\|_{s-m+2, l} \\
        & \qquad \hspace{3em} + \| G \tau^{\mu} P v\|_{s-m+2, l} + \| G P v\|_{s-m+1, l} + \| E v \|_{s, l}) + C_E \| E v \|_{s-1/2, l} \\
        & \qquad \lesssim\sum_j \| G \Phi_{u, j} P v\|_{s-m+2, l} + \| G\tau^{\mu} P v\|_{s-m+2, l} + \| G P v\|_{s-m+1, l} + \\
      & \qquad \hspace{3em} \delta^{1/2 - \beta} (\| \Op(\chi_{\delta}) v \|_{s, l} + \| \Op(\chi_{\phi_u \neq 0}) v \|_{s, l} + \| \Op(\chi_{\tau \neq 0}) v \|_{s, l} + \| G Pv \|_{s-m+1, l} \\ 
      & \qquad \hspace{6em}  + C_E \| E v \|_{s-1/2, l}) + \| \Op(\chi_{\phi_u \neq 0}) v \|_{s+1, l} + \| \Op(\chi_{\tau \neq 0}) v \|_{s+1, l} + C_E \| E v \|_{s-1/2, l}
\end{align*}
Taking $\beta<1/2$ and $\delta$ small enough, we can absorb the $\| \chi_{\delta} v \|_{s-1}$ term to the left. As usual, we can improve the $s-1/2$ reminder error to any $-N$ by iterating the argument. This proves the estimate for $v$ sufficiently regular so that the pairings done in the steps above makes sense. To obtain the strong form of the inequality, i.e. conclude regularity of $v \in H^{r, l}$ at the trapped set from the right hand side being finite, one must thus regularize $v$ and then take a limit to obtain the result. As usual, this ammounts to adding a smoothing weight $(1 + \epsilon \rho)^{-k}$ in the commutators, and then taking the $\epsilon \rightarrow 0$ limit. Due to $H_p \rho \neq 0$ in general, there is some restriction on how much regularization we can do. More precisely, such a smoothing weight adds an extra term
\begin{equation}
    -k \frac{H_p(\rho)}{\rho} \epsilon (1+\epsilon \rho)^{-k} \times a^2
\end{equation}
to $H_p((1 + \epsilon \rho)^{-k} a^2)$, and similarly for $H_{\phi_u}((1 + \epsilon \rho)^{-k} a^2)$. A quick review of the argument shows that it can be absorbed as long as $k$ satisfies 

\begin{align*}
    &((2s - m + 1 - k)g - 2lf - 2p_1 + \nu_s)|_{\Gamma} > 0 \\
    &((2s - m + 3 - k)g - 2lf - 2p_1 + 2 \nu_u)|_{\Gamma} > 0
\end{align*}
We get no conditions from the propagation along $H_{\phi_u}$ since by picking $\delta$ small enough we can dominate the extra term. In order for the pairings to be defined, this corresponds to assuming $v \in H^{t, l}$ microlocally for $t > s_0$, $s_0$ being
\begin{equation}
\begin{split}
    & s_0 = inf_{s \in \mathbb{R}} \{ ((2s - m + 1)g - 2lf - 2p_1 + \nu_s)|_{\Gamma} > 0, \\
    & \qquad ((2s - m + 3)g - 2lf - 2p_1 + 2 \nu_u)|_{\Gamma} > 0
    \}
\end{split}
\end{equation}
as was claimed in the proposition. Then, using uniform estimates over $\epsilon$ and taking the limit $\epsilon \rightarrow 0$ gives the estimate. We refer to Hintz's notes on microlocal analysis [\cite{PeterNote} Section 9.5] for a more thorough description of the regulization process. Similarly, the claim for $P^{\dag}$ follows by the exact same proof applied to $\phi_s$ and propagating backward along the $H_p$ flow, except now we can neglect the $\chi_{\tau}'$ derivative terms as they share the same sign as the control terms. $\Box$

\section{The Fredholm Theory and Examples}\label{sec3}

Equipped with this inequality, we can develop the Fredholm theory of a non-elliptic operator $P$ in the presence of normally hyperbolic trapping. Contrary to the elliptic case, obtaining Fredholm estimate relies on global properties of the flow. More precisely, one pairs propagation estimates along with local a priori control of solutions $\| v \| \lesssim \| Pv \| + \| Rv \|$, $R$ a compact error, to get a global Fredholm estimate. As long as the flow $\varphi_t$ of $H_p$ reaches the region of a priori control $\mathcal{Q}$ in which we have this local a priori control, propagation results are sufficient to obtain an estimate on the whole phase space. Unfortunately, this is not always the case. In general, we expect the presence of a \textit{trapped set} $\Gamma$, which is exactly the points which do not reach $\mathcal{Q}$. It decomposes into backward/forward trapped sets $\Gamma_{u/s}$, which are roughly defined as follows:

\begin{equation}
\begin{split}
    \Gamma_{u/s} &= \{ (x, \xi) \in T^{\ast}M: \varphi_t(x, \xi) \not\to \mathcal{Q} \quad \text{as $t \rightarrow \mp \infty$} \} \\
    \Gamma &= \Gamma_u \cap \Gamma_s
\end{split}
\end{equation}
In order to obtain a Fredholm estimate, one must also control the solution at $\Gamma$. This is precisely the content of theorems \ref{TheoremClosed} and \ref{TheoremB}, which tells us that this is possible if $\Gamma$ is NHT and we localize away from the backward trapped set $\Gamma_u$ using defining functions $\phi_u$ for it. This localization procedure can be directly included into our function spaces by defining for the $b$-case (the compact case is similar except we don't have the $\tau^{\mu}$ term nor the decay order $l$)

\begin{equation}
\label{DefinitionSpaces}
    \begin{split}
    \mathcal{Y}^{s, l} &= \{ v \in H_b^{s, l}: \tau^{\mu}v, \forall j \text{ } \Phi_{u, j} v \in H_b^{s+\alpha, l} \quad (P-iQ)v \in H^{s-m+\alpha, l}_b\} \\
    \mathcal{X}^{s, l} &= \{ v \in H_b^{s, l}: \tau^{\mu}v, \forall j \text{ } \Phi_{u, j} v \in H_b^{s+\alpha, l} \quad (P-iQ)v \in \mathcal{Y}^{s-m+1, l}\} 
    \end{split}
\end{equation}
where $0 \leq \alpha \leq 1$ is a smooth function on $S_b^{\ast}M$ equal to 1 near $\Gamma$ and should be thought off as the improvement in regularity of $\Phi_u v$ relative to $v$, and $Q$ is a pseudodifferential operator of order $m$ representing the imaginary part of the operator. In the case of \textit{stationary} operators $P = P|_{\tau = 0}$, we also drop the $\tau^{\mu} v$ term in these spaces since then one does not need to choose extensions of $\phi_u$. These spaces are equipped with the graph norms \footnote{The $\| (P-iQ) v \|_{s-m+\alpha}^2$ is added to the norm of $\mathcal{Y}^{s,l}$ in order to impose weaker regularity at $q^{-1}(0) \cap p^{-1}(0) \cap \phi_u^{-1}(0) = \Gamma_u$ rather than just $\phi_{u}^{-1}(0)$. It guarantees that these spaces are independent of the extension $\phi_u$ outside the characteristic set.} 

\begin{equation}
    \begin{split}
    \| v \|_{\mathcal{Y}^{s, l}}^2 &= \| v \|_{s, l}^2 + \| \tau^{\mu} v \|_{s+\alpha, l}^2 + \sum_j \| \Phi_{u, j} v \|_{s+\alpha, l}^2 + \| (P-iQ) v \|_{s+\alpha-m, l}^2 \\
    \| v \|_{\mathcal{X}^{s, l}}^2 &= \| v \|_{s, l}^2 + \| \tau^{\mu} v \|_{s+\alpha, l}^2 + \sum_j \| \Phi_{u, j} v \|_{s+\alpha, l}^2 + \| (P-iQ)v \|_{\mathcal{Y}^{s-m+1, l}}^2
    \end{split}
\end{equation}
Spaces of this form have been studied before in the literature and are rather standard. For example, they show up in the study of diffraction, see \cite{melrose2011diffractionsingularitieswaveequation}. They are called \textit{coisotropic} \footnote{The name comes from the fact that these are usually considered when $\Phi_{u, j}$ are defining functions for coisotropic submanifolds in $T^{\ast}M$, as is the case here.} function spaces, and the property that $\Phi_u v$ is more regular than $v$ is called \textit{coisotropic regularity}. While slightly exotic, they are quite tame. For example, standard PSDOs act on them in the expected way $A: \mathcal{Y}^{s, l} \rightarrow \mathcal{Y}^{s-ord(A), l}$ and they contain smooth functions as a dense subspace. This is discussed in the appendix \ref{AppendixA}. By using our estimate \ref{TheoremClosed} at the trapped set and propagating it along the whole of $\Gamma_u$, we would obtain an inequality of the form

\begin{equation}
    \| v \|_{\mathcal{X}^{s, l}} \lesssim \| (P-iQ) v \|_{\mathcal{Y}^{s-m+1, l}} + \| v \|_{-N, l}
\end{equation}
which is exactly a Fredholm estimate for $P$ modulo improving the decay order of $\| v \|_{-N, l}$, which is done by inverting the normal operator as described in section \ref{bCalculusIntro}. The adjoint estimate for $P^{\dag}$ requires a bit more care as taking dual does not switch the role of $\phi_u$ and $\phi_s$ in these spaces, but can be handled by a regularity argument. All and all, this suggest the rough statement in theorem \ref{AlmostFredholmEstimates}. In the following subsection, we make this precise by defining what we mean by a 'region of a priori control $\mathcal{Q}$' and the required dynamics. Roughly speaking, $\mathcal{Q}$ will be somewhere where we can estimate $v$ by $Pv$ modulo a compact error. The three main ways this can be done is by using boundary conditions, complex absorption, i.e. a non-trivial imaginary part of the principal symbol, and radial point estimates.
 
\subsection{Setup of the Fredholm Theory}
\label{FredholmFunctionSpaces}

As described in the previous section, we will now describe precisely the Fredholm theory of non-elliptic operators with normally hyperbolic trapping by listing the assumptions required. We will consider operators of the forms $P - iQ$, where $P$ satisfies the assumptions of theorems \ref{TheoremClosed} or \ref{TheoremB} and $Q$ is a pseudo-differential operator with non-negative real principal symbol of order $m$. We start by describing the \textit{region of a priori control} $\mathcal{Q}$ where we have locally an estimate of the form $\| v \| \lesssim \| Pv \| + \| Rv \|$. We do so for the $b$-pseudo differential operators, but the definition is completely analogous in the case of a closed manifold, except we don't work at the boundary and there is no normal operator $N(P-iQ)$. One of the ways we will obtain a priori control is through radial sets \ref{RadialPointsB} which impose some constraints on the regularity order $s$ at the radial set. As one might have different radial sets in phase space, each with different constraints on $s$, it often happens that a constant $s$ cannot satisfy all these constraints simultaneously. Thus, one must let the regularity order $s$ vary for the radial set constraints to be all satisfied, which forces us to use \textit{variable order Sobolev spaces} $H^{s,l}_b$. These are quite common in non-elliptic problems; see [Section 5.3.9 \cite{VariableOrderSobolev}] for a reference on how they are constructed. The general idea is to replace the usual Sobolev norm by

\begin{equation}
    \| v \|_{s, l} \| \tau^{-l} \Op(\rho^{s}) v \|_{L^2_b} 
\end{equation}
for some smooth function $s$ on $S^{\ast}_b M$. One can show that this still forms a calculus, with the main caveat being that we lose a logarithmic factor when taking commutators. Furthermore, the propagation of singularity theorem \ref{PropagationOfSingularities} will still hold as long as $H_p(s) \leq 0$ and we propagate forward along the $H_p$ flow. If $H_p(s) \geq 0$, then one must propagate backward. So from now on, we assume that $s$ varies on $S^{\ast}_b M$. Furthermore, we fix $l \in \mathbb{R}$ and $z \in \mathbb{C}$ with $Im(z) = O(h)$ for the reminder of the discussion. As for the treshold regularity $s_0$, since it will be different at each possible radial sets, we will also allow for $s_0$ to be a function on $S^{\ast}_b M$, with the understanding that it equals the required (constant) threshold at each radial set.

\begin{definition}\textbf{A priori Control Region $\mathcal{Q}$:}
 We denote by $\mathcal{Q} = \mathcal{Q}_{P} \cup \mathcal{Q}_{P^{\dagger}}$ the region of a priori control, where $\mathcal{Q}_{P}$ is an open set in $T_b^{\ast}M$ such that for all $(x, \xi) \in \mathcal{Q}_{P}$, all $N \in \mathbb{R}$ such that $-N + 2 < s$ and all $s$ such that $s(x, \xi) > s_0(x, \xi)$, there exists some $0^{th}$ order pseudo-differential operator $A$ elliptic at $(x, \xi)$ satisfying

\begin{align}
    \| Av \|_{s+\alpha, l} &\lesssim \| (P - i Q) v \|_{s+\alpha-m+1, l} + \| v \|_{-N, l} \\
    \| Av \|_{H^{s+1}_h} &\lesssim h^{-1} \| \widehat{N(P - i Q)}_h v \|_{H^{s-m+2}_h} + o(1) \| v \|_{H^{-N}_h}    
\end{align}
 and $\mathcal{Q}_{P^{\dagger}}$ is an open set in $T_b^{\ast}M$ such that for all $(x, \xi) \in \mathcal{Q}_{P^{\dagger}}$, all $N \in \mathbb{R}$ such that $-N+2 < -s+m$ and all $s$ such that $s(x, \xi) < s_0(x, \xi)$, there exists some $0^{th}$ order pseudo-differential operator $B$ elliptic at $(x, \xi)$ satisfying

\begin{align}
    \| Bv \|_{-s+m-1, -l} &\lesssim \| (P^{\dag} + i Q) v \|_{-s+1, -l} + \| v \|_{-N, -l} \\
    \| Bv \|_{H^{-s+m-1}_h} &\lesssim h^{-1} \| \widehat{N(P - i Q)}_h^{\dagger} v \|_{H^{-s+1}_h} + o(1) \| v \|_{H^{-N}_h}
\end{align}
\end{definition}
Note that in the semiclassical estimates above, it should be understood that $A$ and $B$ are now semiclassical pseudo-differential operators with the aforementioned properties. The reason for the asymmetry between $P$ and $P^{\dagger}$ in these definitions is because the dual of our coisotropic spaces defined above cannot be expressed in terms of $\phi_s$, hence we obtain a lossy estimate for $P^{\dagger}$ relative to those spaces. Furthermore, we ask for $-N +2 < s$ or $-N + 2 < -s + m$ because we will be losing two derivatives in the reminder term $\| v \|_{-N, -l}$ when using the normal operator to improve it. Note in particular that these assumptions are satisfied by elliptic estimates for

\begin{equation}
    \mathcal{Q} = \{ q \neq 0 \}
\end{equation}
for all $s_0$ or for neighborhoods of radial sets \ref{RadialPointsB} acting as sources, where the threshold $s_0$ is given by $(m-1)/2 + l\frac{f}{g} \pm \frac{p_1}{g}$. In the presence of vanishing boundary conditions, one can also take $\mathcal{Q}$ to be the region past the boundary since $v$ vanishes there. Furthermore, in practice, the Hamiltonian flow of $H_p$ will contain different type of radial sets. At the sources, we will have unconditional control of $v$ by $Pv$ hence they will be included in $\mathcal{Q}_P$. On the other end, at the sinks, we will be using the low regularity versions of theorem \ref{RadialPointsB} or \ref{HyperbolicRadialSets} to propagate estimates to it. The sinks are thus not a part of $\mathcal{Q}_P$ hence must be treated separately. Similarly, when deriving the estimate for $P^{\dagger}$, the reverse is true since we will be propagating backward along the flow of $H_p$: the sinks of $H_p$ will be included in $\mathcal{Q}_P^{\dagger}$ but not the source. The purpose of the next definition is thus to include radial sets within our abstract framework, where one should think of the set $\mathcal{R}$ as an abstract definition of a radial set:

\begin{definition}\textbf{Abstract definition of radial sets:} We denote by $\mathcal{R} = \mathcal{R}_{P} \cup \mathcal{R}_{P^{\dagger}}$ the abstract radial sets, where $\mathcal{R}_{P}$ is an $H_p$-invariant closed homogeneous submanifold of $T_b^{\ast}M|_{\partial M}$ such that for all neighborhoods $U$ of $\mathcal{R}_{P}$, all s satisfying $s|_{\mathcal{R}_P} < s_0|_{\mathcal{R}_P}$ and all $N \in \mathbb{R}$ with $-N+2 < s$, there exists some $0^{th}$ order pseudo-differential operator $A$ elliptic at $\mathcal{R}_{P}$ and some $0^{th}$ order pseudo-differential operator $\Op(\chi_{\phi \neq 0})$ with wavefront set disjoint from $\mathcal{R}_P$ but contained in $U$ satisfying

\begin{align}
    \| Av \|_{s+\alpha, l} &\lesssim \| (P - i Q) v \|_{s+\alpha-m+1, l} + \| \Op(\chi_{\phi \neq 0}) v \|_{s+\alpha, l} + \| v \|_{-N, l} \\
    \| Av \|_{H^{s+1}_h} &\lesssim h^{-1} \| \widehat{N(P - i Q)}_h v \|_{H^{s-m+2}_h} + \| \Op(\chi_{\phi \neq 0}) v \|_{s+1, l} + o(1) \| v \|_{H^{-N}_h}    
\end{align}
and $\mathcal{R}_{P^{\dagger}}$ is an $H_p$-invariant closed homogeneous submanifold of $T_b^{\ast}M|_{\partial M}$ such that for all neighborhoods $U$ of $\mathcal{R}_{P^{\dagger}}$, all s satisfying $s|_{\mathcal{R}_{P^{\dagger}}} > s_0|_{\mathcal{R}_{P^{\dagger}}}$ and all $N \in \mathbb{R}$ with $-N+2 < -s+m$, there exists some $0^{th}$ order pseudo-differential operator $B$ elliptic at $\mathcal{R}_{P^{\dagger}}$ and some $0^{th}$ order pseudo-differential operator $\Op(\chi_{\phi \neq 0})$ with wavefront set disjoint from $\mathcal{R}_{P^{\dagger}}$ but contained in $U$ satisfying

\begin{align}
    \| Bv \|_{-s+m-1, -l} &\lesssim \| (P^{\dag} + i Q) v \|_{-s+1, -l} + \| \Op(\chi_{\phi \neq 0}) v \|_{-s+m+1, l} + \| v \|_{-N, -l} \\
    \| Bv \|_{H^{-s+m-1}_h} &\lesssim h^{-1} \| \widehat{N(P - i Q)}_h^{\dagger} v \|_{H^{-s+1}_h} + \| \Op(\chi_{\phi \neq 0}) v \|_{-s+m+1, l} + o(1) \| v \|_{H^{-N}_h}
\end{align}
\end{definition}
One should think of $\mathcal{R}_P$ as the sinks of the $H_p$ flow and $\mathcal{R}_{P^{\dagger}}$ as its sources. Finally, the last possibility we wish to include in our theory is the presence of \textit{saddle} radial sets \ref{HyperbolicRadialSets}, which have opposite behavior along the boundary and transverse to it. These come up in the Kerr-de Sitter family of black holes for example which will be discussed in section \ref{SectionExamples}. They play a similar role as $\mathcal{Q}$, in the sense that they offer some a priori control. However, they are different since they require some control outside the boundary to play that role, as the estimate we have there is theorem \ref{HyperbolicRadialSets}

\begin{equation}
    \| \chi v \|_{s, l} \lesssim \| Pv \|_{s-m+1, l} + \| \Op(\chi_{\tau \neq 0}) v \|_{s, l} + \| v \|_{-N, l}
\end{equation}
Thus, one must make an additional propagation estimate in order to estimate $ \| \Op(\chi_{\tau \neq 0}) v \|_{s, l}$ hence get something of the form $\| v \| \lesssim \| Pv \| + \| R v \|$. This step can be thought off as propagating 'through' the saddle point, since we are basically ignoring its presence. Abstractly, what we need is that one can propagate the estimate on $\mathcal{Q}$ to these saddle points, and then we can propagate it everywhere else. This is the content of the next definition:

\begin{definition}\textbf{Abstract definition of saddle radial sets:} \label{AbstractSaddleRadialSets} We denote by $\mathcal{L} = \mathcal{L}_{P} \cup \mathcal{L}_{P^{\dagger}}$ the abstract saddle radial sets, where $\mathcal{L}_{P}$ is an $H_p$-invariant closed homogeneous submanifold of $T_b^{\ast}M|_{\partial M}$ such that for all neighborhoods $U$ of $\mathcal{L}_{P}$, all s satisfying $s|_{\mathcal{L}_P} > s_0|_{\mathcal{L}_P}$ and all $N \in \mathbb{R}$ with $-N+2 < s$, there exists some $0^{th}$ order pseudo-differential operator $A$ elliptic at $\mathcal{L}_{P}$ and some $0^{th}$ order pseudo-differential operator $\Op(\chi_{\phi \neq 0})$ with wavefront set disjoint from $\mathcal{L}_P$ but contained in $U$ satisfying

\begin{align}
    \| Av \|_{s+\alpha, l} &\lesssim \| (P - i Q) v \|_{s+\alpha-m+1, l} + \| \Op(\chi_{\phi \neq 0}) v \|_{s+\alpha, l} + \| v \|_{-N, l} \\
    \| Av \|_{H^{s+1}_h} &\lesssim h^{-1} \| \widehat{N(P - i Q)}_h v \|_{H^{s-m+2}_h} + \| \Op(\chi_{\phi \neq 0}) v \|_{s+1, l} + o(1) \| v \|_{H^{-N}_h}    
\end{align}
Furthermore, we assume that the backward integral curves starting from $WF(\Op(\chi_{\phi \neq 0}))$ reach $\mathcal{Q}_{P}$ in finite time. Similarly, $\mathcal{L}_{P^{\dagger}}$ is an $H_p$-invariant closed homogeneous submanifold of $T_b^{\ast}M|_{\partial M}$ such that for all neighborhoods $U$ of $\mathcal{L}_{P^{\dagger}}$, all s satisfying $s|_{\mathcal{L}_{P^{\dagger}}} < s_0|_{\mathcal{L}_{P^{\dagger}}}$ and all $N \in \mathbb{R}$ with $-N + 2 < -s+m$, there exists some $0^{th}$ order pseudo-differential operator $B$ elliptic at $\mathcal{L}_{P^{\dagger}}$ and some $0^{th}$ order pseudo-differential operator $\Op(\chi_{\phi \neq 0})$ with wavefront set disjoint from $\mathcal{L}_{P^{\dagger}}$ but contained in $U$ satisfying

\begin{align}
    \| Bv \|_{-s+m-1, -l} &\lesssim \| (P^{\dag} + i Q) v \|_{-s+1, -l} + \| \Op(\chi_{\phi \neq 0}) v \|_{-s+m+1, l} + \| v \|_{-N, -l} \\
    \| Bv \|_{H^{-s+m-1}_h} &\lesssim h^{-1} \| \widehat{N(P - i Q)}_h^{\dagger} v \|_{H^{-s+1}_h} + \| \Op(\chi_{\phi \neq 0}) v \|_{-s+m+1, l} + o(1) \| v \|_{H^{-N}_h}
\end{align}
and for which the forward integral curves starting from $WF(\Op(\chi_{\phi \neq 0}))$ reach $\mathcal{Q}_{P^{\dagger}}$ in finite time.
\end{definition}
So effectively, after propagating the estimate from $\mathcal{Q}$ to $\mathcal{L}$, one can think of a neighborhood of $\mathcal{L}$ as also being a region of a priori control where we have $\| v \| \lesssim \| Pv \| + \| Rv \|$. With that, the \textit{global} dynamical assumptions that we must assume for the Fredholm theory are the following. 

\begin{definition}\textbf{Global Normally Hyperbolic Trapping:}\label{GlobalNHT} Let $\varphi_t$ be the Hamiltonian flow of $H_p$ in $T^{\ast}_b M$ and define the backward and forward trapped set $\Gamma_{u/s} \subset T^{\ast}_b M|_{\partial M}$ by

\begin{equation}
\begin{split}
    \Gamma_{u/s} = &\{ (x, \xi) \in p^{-1}(0) \cap q^{-1}(0) \cap T_b^{\ast}M|_{\partial M}: \\
    &\varphi_t(x, \xi) \text{ does not reach $\mathcal{Q}_{P/P^{\dagger}}$ or a neighborhood of $\mathcal{L}_{P/P^{\dagger}}$ in a finite negative/positive time} \} \backslash (\mathcal{R} \cup \mathcal{L})
\end{split}
\end{equation}
Suppose that $\Gamma = \Gamma_u \cap \Gamma_s$ is a normally hyperbolic invariant manifold in the sense of section \ref{IntroTrapping} and that the following global dynamical assumptions are satisfied:
    
    \begin{itemize}
         
        \item All backward integral curves of $H_p$ from $p^{-1}(0) \cap q^{-1}(0) \backslash (\Gamma_{u} \cup \mathcal{R}_P \cup \mathcal{L}_{P})$ reach $\mathcal{Q}_{P}$ in finite time or tend to $\mathcal{L}_P$ while all forward integral curves of $H_p$ from $p^{-1}(0) \cap q^{-1}(0) \backslash (\Gamma_{s} \cup \mathcal{R}_{P^{\dagger}} \cup \mathcal{L}_{P^{\dagger}})$ reach $\mathcal{Q}_{P^{\dagger}}$ in finite time or tend to $\mathcal{L}_{P^{\dagger}}$. 

        \item The backward/forward integral curves of $H_p$ along $\Gamma_{u/s}$ tend to $\Gamma$.   

        \item There is a smooth function $0 \leq \alpha \leq 1$ on $S_b^{\ast}M$ which is equal to 1 in a neighborhood of $\Gamma$, is non-increasing along the Hamiltonian flow $H_p(\alpha) \leq 0$ and such that $\alpha|_{q^{-1}(0, \infty) \cap \Gamma_u} < 1/2$. In words, $\alpha$ has decreased along $\Gamma_u$ by at least $1/2$ once it reaches the absorbing region $\{ q > 0 \}$ \footnote{This assumption may seem strange. It is due to a loss of coisotropic regularity when crossing the boundary of $q = 0$. See appendix \ref{AppendixA}. Furthermore, it allows us to 'cutoff' $\Gamma_u$ before it reaches any problematic regions in phase space, say $\mathcal{R}$ or $\mathcal{L}$, where the $\Phi_u v$ term could affect the radial point estimates.}.

        \item The unstable manifold $\Gamma_u$ has defining functions $\phi_u$ in the set $\{ \alpha \neq 0 \} \cap q^{-1}(0)$. More precisely, there are global functions $\phi_u$ on $S^{\ast}_b M|_{\partial M}$ such that 
        
        \begin{equation}
            \phi_{u}^{-1}(0) \cap \{ q < \epsilon \} \cap \{ \alpha \neq 0 \} \cap p^{-1}(0) \cap S_b^{\ast}M|_{\partial M}  = \Gamma_{u/s}
        \end{equation}

        for some $\epsilon > 0$.

       \item The flow is non-trapping in $\{ \sigma = 0 \}$, that is $\sigma$ does not vanish on $\Gamma_{u/s}$ \footnote{One can do without this assumption by using a different homogeneous of degree 1 function $\rho$ then $\sigma$ in the trapping estimate. We do not expand on that because in practice this is satisfied.}.
       
    \end{itemize} 
    
\end{definition}
\textit{Remarks:} For our applications to waves, the characteristic set $p^{-1}(0)$ of the operator will have multiple (2 in our case) connected components. For such operators, one must make a choice in each connected components about which direction along $H_p$ one wants to propagate estimates: either along $+ H_p$ or $-H_p$. For example, in the case of initial value problems for waves, one wants to propagate forward along $H_p$ for future-directed covectors, and backward along $H_p$ for past-directed covectors. For those operators, the assumptions above should then be satisfied for each connected components of $p^{-1}(0)$, where $H_p$ is replaced by $\pm H_p$ depending on which sign we chose for that component. Furthermore, one must take $H_p(s) \geq 0$ in the connected components in which we are propagating backward along $H_p$. Once this is done, the rest carries through.

\vspace{0.2in}
Note that in order to define the coisotropic spaces in \ref{DefinitionSpaces} we need to quantize global defining functions $\phi_u$ for $\Gamma_u$, which are now global submanifolds of $T^{\ast}_bM|_{\partial M}$ rather than defined in a neighborhood of $\Gamma$ only. The issue is that the full $\Gamma_u$ may fail to have defining functions. For example, it might have a boundary once we take its closure. Thus we must work with a localized version usually. Furthermore, one must have techniques to check the global assumptions above without having to solve explicitely for the flow, which can be quite difficult. This is adressed in appendix \ref{AppendixB} and is the reason for the last assumption in definition \ref{GlobalNHT} since it can be shown that in reasonable cases we can construct $\alpha$ so that $\Gamma_u$ has a defining function in $\{ \alpha \neq 0 \}$. Roughly speaking, this is done by propagating the local unstable manifold in a neighborhood of $\Gamma$ by some finite time $T$ along the flow of $H_p$, using the fact that finite time propagation preserves the existence of defining functions. See proposition \ref{TrappingConclusion} for the details. For $q \geq \epsilon$, we may take arbitrary (homogeneous of degree 0) extensions, which ammounts to extending $\Gamma_u$ in an arbitrary way. Furthermore, once again, one should think of the finite frequency regions $\mathcal{Q} \cap \{ \sigma = Re(z) \}$ and $\Gamma \cap \{ \sigma = Re(z) \}$ as being relevant to the semiclassical problem, while the regions at 'infinite' frequencies $\mathcal{Q} \cap S^{\ast}_b M$ and $\Gamma \cap S^{\ast}_b M$ corresponds to the dynamics for $P-iQ$. Once we have made these assumptions, they are sufficient to show that $P-iQ$ is Fredholm on coisotropic spaces whose decay  order $l$ is bounded. 

\begin{theorem}
\label{FredholmComplexPotential}
    Let M be a compact manifold with boundary and $P - iQ \in \Psi_b^{m}$ a $b$-pseudo-differential operator whose principal symbol has a non-degenerate homogeneous of degree $m$ real part $p$ and whose imaginary part $q$ is $\geq 0$ with $Q = Q^{\dag}$. Assume furthermore that it satisfies global hyperbolic trapping as in definition \ref{GlobalNHT} and that $P - iQ$ has a normal operator $(P - iQ)|_{\tau = 0}$ for which $(P - iQ) - (P - iQ)|_{\tau = 0} = O(\tau^{\mu})$. Suppose furthermore that the decay order $l$ satisfies at $\Gamma$

    \begin{equation}
        l < \min(\mu, \nu/2) - \frac{p_1}{f}|_{\Gamma}
    \end{equation}
    where $f = \Tilde{H}_p(\tau) / \tau$ and $p_1$ is the subprincipal symbol of $P$ divided by $\rho^{m-1}$. Suppose also that the indicial family $\widehat{N(P-iQ)}(\sigma)$ of the normal operator has no poles on the line $Im(\sigma) = -l$. Then $P - iQ: \mathcal{X}^{s, l} \rightarrow \mathcal{Y}^{s-m+1, l}$ is a Fredholm operator for $s$ and $s-2$ satisfying the constraints of $\mathcal{Q}$, $\mathcal{R}$, $\mathcal{L}$ with $H_p(s) \leq 0$.
\end{theorem}
\textit{Remarks:} The semiclassical version of the estimate, equation \ref{SemiclassicalNHT_Estimate}, and the non-trapping assumption in $\sigma = 0$ will show that the poles of the indicial family are discrete hence for all $l$ but a discrete set there is no poles on the line $Im(\sigma) = -l$. Furthermore, there is no restriction on $s$ coming from the NHT since we are assuming that $\sigma \neq 0$ on $\Gamma$ hence we can take $g = 0$ in our trapping estimate \ref{TheoremB} by choosing $\rho = \sigma$.
\vspace{0.2in}

\textbf{Proof:}  
We start by describing how to obtain the estimate for $P-iQ$ then briefly mentions how to show the claim that the indicial family $\widehat{N(P-iQ)(\sigma)}$ is a meromorphic family of Fredholm operator. Outside of the characteristic set $p^{-1}(0) \cap q^{-1}(0)$, the Fredholm estimate follows from elliptic regularity. Thus we focus on a neighborhood of $p^{-1}(0) \cap q^{-1}(0)$. Since $M$ is compact, so is $p^{-1}(0) \cap q^{-1}(0) \subset S^{\ast}M$ so let $\psi_i$ be a finite partition of unity for $p^{-1}([-\epsilon, \epsilon]) \cap q^{-1}[0, \epsilon) \backslash (V_u \cup V_{\mathcal{R}_P} \cup V_{\mathcal{L}_P})$ where $V_u$ is an open neighborhood of $\Gamma_u$, $V_{\mathcal{R}_P}$ an open neighborhood of $\mathcal{R}_P$ and $V_{\mathcal{L}_P}$ an open neighborhood of $\mathcal{L}_P$. First note that on $\mathcal{L}_P$, one can obtain unconditional control of $v$ by $(P-iQ)v$ by assumption since we can propagate the control of $\mathcal{Q}_P$ to it. So for all purposes, one can think of $V_{\mathcal{L}_P}$ as being part of $\mathcal{Q}_P$. By assumption, for small enough $\epsilon$, all backward integral curves from $supp(\psi_i)$ reach $\mathcal{Q}_P \cup \mathcal{L}_P$ in finite time. Since $q \geq 0$ and $H_p(\alpha) \leq 0$, propagation of regularity allows us to propagate estimates forward along the Hamiltonian flow of $p$ hence we obtain for some $WF(A) \subset \mathcal{Q}_P$ or $WF(A) \subset \mathcal{L}_P$

\begin{align*}
    \| \Op(\psi_i) v \|_{s+\alpha, l} &\lesssim \| A v \|_{s+\alpha, l} + \| (P - iQ) v \|_{s+\alpha-m+1, l} + \| v \|_{-N, l} \\
    &\lesssim  \| (P - iQ) v \|_{s +\alpha - m + 1, l} + \| v \|_{-N, l}
\end{align*}
To estimate $v$ on $V_u$, we will use that the backward flow on $\Gamma_u$ tends to $\Gamma$ hence we can propagate the estimate at $\Gamma$ of theorem \ref{TheoremB}. To do so, one must be able to propagate coisotropic regularity for operators with an imaginary part. This is shown to be true in the appendix proposition \ref{PropagationCoisotropic} as long as we stay within $q = 0$, while if we cross to $q > 0$ we lose a 1/2 regularity for $\Phi_u v$. By our assumption that the regularity order $s + \alpha$ of $\Phi_u v$ decreases by at least 1/2 once we reach $q > 0$, we can propagate to that region as well. Thus, given a partition of unity $\chi_j$ for $\Bar{V}_u$, backward integral curves starting from $supp(\chi_j)$ go to a neighborhood $U$ of $\Gamma$ in finite time hence we have for $WF(A) \subset U$

\begin{align*}
    & \| \Op(\chi_j) v \|_{s, l} + \sum_j \| \Op(\chi_j) \Phi_{u, j} v \|_{s+\alpha, l} + \| \Op(\chi_j) \tau^{\mu} v \|_{s + \alpha, l} \\
    & \qquad \lesssim \| A v \|_{s, l} + \| (P- iQ) v \|_{s - m + 1, l} + \sum_j \| \Phi_{u, j} A v \|_{s+\alpha, l} + \\
    & \qquad \hspace{3em} \| \tau^{\mu} A v \|_{s + \alpha, l} + \| \Phi_{u, j} (P - iQ) v \|_{s+\alpha - m + 1, l} + \| \tau^{\mu} (P - iQ) v \|_{s - m + 1 + \alpha, l} + \| v \|_{-N, l} 
     \\
     & \qquad \lesssim \| (P - iQ) v\|_{s - m + 1, l} + \sum_j \| \Phi_{u, j} (P - iQ) v\|_{s + \alpha - m + 1, l} \\
     & \qquad \hspace{3em} + \| \tau^{\mu} (P - iQ) v\|_{s + \alpha - m + 1, l} + \| \Op(\chi_{\phi_u \neq 0}) v \|_{s + \alpha, l} + \| \Op(\chi_{\tau \neq 0}) v \|_{s + \alpha, l} + \| v \|_{-N, l}
\end{align*}
where we've also used that $\alpha = 1$ near $\Gamma$. The terms $\| \Op(\chi_{\phi_u \neq 0}) v \|_{s+\alpha} + \| \Op(\chi_{\tau \neq 0}) v \|_{s + \alpha, l}$ can be handled (by shrinking $V_u$ if required) using the $\psi_i$ estimate since it is supported away from $\Gamma_u$. As for estimating $v$ near the radial sets on $V_{\mathcal{R}_P}$, we use that at those one can propagate estimates from a punctured neighborhood of $\mathcal{R}_P$ to it. So take a partition of unity $\chi_j$ covering $\Bar{V}_{\mathcal{R}_P}$. On each of those, we have the estimate 

\begin{equation}
    \| \Op(\chi_j) v \|_{s+\alpha, l} \lesssim \| (P - i Q) v \|_{s+\alpha-m+1, l} + \| \Op(\chi_{\phi \neq 0}) v \|_{s+\alpha, l} + \| v \|_{-N, l}
\end{equation}
with $WF(\Op(\chi_{\phi \neq 0}))$ in $V_{\mathcal{R}_P}$ and disjoint from $\mathcal{R}_P$. By our global NHT assumptions, the backward integral curves from $WF(\Op(\chi_{\phi \neq 0}))$ reach $\mathcal{Q}_P$ in finite time hence one can propagate the control of $\mathcal{Q}_P$ to it. Thus, we have

\begin{equation}
    \| \Op(\chi_{\phi \neq 0}) v \|_{s+\alpha, l} \lesssim \| (P - i Q) v \|_{s+\alpha-m+1, l} + \| v \|_{-N, l}
\end{equation}
hence plugging that in the equation above gives the required estimate on $V_{\mathcal{R}_P}$ too. Then, patching all of these estimates together yields

\begin{equation}
\begin{split}
    & \| v \|_{s, l} + \sum_j \| \Phi_{u, j} v \|_{s+\alpha, l} + \| \tau^{\mu} v \|_{s+\alpha, l} \lesssim \| (P - iQ) v \|_{s - m + 1, l} \\
    & \qquad + \sum_j \| \Phi_{u, j} (P - iQ)v \|_{s + \alpha - m + 1, l} + \| \tau^{\mu} (P - iQ)v \|_{s + \alpha - m + 1, l} + \| v \|_{-N, l}
\end{split}
\end{equation}
which is exactly the global estimate $\| v \|_{\mathcal{X}^{s, l}} \lesssim \| (P - iQ) v \|_{\mathcal{Y}^{s - m + 1, l}} + \| v \|_{-N, l}$. In order to improve the decay order, we will use that the indicial family $\widehat{N(P-iQ)}(\sigma)$ of the normal operator has no poles on the line $Im(\sigma) = -l$. This implies that the normal operator $N(P-iQ)$ is invertible for that $l$ hence it satisfies an estimate of the form \footnote{Here we must take $-N$ to also satisfy the constraints of $\mathcal{Q}$, $\mathcal{R}$ and $\mathcal{L}$. This can be done while keeping $-N < s - 2$ by our assumption that $s-2$ satisfy the constraints too.}

\begin{equation}
    \| \chi v \|_{H^{-N, l}} \lesssim \| N(P - iQ) \chi v \|_{H^{-N - m + 2, l}}
\end{equation}
where $\chi$ is a cutoff equal to 1 near the boundary $\partial M$. Thus, we estimate the reminder term as

\begin{align*}
    \| v \|_{-N, l} & \lesssim \| \chi v \|_{-N, l} + \| (1 - \chi) v \|_{-N, l} \\
    & \lesssim \|N(P-iQ) \chi v \|_{H^{-N - m + 2, l}} + \| (1 - \chi) v \|_{-N, -M} \\
    & \lesssim \| (P - iQ) \chi v \|_{H^{-N - m + 2, l}} + \| ((P - iQ) - (P - iQ)|_{\tau = 0}) \chi v \|_{H^{-N - m + 2, l}} +  \| v \|_{-N, -M} \\
    & \lesssim \| (P - iQ) v \|_{H^{-N - m + 2, l}} + \| v \|_{H^{-N + 2, l + \mu}} 
\end{align*}
Combining both results, we get the desired Fredholm estimate

\begin{equation}
\label{FredholmEstimateComplexPot}
    \| v \|_{\mathcal{X}^{s, l}} \lesssim \| (P - iQ) v \|_{\mathcal{Y}^{s - m + 1, l}} + \| v \|_{-N+2, l + \mu}
\end{equation}
This shows that $P - iQ$ has closed ranged and finite dimensional kernel. In order to obtain that it has finite dimensional cokernel, note that by the exact same reasoning as above but now propagating backward along the $H_p$ flow and working with $\Gamma_s$ rather than $\Gamma_u$, we have the following estimate for the adjoint $P^{\dag} + iQ$

\begin{equation}
    \| v \|_{-s+m-1, -l} \lesssim \| (P^{\dag} + iQ) v \|_{-s+1, -l} + \| v \|_{-N+2, -l + \mu}
\end{equation}
This implies that 

\begin{equation}
    P - iQ : \{ v \in H_b^{s-1, l}: (P-iQ)v \in H_b^{s-m+1, l}\} \rightarrow H_b^{s-m+1, l}
\end{equation}
has a finite dimensional cokernel. It follows that $P-iQ: \mathcal{X}^{s, l} \rightarrow \mathcal{Y}^{s-m+1, l}$ also has finite dimensional cokernel since we can write for $f \in \mathcal{Y}^{s-m+1, l}$

\begin{equation}
    f = (P-iQ)v + g
\end{equation}
with $v \in H^{s-1, l}$ and $g \in \text{Coker}(P-iQ) \simeq (\text{Ker}(P^{\dag}+iQ))^{\dagger}$. Since $\text{Ker}(P^{\dag}+iQ) \subset H^{-s+m}$ the elements of the cokernel have representative in $H^{s-m+2} \supset \mathcal{Y}^{s-m+1}$ hence if $f \in \mathcal{Y}^{s-m+1, l}$ then so is $(P-iQ)v$. Now note that if $(P-iQ)v \in \mathcal{Y}^{s-m+1}$ the fact that inequality \ref{FredholmEstimateComplexPot} holds strongly means that $v$ is actually in $\mathcal{X}^{s, l}$. This shows that $Im(P-iQ) \subset \mathcal{Y}^{s-m+1}$ has finite codimension hence $P-iQ: \mathcal{X}^{s-m+1, l} \rightarrow \mathcal{Y}^{s-m+1, l}$ is Fredholm as desired. Finally, we conclude by quickly explaining why the indicial family $\widehat{N(P-iQ)}(\sigma)$ is a meromorphic family of Fredholm operator. To show that, one must check the assumptions of the analytic Fredholm theorem, that is check that $\widehat{N(P-iQ)}$ is a family of Fredholm operator which is invertible somewhere. As explained in section \ref{ThreeDynamicalSystems}, the dynamics for $\widehat{N(P-iQ)}(\sigma)$ are the ones of $H_p$ in the $\{ \sigma = 0 \}$ hyperplane hence because we assume that $\Gamma_{u/s}$ does not intersect it, we obtain non-trapping estimate for it by propagating the control on $\mathcal{Q}$ everywhere:

\begin{align}
    &\| v \|_{s} \lesssim \| \widehat{N(P - iQ)}(\sigma) v \|_{s - m + 1} + \| v \|_{-N} \\
    & \| v \|_{-s + m - 1} \lesssim \| \widehat{N(P - iQ)}^{\dagger}(\sigma) v \|_{-s} + \| v \|_{-N}
\end{align}
which shows that $\widehat{N(P-iQ)}(\sigma)$ is Fredholm for any $\sigma$ between

\begin{equation}
    \{ v \in H^{s}: \widehat{N(P-iQ)}(\sigma)v \in H^{s-m+1}\} \rightarrow H^{s-m+1}
\end{equation}
where the distributions $v$ should be thought off as defined on a tubular neighborhood $[0, 1) \times \partial M$ of the boundary. Note also that the space above is independent of $\sigma$ since it only depends on the principal symbol of $\widehat{N(P-iQ)}(\sigma)$, which is the same for all $\sigma$. As for the high-energy estimate, similarly as we did for $P-iQ$, propagating the control on $\mathcal{Q}$, now at finite frequency $\sigma = Re(z)$, and using the trapping estimate yields:

\begin{align}
    &\| v \|_{H^s_h} \lesssim h^{-2} \| \widehat{N(P - iQ)}_h v \|_{H_h^{s - m + 2}} + o(1) \| v \|_{H^{-N}_h} \\
    & \| v \|_{H^{-s + m - 1}_h} \lesssim h^{-2} \| \widehat{N(P - iQ)}^{\dagger}_h v \|_{H^{-s}_h} + o(1)\| v \|_{H^{-N}_h}
\end{align}
By taking $h \rightarrow 0$, which ammounts to letting $|Re(\sigma)| \rightarrow \infty$, the reminder term can be made negligible compared to the left-hand side which shows that for large $|Re(\sigma)|$, the operator $\widehat{N(P-iQ)}(\sigma)$ is invertible between the spaces

\begin{equation}
    \{ v \in H^{s}: \widehat{N(P-iQ)}(\sigma)v \in H^{s-m+2}\} \rightarrow H^{s-m+2}
\end{equation}
hence one can apply the analytic Fredholm theorem on that space, showing meromorphicity.

$\Box$

\vspace{0.2in}
It should be said that for the sake of simplifying the exposition, we did not cover the most general situation. For example, something that might happen is that there are multiple saddle points following one another along the flow, and we must propagate through each of them one at a time. Or also there could be more exotic 'radial sets' then the ones we've defined, which are not regular (this occurs for example in the study of gradient flows in \cite{dang2018spectralanalysismorsesmalegradient} and generally when studying vector fields). Colloquially, the Fredholm theory works in all cases where we can propagate the estimates from $\mathcal{Q}$ to everywhere else besides $\Gamma_{u/s}$, however this is done. Similarly, we have a Fredholm theory on closed manifold as well. The spaces are now defined as

\begin{equation}
    \begin{split}
    &\mathcal{Y}^{s} = \{ v \in H^{s}: \forall j \text{ } \Phi_{u, j} v \in H^{s+\alpha} \quad (P - iQ) v \in H^{s - m + \alpha} \} \\
    &\mathcal{X}^s = \{ v \in H^{s}: \forall j \text{ } \Phi_{u, j} v \in H^{s+\alpha} \quad (P - iQ) v \in \mathcal{Y}^{s-m+1} \} 
    \end{split}
\end{equation}
with the graph norms

\begin{equation}
    \begin{split}
    \| v \|_{\mathcal{Y}^{s}}^2 &= \| v \|_{s}^2 + \sum_j \| \Phi_{u, j} v \|_{s + \alpha}^2 + \| (P - iQ) \|_{s-m+\alpha}^2 \\
    \| v \|_{\mathcal{X}^{s}}^2 &= \| v \|_{s}^2 + \sum_j \| \Phi_{u, j} v \|_{s+\alpha}^2 + \| (P - iQ) v \|_{\mathcal{Y}^{s-m+1}}^2
    \end{split}
\end{equation}
As before, the global version of theorem \ref{TheoremClosed} is exactly the Fredholm estimate we need. Using the exact same reasoning as above, we obtain the following statement. Note that there is now a constraint on the regularity order $s$ coming from the trapping because there is no natural $\rho$ to choose from when using the trapping estimate \ref{TheoremClosed}, so $g = \Tilde{H}_p (\rho) / \rho$ is non-zero in general. 
 
\begin{theorem}
\label{FredholmClosedManifold}

    Let M be a closed manifold and $P - iQ \in \Psi^m(M)$ have principal symbol whose real part $p$ is non-degenerate and homogeneous of degree $m$ while its imaginary part $q$ is $\geq 0$. Suppose furthermore that it satisfies global normal hyperbolicity as in definition \ref{GlobalNHT} and that the subprincipal symbol $p_1$ of $P$  and the regularity order $s$ satisfies at $\Gamma$
    \begin{align}
         (p_1 + (s - \frac{(m-1)}{2})g)|_{\Gamma} &< \frac{\nu}{2} \\
         (p_1 + (s - \frac{(m-3)}{2})g)|_{\Gamma} &< \nu
    \end{align}
    Then $P - iQ: \mathcal{X}^{s} \rightarrow \mathcal{Y}^{s-m+1}$ is a Fredholm operator for s and s-1 satisfying the constraints of $\mathcal{Q}$, $\mathcal{R}$ and $\mathcal{L}$ with $H_p(s) \leq 0$.
\end{theorem}
\textit{Remark:} If the dynamics at $\Gamma$ are at most exponential, that is $\rho(\varphi_{|t|})|_{\Gamma} \lesssim e^{\epsilon |t|}$ for some $\epsilon > 0$ and some non-vanishing homogeneous of degree 1 function $\rho$, then an averaging argument produces a $g$ which satisfies $|g| \leq \epsilon$. Thus for the subexponential cases where this holds for any $\epsilon$, the restriction on $s$ is not significant.

\vspace{0.2in}
While it is possible to construct examples of NHT on closed manifolds relatively easily by gluing local dynamical systems exhibiting NHT to a closed manifold, natural examples of operators on closed manifolds satisfying these assumptions is harder to obtain because the typical dynamical systems do not have smooth trapped set. For example, the main class of operators one might consider are vector fields $-iX$. The requirement for them to be Fredholm in the non-trapping case is then that its flow has global hyperbolic attractor and repellors, with the key property being the exponentially growth/contraction of frequencies near the attractor/repellor. This is for example the case for Anosov flows \cite{dyatlov2016dynamical} with the stable and unstable bundle playing the role of attractor and repellor, and for the gradients of Morse-Smale functions \cite{dang2018spectralanalysismorsesmalegradient} where the union over the critical points of the stable and unstable bundles play the same role. In this context, a trapped set would then be a \textit{center bundle} where the pushforward of the flow $d \varphi_t$ is subexponential, say polynomial. For example, a bundle $E_c$ over some submanifold $N$ of $M$ which is stable under $\varphi_t$ and such that 

\begin{equation}
    \| d \varphi_t|_{E_c} \| \lesssim O((1+|t|)^N)
\end{equation}
Examples of these would be gradient flow with degenerate equilibrium point, where the center manifold theorem provided the existence of a non-unique center manifold [Section 3.2 \cite{CenterManifold}] on which the dynamics are determined by non-linear ODEs and subexponentials. The issue of course will be its regularity. In the gradient case for example, we'd have to take the union over all critical points of all the center manifolds, which isn't regular. It would thus be quite interesting to understand how to fit such vector fields in this framework and to come up with other natural class of operators on closed manifolds exhibiting NHT.

\subsection{Fredholmness of Wave Operators with Trapping}
\label{SectionExamples}

We now present our application of the Fredholm theory, which will highlight the global approach to wave equations on spacetimes. For a good introduction to the method, see for example \cite{baskin2014asymptotics} and its generalization in \cite{baskin2018asymptotics} where it was used to compute the asymptotics of waves on Minkowski-like spacetimes. In the case of black hole spacetimes, see the proof of the stability of Kerr deSitter \cite{KerrStab} which crucially used Fredholm estimates. We'd like to show that $O(\tau^{\mu})$ perturbations of wave operators $P$ associated to a pseudo-Riemannian manifold $(M, g)$ is Fredholm on weighted (in $\sim (t, r)$) Sobolev spaces. A pseudo-differential operator is called a \textit{wave operator} if its principal symbol $p$ is the inverse metric $g^{-1}(\xi, \xi)$, or $g^{-1}(\xi, \xi) Id$ in the case of vectorial wave equations. In order to put ourselves in a compact setting, we study the problem on a suitable \textit{compactification} $\Bar{M}$ of $M$. By doing so, we introduce a boundary to the problem hence we end up working with b-operators. We will cover two cases: the Kerr de-Sitter family and an example of a scattering spacetime exhibiting NHT. We start with the most important example, the Kerr-de Sitter family. The r-normally hyperbolic nature of the photon sphere for slowly rotating Kerr black holes was first described in \cite{Wunsch_2011}. It ended up being a key part in its study and was a major motivation for the development of the theory of normally hyperbolic trapping. For the global treatment of normally hyperbolic trapping in the full subextremal range of the Kerr-de Sitter metric, we refer to \cite{petersen2023waveequationskerrdesitter} for the details.

\begin{figure}[ht!]
    \centering
    \includegraphics[width=0.5\linewidth]{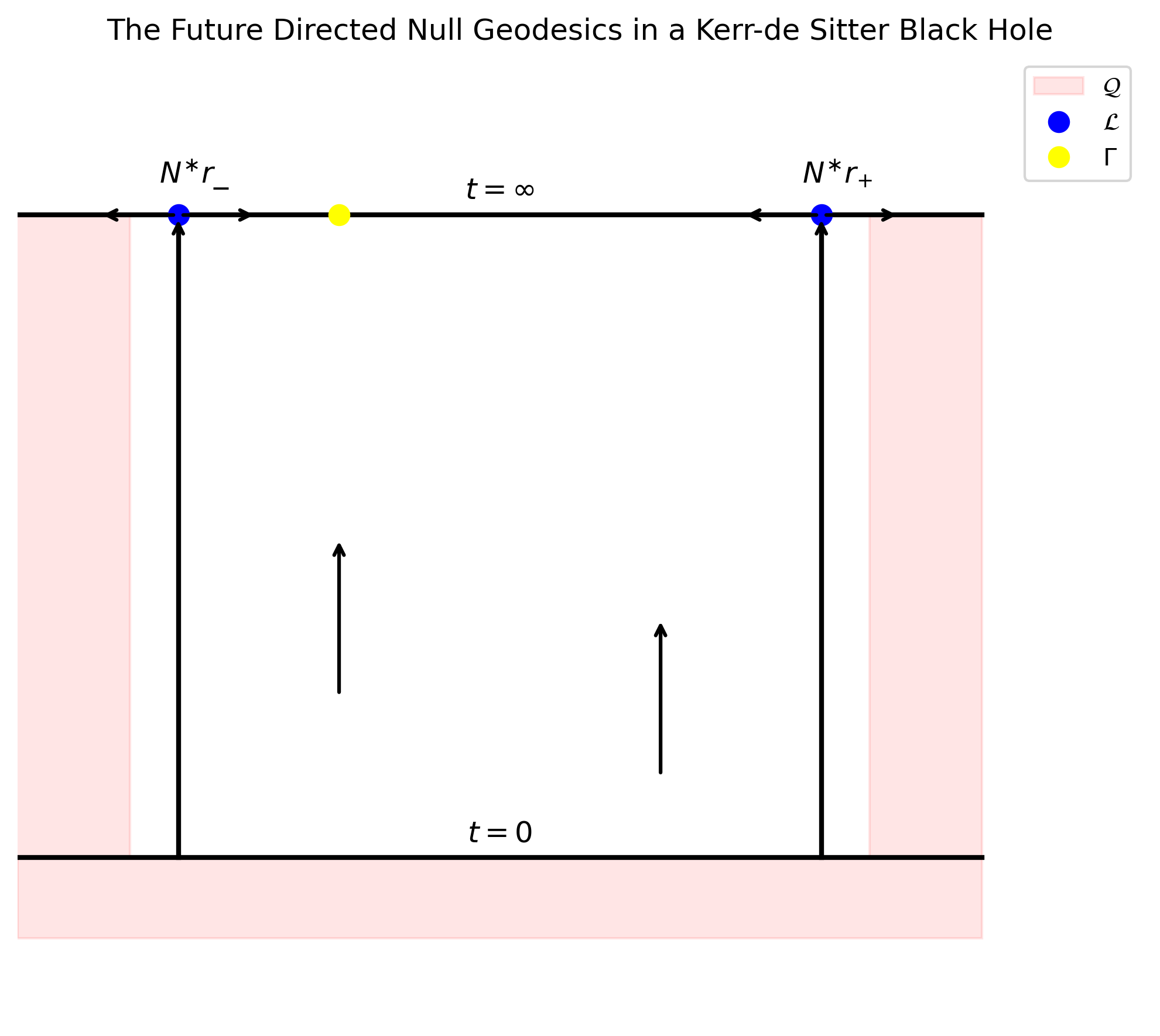}
    \caption{This is a 2D sketch of the dynamics within the $(r,t)$ plane of $T^{\ast}_b \Bar{M}$ for a subextremal Kerr-de Sitter black hole. For ease of visualization, the trapped set $\Gamma$ was represented as taking place on a singular radius $r$, here represented as a dot, but in the case of a non-zero angular momentum it projects to physical space as a range of radii.}
    \label{fig:KerrDynamics}
\end{figure}

\vspace{0.2in}
\textbf{Kerr-de Sitter Example:}
Let $(M, g)$ be the Kerr-de Sitter black hole space time, that is (in \textit{Boyer–Lindquist} coordinates)
\begin{equation}
\begin{split}
    &M = \mathbb{R}_t \times (r_{-}, r_{+}) \times \mathbb{S}^2 \\
    &g = \Sigma^2 \left( \frac{dr^2}{\mu} + \frac{d\theta^2}{\kappa} \right) + \frac{\kappa \sin^2(\theta)}{(1 + \lambda)^2 \Sigma^2}(a dt - (r^2 + a^2) d\phi)^2 \\
    & \qquad - \frac{{\mu_b}}{(1 + \lambda)^2 \Sigma^2} (dt - a \sin^2(\theta) d\phi)^2 \\
    &\mu = (r^2 + a^2)(1 - \frac{\Lambda r^2}{3}) - 2 M r \\
    &\Sigma^2 = r^2 + a^2 \cos^2(\theta), \quad \lambda = \frac{\Lambda a^2}{3}, \quad \kappa = 1 + \lambda \cos^2(\theta)
\end{split}
\end{equation}
where $(M, a)$ are parameters denoting the mass and angular momentum of the black hole, $\Lambda$ is the cosmological constant and  $r_{\pm}$ are the event horizons corresponding to two of the roots of $\mu_b$. We furthermore assumed that the latter has four distinct roots with $\mp \frac{\partial \mu}{\partial r} > 0$ at $r_{\pm}$ (this is called the \textit{subextremal case}). Just as was discussed in [Section 6 \cite{vasy2011microlocal}] (or [Section 3 \cite{KerrStab}]), we compactify this spacetime to a $b$-manifold in order to put ourselves in the setting of this paper. Thus, we instead work on the compact manifold with boundary

\begin{equation}
    \Bar{M} = [0, 1]_{\tau} \times \mathbb{S}^1_r \times \mathbb{S}^2 \quad \tau = e^{-t_{\ast}} 
\end{equation}
where the $\mathbb{S}^1_r$ factor is thought of as $[r_{-} - 2 \epsilon, r_{+} + 2 \epsilon]$ with the endpoints identified (and the metric extended periodically to it) \footnote{This is done to avoid dealing with the boundaries $r_{\pm} \mp 2\epsilon$ which would require $\Phi_u$ to preserve boundary conditions there. It does not affect the solution within $(r_{-}, r_{+})$ since $r = r_{\pm}$ are event horizons.} while $t_{\ast}$ is a suitable modification of $t$ in order to extend $g$ beyond the event horizons $r_{\pm}$, as described in \cite{petersen2023waveequationskerrdesitter}. Then, the dynamic of the null geodesics is as needed for theorem \ref{FredholmComplexPotential} to apply. We summarize the important points here, but refer to \cite{vasy2011microlocal} or \cite{KerrStab} for more details and to figure \ref{fig:KerrDynamics} for an illustration of it:

\begin{itemize}
    \item The characteristic set $\Sigma$ of non-zero null covectors $g^{-1}(\xi, \xi) = 0$ decomposes into two connected component $\Sigma^{\pm} = \{ g^{-1}(\xi, \xi) = 0, \pm g^{-1}(\xi, dt_{\ast}) > 0 \}$ corresponding to the two time orientations of $M$. 

    \item Within $\Sigma^{\pm}$, the conormal bundle $\mathcal{L}^{\pm} = \tau^{-1}(0) \cap \Sigma^{\pm} \cap (N^{\ast}\{ r = r_{-}\} \cup N^{\ast}\{ r = r_{+}\})$ of the horizons are invariant submanifolds. Furthermore, they are saddle radial sets, in the sense that they are source(+)/sink(-) within $S^{\ast}_{b}(\Bar{M})|_{\partial \Bar{M}}$ and saddle points if we include the $\tau$ behavior, i.e. $\mathcal{L}^{\pm}$ are sink(+)/source(-) in the $\tau$ direction.

    \item There is another invariant manifold, the \textit{trapped set} $\Gamma \subset \Sigma$, which is smooth of codimension 2 in $\Sigma$, symplectic and $r$-normally hyperbolic for every $r$. It is given by

    \begin{equation}
        \Gamma = \{ H_p(r) = H_p^2(r) = 0 \} \cap \tau^{-1}(0) \cap \Sigma
    \end{equation}
    
    In the case of static black holes, i.e. $a = 0$, it is the cotangent bundle of the photon sphere $r = 3M$ while for rotating black holes, it is more like a deformed sphere in the sense that the radius depends on the angular velocity: $r = r_{\xi_t, \xi_{\phi}}$. The unstable and stable manifold $\Gamma_{u/s}$ are smooth with explicit defining functions. This is discussed in great lengths in [Section 3.2 \cite{Dyatlov_2015}] and was proved for the full subextremal range in [Section 3 \cite{petersen2023waveequationskerrdesitter}]. For reference, here are the defining functions (within $\Sigma$) in Boyer–Lindquist coordinates 

    \begin{equation}
    \begin{split}
        \phi_{u/s} = \xi_r \mp sgn(r - r_{\xi_t, \xi_{\phi}}) \sqrt{\frac{F(r) - F(r_{\xi_t, \xi_{\phi}})}{\mu}}
    \end{split}
    \end{equation}

    where $r_{\xi_t, \xi_{\phi}}$ is the unique solution of the $(\xi_t, \xi_{\phi})$ dependent polynomial equation

    \begin{equation}
        F'(r_{\xi_t, \xi_{\phi}}) = 0 \quad F(r) = \frac{1}{\mu} ((r^2 + a^2)\xi_t + a \xi_{\phi})^2
    \end{equation}
    
    and $(-)$ is for $(u)$ and $(+)$ for $(s)$. As we discuss next, the (homogeneous of degree 0) geodesic flow crosses the event horizon $r_{\pm}$ in finite time along $\Gamma_{u/s}$ hence these submanifolds extend to global submanifold on $S^{\ast}\Bar{M}_b|_{\partial \Bar{M}}$, i.e. past the event horizons. Furthermore, they do not intersect $\{ \sigma = 0 \}$ hence the flow is non-trapping there. 
    
    \item All null geodesic in $\Sigma^{\pm}$ flows backward(+)/forward(-) to either $\mathcal{L}^{\pm}$, $\Gamma^{\pm}$ (in which case it is in $\Gamma_{u/s}$) or $\tau = 1$. In the latter case, it crosses $\tau = 1$ in finite time. In the forward(+)/backward(-) direction, they either start at $\mathcal{L}^{\pm}$, cross the event horizons in finite time, i.e. leave $[r_{-} - \epsilon, r_{+} + \epsilon]$, or tend to $\Gamma^{\pm}$ (in which case it is in $\Gamma_{s/u}$).
    
\end{itemize}

\vspace{0.2in}

As explained in the remarks after our definition of global NHT \ref{GlobalNHT}, we will have to make a choice of direction of propagation in each of the connected components $\Sigma^{\pm}$. Since we will be solving the initial value problem, we wish to propagate forward in time. This ammounts to propagating estimates forward along $H_p$ for future-directed covectors $\Sigma^{+}$ and backward along $H_p$ for past-directed covectors $\Sigma^{-}$. Equivalently, one can think of our dynamical system being

\begin{equation}
    \Tilde{H}_p = \frac{H_p}{g^{-1}(\xi, dt_{\ast})}
\end{equation}
Now let us start by discussing how to obtain unconditional control beyond the horizons $r_{\pm}$ and the intitial hypersurface $\tau = 1$. At the boundary $\tau = 0$ and beyond the horizons, we will apply complex absorption to make the operator elliptic there. This can be done with any $q$ a homogeneous of order 2 function on $S^{\ast}_{b} \Bar{M}|_{\partial \Bar{M}}$ which satisfies

\begin{align}
    & \pm q \geq 0 \text{ on } \Sigma^{\pm} \\
    & \pm q > 0 \text{ on } ([r_{-} - 2\epsilon, r_{-} - \epsilon) \cup (r_{+} + 2\epsilon, r_{+} + \epsilon]) \cap \Sigma^{\pm} 
\end{align}
An example of such a function is

\begin{equation}
    q = \psi(p) \eta(r) g^{-1}(\xi, dt_{\ast}) \sqrt{\xi_{t_{\ast}}^2 + \xi_r^2 + \xi_{\phi}^2 + \xi_{\theta}^2 }
\end{equation}
where $\eta(r)$ and $\psi(p)$ are positive cutoffs in $([r_{-} - 2\epsilon, r_{-} - \epsilon) \cup (r_{+} + \epsilon, r_{+} + 2\epsilon])$ and $|p| < \epsilon$ respectively. $\Tilde{Q}$ then denotes a quantization of $q$ and we set $Q = \chi(\tau) \Tilde{Q} \chi(\tau)$ for some cutoff $\chi$ equal to 1 near $\tau = 0$ and decreasing to 0 towards $\tau = 1$. As for the initial hypersurface $\tau = 1$, we will impose vanishing boundary conditions there in order to solve the initial value problem. We will thus be working on the spaces
\begin{equation}
    \Dot{H}^{s, l}_b = \{ v \in H^{s, l}_b, supp(v) \subset \{ \tau \geq 1 \} \}
\end{equation}
By applying cutoffs $\chi(\tau)$ on both sides of our operators $\Phi_u$ and $\Tilde{Q}$, we make sure that they preserve this boundary condition hence act on these spaces continuously. Thus, the control region $\mathcal{Q}$ will be
\begin{equation}
    \mathcal{Q}_{P} = \{ \tau > 1 \} \quad \mathcal{Q}_{P^{\dagger}} = ([r_{-} - 2\epsilon, r_{-} - \epsilon) \cup (r_{+} + 2\epsilon, r_{+} + \epsilon])
\end{equation}
 As for the radial sets, we only have saddle points $\mathcal{L}$. As explained above, these are
 \begin{equation}
     \mathcal{L}_{P} = \mathcal{L}_{P^{\dagger}} = \mathcal{L}^{-} \cup \mathcal{L}^{+}
 \end{equation}
Note that they satisfy the assumptions of definition \ref{AbstractSaddleRadialSets} by the dynamics of geodesics discussed above. This puts us in the setting of theorem \ref{FredholmComplexPotential} hence we can apply it. The operators $\Phi_u$ and the improvement order $\alpha$ for $\Phi_u v$ are constructed as described in Appendix \ref{AppendixB}, using the function $r$ to construct a monotone function $F$. Roughly speaking, this amounts to working with a finite time propagation of $\Gamma_u$ rather than the full backward trapped set. Then, the spaces $\mathcal{\Dot{X}}^{s, l}$ and $\mathcal{\Dot{Y}}^{s-m+1, l}$ are then defined analogously to above \ref{DefinitionSpaces}, except we replace $\Phi_u$ by $\chi(\tau) \Phi_u \chi(\tau)$. Thus, we conclude that for $s, l$ satisfying

\begin{align}
    & s > (5/2 + \frac{2(1+\lambda)(r^2 + a^2)l}{|\mu'(r)|} \mp \frac{(r^2 + a^2 \cos^2(\theta))p_1}{|\mu'(r)|})|_{\mathcal{L}^{\pm}} \\
    & l < \frac{\min(\mu, \nu/2) - p_1}{f}|_{\Gamma}
\end{align}
and for $l = -Im(\sigma)$ not having a resonance of the indicial family $\widehat{N(P-iQ)}(\sigma)$,

\begin{equation}
    P - iQ : \mathcal{\Dot{X}}^{s, l} \rightarrow \mathcal{\Dot{Y}}^{s-m+1, l}
\end{equation}
 is Fredholm. Note that one can take $s$ to be constant here by taking it to be sufficiently large. In the case of scalar wave equations, one can assume $p_1 = 0$ by choosing the inner-product induced by $g$ itself hence we get a Fredholm estimate for slightly decaying spaces. For the wave operator between vector bundles, $p_1$ is in general arbitrary and depends on the inner-product chosen, but dynamical arguments can adjust the inner-product so that $p_1$ is as small as desired (currently only written for small angular momentum $a$ in [Theorem 4.8 \cite{Hintz_2017}]) hence we again obtain a Fredholm theory for slightly decaying spaces. Compare this to [Theorem 5.4 \cite{KerrStab}] where the Fredholm theory obtained for $\Box$ had a loss of 2 derivatives between $\Box v$ and $v$. This loss is due to the trapping, and was avoided here by the inclusions of the defining functions $\Phi_u$ in the function spaces. 

\vspace{0.2in}

\begin{figure}[ht!]
    \centering
    \includegraphics[width=0.5\linewidth]{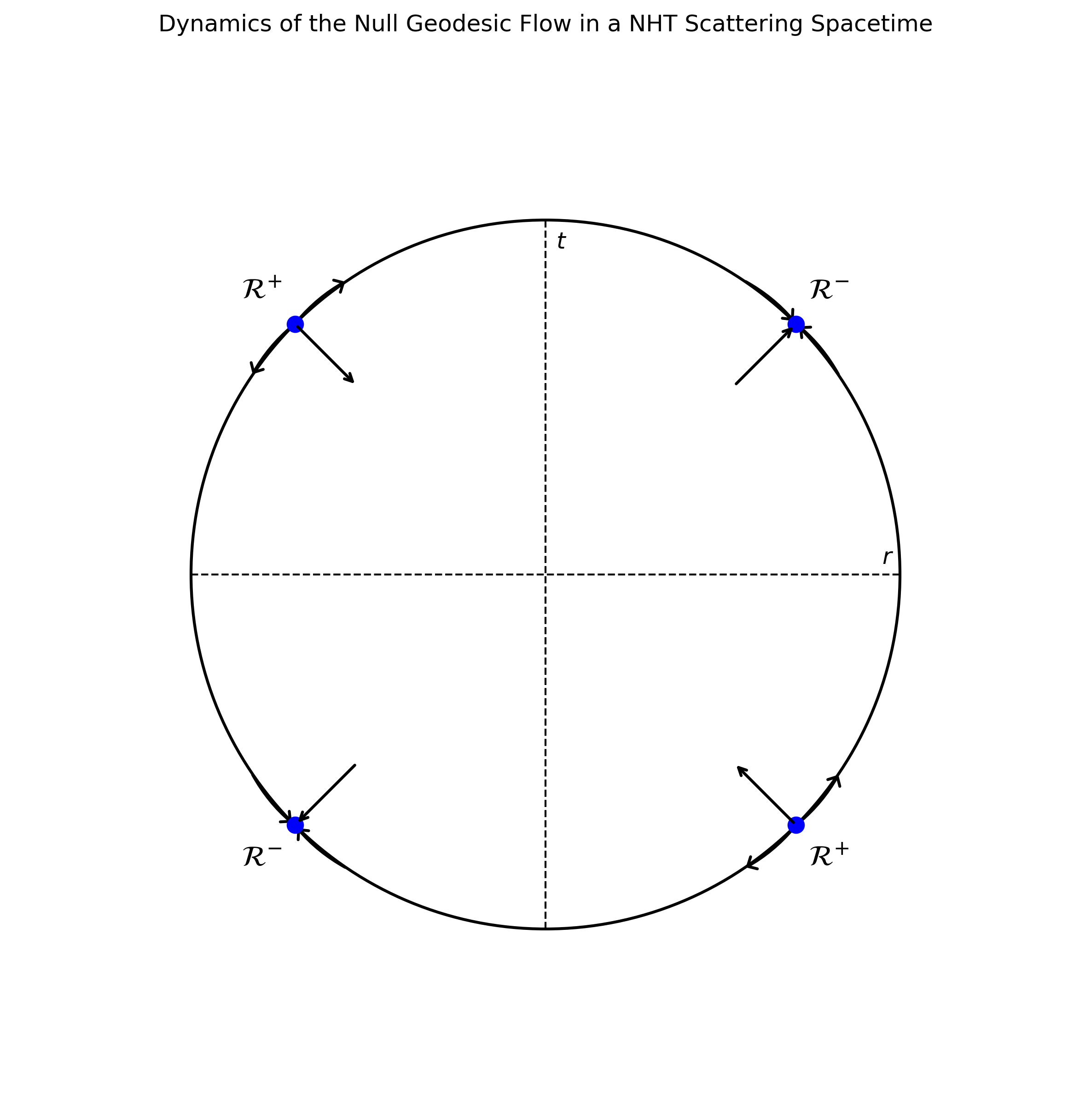}
    \caption{This is a 2D sketch of the null geodesic flow in the radial compactification of a scattering spacetime for $\eta < 0$. The null cones $r^2 = t^2$ at infinity acts as radial sets $\mathcal{R}$ and are source or sink depending on which quadrants of the $(r,t)$ plane they are in. The trapped set $\Gamma$ is not depicted here because it does not project to physical space nicely.}
    \label{fig:ScatteringExample}
\end{figure}

\textbf{Scattering Spacetime Example:} The next example we give is within the framework of scattering spacetimes. These are a generalization of the radial compactification $\tau = \sqrt{1 + r^2 + t^2}$ of Minkowski which have similar dynamical features. In particular, there is a 'light cone at infinity' which serves as radial points. We refer to \cite{baskin2014asymptotics} and [Definition 3.1 \cite{baskin2018asymptotics}] for the definition of scattering spacetimes, and to \cite{Gell_Redman_2016} for a nice application of the theory to the construction of the Feynman propagator. The advantage of our approach here is that we can now treat them when they have normally hyperbolic trapping. We give a simple example in order to illustrate that. Consider the following Lorentzian variation of the standard Euclidean metric
\begin{equation}
    ds^2 = -dt^2 + dr^2 + (1 + r^2 + t^2) d\omega^2 \quad \text{on $\mathbb{R}^2 \times S$}
\end{equation}
where $d\omega^2$ is the metric on a compact surface $S$. We will require of this closed surface $S$ that its geodesic flow $\varphi^S_t$ has subexponential differential, meaning that
\begin{equation}
    \| d\varphi^S_{t} \| \lesssim e^{\epsilon |t|}
\end{equation}
    for any $\epsilon > 0$. This is to ensure that the trapping we will obtain is $r-$normally hyperbolic for any $r$. It is satisfied by the sphere $\mathbb{S}^2$ or the Torus $\mathbb{T}^2$ for example, but not for most closed surfaces. More precisely, generically, a closed surface will contain a nontrivial hyperbolic set in which there is exponential growth, see [Theorem A \cite{GenericGeodesicFlow}], so the assumptions we are making are non-generic. Small enough exponential growth would also works, as it would corresponds to $r$-normally hyperbolic trapping for large $r$, which corresponds to working with finite but high regularity $\Gamma_{u/s}$ . We start by showing that all the null geodesics of this metric leave any compact set in $\mathbb{R}^2 \times S$. Once we have done the compactification, this will show that they all tend to the boundary, hence will allow us to restrict our analysis to there by using the propagation of singularities theorem. In these Cartesian coordinates, the equations for the geodesic flow in the $(r,t)$ plane are
\begin{align}
    &\frac{dt}{ds} = - 2 \xi_t  \\
    &\frac{dr}{ds} = 2 \xi_r \\
    &\frac{d\xi_t}{ds} = \frac{2t}{(1+r^2+t^2)^2} h(\lambda, \lambda) \\
    &\frac{d\xi_r}{ds} = \frac{2r}{(1+r^2+t^2)^2} h(\lambda, \lambda)
\end{align}
where $h(\lambda, \lambda)$ is the norm of the cotangent vectors in $T^{\ast}S$ and is conserved. If $h=0$, these are just the equations of straight lines hence the integral curve eventually leaves any compact set since $\xi_t$ and $\xi_r$ cannot be both 0 for null covectors. Now note that within the characteristic set, we have
\begin{equation}
    \xi_t^2 = \xi_r^2 + \frac{h}{1+r^2+t^2}
\end{equation}
If $h \neq 0$, then this gives
\begin{equation}
    1+r^2+t^2 \geq \frac{h}{\xi_t^2}
\end{equation}
But then, suppose $\xi_t^2$ stays bounded away from 0. Then $\frac{dt}{ds}$ does too hence $|t| \rightarrow \infty$. If it does not, the inequality above says that $1+r^2+t^2$ grows to $\infty$. Thus, in all cases, integral curves leave any compact sets, and hence we can restrict our attention to the boundary of our compactification. After radially compactifying the $\mathbb{R}^2$ factor, we obtain the following b-metric at the boundary $\tau = 0$
\begin{equation}
\begin{split}
    & \tau^2 ds^2 =  -\cos(2\theta) d\theta^2 + \cos(2\theta) \frac{d\tau^2} {\tau^2} + 2\sin(2\theta) d\theta \frac{d\tau}{\tau} + d\omega^2 \\
    & \tau = \frac{1}{\sqrt{r^2+t^2}} \\
    & \cos(\theta) = \frac{r}{\sqrt{r^2+t^2}} \quad \sin(\theta) = \frac{t}{\sqrt{r^2 + t^2}}
\end{split}
\end{equation}
We will now show that this satisfies the assumptions of theorem \ref{FredholmComplexPotential}. This will also illustrate how one checks them in practice. If we denote the dual variables by
\begin{equation}
    \sigma \frac{d\tau}{\tau} + \eta d\theta + \lambda
\end{equation}
then the inverse metric is 
\begin{equation}
    \tau^{-2} g^{-1} = \cos(2\theta)(\sigma^2 - \eta^2) + 2 \sin(2\theta) \sigma \eta + h(\lambda, \lambda)
\end{equation}
which leads to the following equations for geodesics in the $(\tau, \theta, \sigma, \eta)$ variables; the behavior in the $S$ factor is the regular geodesic flow there. Furthermore, a factor of $1/2$ has been factored out, so we are considering $\frac{1}{2} H^2_{\tau^{-2} g^{-1}}$, which corresponds to working with the wave operator $P = \frac{1}{2} \tau^{-2} \Box$.
\begin{align}
    &\frac{d\tau}{ds} = \tau (\cos(2\theta) \sigma + \sin(2\theta) \eta) \\
    &\frac{d\theta}{ds} = -\cos(2\theta) \eta + \sin(2\theta) \sigma \\
    &\frac{d\sigma}{ds} = 0 \\
    &\frac{d\eta}{ds} = - \sin(2\theta) \eta^2 + \sin(2\theta) \sigma^2 - 2 \cos(2\theta) \sigma \eta
\end{align}
To determine if the dynamics have normal hyperbolic trapping, we proceed as described in Appendix \ref{AppendixB} by finding a function $f$ such that it's second derivative $H_{g^{-1}}^2(f)$ is non-zero when $H(f) = 0$ and we are away from the trapping. The trapped set will then be $H(f)^{-1}(0) \cap H^2(f)^{-1}(0)$. To do so, just note that 
\begin{equation}
\begin{split}
    &\frac{d\theta}{ds} = -\cos(2\theta) \eta + \sin(2\theta) \sigma \\
    &\frac{d^2\theta}{ds^2} = 2\sin(2\theta) \eta \frac{d \theta}{ds} + \cos(2\theta)\sin(2\theta) (\eta^2 + \sigma^2)
\end{split}
\end{equation}
hence when $\frac{d\theta}{ds} = 0$
\begin{equation}
    \frac{d^2\theta}{ds^2} = \cos(2\theta) \sin(2\theta) (\eta^2 + \sigma^2)
\end{equation}
which is non-zero in the characteristic set away from the following invariant submanifolds
\begin{equation}
\begin{split}
    & \mathcal{R} = \{ \cos(2\theta) = 0, \sigma = 0 \} \cap (\tau^2 g)^{-1}(0) \\
    & \Gamma = \{ \sin(2\theta) = 0, \eta = 0 \} \cap (\tau^2 g)^{-1}(0)
\end{split}
\end{equation}
Note that $\mathcal{R}$ is the conormal bundle of the light cones at infinity $r^2 - t^2 = 0$. This presence of a light cone at infinity acting like a radial set is the characterizing feature of scattering spacetimes. We now show that it is a radial set, in the sense that frequencies are expanding/contracting exponentially fast. To do so, note that the derivatives of its defining functions in $S^{\ast}_b M|_{\partial M}$ are
\begin{equation}
    \frac{1}{\eta} \frac{d}{ds} \begin{pmatrix}
     \cos(2\theta) \\
    \frac{\sigma}{\eta}
    \end{pmatrix} = 
    \begin{pmatrix}
     2\sin(2\theta) & -2\sin^2(2\theta)  \\
     2 \frac{\sigma^2}{\eta^2} & \sin(2\theta) - \sin(2\theta) \frac{\sigma}{\eta}
    \end{pmatrix}
    \begin{pmatrix}
     \cos(2\theta) \\
    \frac{\sigma}{\eta}
    \end{pmatrix}
\end{equation}
At $\cos(2\theta) = 0, \sigma = 0$, this matrix has spectrum
\begin{equation}
    sgn(\sin(2\theta)) \{ 2, 1 \}
\end{equation}
hence is a source/sink according to the sign of $\sin(2\theta) \eta$. Furthermore, the frequencies have the same behavior since 
\begin{equation}
    \frac{1}{\eta} \frac{d}{ds} \frac{1}{\eta} = \sin(2\theta) \frac{1}{\eta} + O(\cos^2(2\theta), \sigma^2)
\end{equation}
hence has a non-zero weight at $\mathcal{R}$ with the same sign as $\sin(2\theta)$ too. Finally, the behavior transversal to the boundary in $\tau$ is also the same since
\begin{equation}
    \frac{1}{\eta} \frac{d}{ds} \tau = \sin(2\theta) \tau
\end{equation}
Thus, we conclude $\mathcal{R}$ is indeed a radial set as described in theorem \ref{RadialPointsB}. It is a source for $\sin(2\theta) \eta > 0$ and a sink for $\sin(2\theta) \eta < 0$. Note that the sign of $\eta$ determines which half of the characteristic set we are at $\mathcal{R}$, with opposite behavior in each half. we will denote the source and sink by
\begin{equation}
    \mathcal{R}^{\pm} = \mathcal{R} \cap \{ \pm \sin(2\theta) \eta > 0 \}
\end{equation}
We thus have unconditional control of $v$ by $Pv$ in a neighborhood of $\mathcal{R}^{+}$ on high enough regularity spaces and can propagate regularity to $\mathcal{R}^{-}$ on low regularity spaces. One issue is that the treshold regularity necessary for that cannot be achieved at both $\mathcal{R}^{+}$ and $\mathcal{R}^{-}$ if $s$ is taken to be constant. To resolve this, one must let $s$ vary so that it is higher than the required treshold at $\mathcal{R}^{+}$ but lower than the one at $\mathcal{R}^{-}$. As long as $s$ decreases along the Hamiltonian flow, we can still propagate estimates forward along the Hamiltonian flow of $H_p$, while if it increases, one must propagate backward. A simple way to do so is to let $s$ decreases below the required treshold in a neighborhood of $\mathcal{R}^{-}$. For example, if $\phi$ is a quadratic defining function for it in $S^{\ast}_b M|_{\partial M}$ so that $\Tilde{H}_p(\phi) \sim -w \phi$ with $w|_{\mathcal{R}^{-}} > 0$, one can simply pick

\begin{equation}
    s = const. - \psi(\phi)
\end{equation}
for some cutoff $\psi$ in $(0, \epsilon)$ which is equal to some large number $C$ near $\phi = 0$ and going monotonically to 0 at $\epsilon$. Note also that this does not affect our trapping estimate since we are modifying $s$ only in a neighborhood of $\mathcal{R}^{-}$, which is away from the trapped set $\Gamma$. For the reminder of this example, $s$ is assumed to vary in such a way. As for $\Gamma$, it is a normally hyperbolic trapped set. Indeed, its the cotangent bundle of $\sin(2\theta) = 0$ hence is symplectic. Furthermore, $\sigma$ is non-vanishing on it and satisfies $\frac{d \sigma}{ds} = 0$. Finally, the derivatives of its defining functions are given by
\begin{equation}
\label{NHT_Scattering}
    \frac{1}{\sigma} \frac{d}{ds} \begin{pmatrix}
     \sin(2\theta) \\
    \frac{\eta}{\sigma}
    \end{pmatrix} = 
    \begin{pmatrix}
     2\cos(2\theta) & -2\cos^2(2\theta)  \\
     1 - \frac{\eta^2}{\sigma^2} & -2\cos(2\theta)
    \end{pmatrix}
    \begin{pmatrix}
     \sin(2\theta) \\
    \frac{\eta}{\sigma}
    \end{pmatrix}
\end{equation}
which at $\sin(2\theta) = 0, \eta = 0$ has eigenvalues $\pm \sqrt{2}$. This corresponds to normally hyperbolic trapping. It is $r$-normally hyperbolic trapping for every $r$ because the flow on this trapped set $\Gamma$ is given by the geodesic flow on $S$ which has subexponential bounds on its differential $d\varphi^S_t$ by assumption. As for the behavior transversal to the boundary, we have
\begin{equation}
    \frac{1}{\sigma} \frac{d}{ds} \tau = \cos(2\theta) \tau
\end{equation}
At $\Gamma$, $\cos(2\theta) < 0$ hence this is indeed of the form $-f \tau$ for $f>0$. Note once again that the behavior of the flow is opposite for the two different connected components of the characteristic set. More precisely, the sign of $\sigma$ determines which half of the characteristic set we are in, with the direction of propagation of our estimates being forward along $H_p$ for $\sigma > 0$ and backward along $H_p$ for $\sigma < 0$. Equivalently, one can consider $\frac{1}{\sigma} H_p$ as the 'correct' dynamical system. Finally, one must check that the unstable and stable bundles  $E_{u/s}$ are trivial over $\Gamma$ in order for $\Gamma_{u/s}$ to have local defining functions. To do so, note that their dual will be spanned by the differentials of the eigenfunctions of the matrix above \ref{NHT_Scattering} hence it suffice to show that they globally exist on $\Gamma$. But at the trapped set $\sin(2\theta) = \eta = 0$ this is just a constant coefficient matrix whose eigenfunctions can be explicitely computed. Thus, $E_{u/s}$ are trivial and their dual are given by
\begin{equation}
    E^{\ast}_{u/s} = \text{span} \{ d(\sin(2\theta)) + (\text{sgn}(\cos(2\theta)) 1 \mp 1/\sqrt{2}) d(\frac{\eta}{\sigma}) \}
\end{equation}
Thus, after choosing some propagation finite time $T$ from a neighborhood of $\Gamma$, we can apply the result Proposition \ref{TrappingConclusion} of Appendix \ref{AppendixB} to construct the operators $\Phi_u$ and improvement order $\alpha$ for our Fredholm theory. In conclusion, the assumptions of \ref{FredholmComplexPotential} are satisfied, with the a priori-control region and the radial sets being
\begin{align}
    &\mathcal{Q}_{P} = \text{Neighborhood of $\mathcal{R}^{+}$} \quad &\mathcal{Q}_{P^{\dagger}} = \text{Neighborhood of $\mathcal{R}^{-}$} \\
    &\mathcal{R}_{P} = \mathcal{R}^{-} \quad &\mathcal{R}_{P^{\dagger}} = \mathcal{R}^{+}
\end{align}
This is summarized in figure \ref{fig:ScatteringExample}, which illustrates the dynamics in the $(r, t)$ plane. 
Thus, for $s, l$ satisfying 
\begin{align}
    & (2s - (m - 1)) + 2l - p_1 )|_{\mathcal{R^{+}}} > 0 \\
    & (2s - (m - 1)) + 2l + p_1 )|_{\mathcal{R^{-}}} < 0 \\
    & l < \min(\mu, 1/\sqrt{2}) - p_1|_{\Gamma}
\end{align}
and $Im(\sigma) = -l$ not containing a resonance of the indicial family $\widehat{N(P)}(\sigma)$, the wave operator $P$ and its $O(\tau^{\mu})$ perturbations are Fredholm between the spaces
\begin{equation}
    P = \frac{1}{\tau^2} \Box : \mathcal{X}^{s, l} \rightarrow \mathcal{Y}^{s-m+1, l}
\end{equation}
Again, this holds with $p_1=0$ for scalar wave equations, while for tensorial waves there will be in general a non-zero subprincipal symbol. Note in particular that no complex absorbing operator $Q$ or boundary conditions are necessary for this example hence in that respect it is more 'natural' then the Kerr-de Sitter one.

\section{Conclusion}\label{sec5}

To summarize, we have shown that one can obtain a Fredholm theory on `reasonable' spaces for non-elliptic operators with normally hyperbolic trapping. The function spaces necessary have weaker regularity at the backward trapped set $\Gamma_u$, which is implemented by localizing using defining functions for $\Gamma_u$. We then applied the setup to show that wave operators on trapping spacetimes are Fredholm. Moving forward, it would be very interesting to exhibit a more varied class of spacetimes with normally hyperbolic trapping. As far as the author knows, the main examples of physical interests known are the black-hole spacetimes. Understanding how normally hyperbolic trapping appears in a dynamical context would nicely complement our current understanding of those flows and would allow us to apply our understanding of NHT to more varied situations.

\section{Acknowledgments}

I first would like to thank my advisor Andras Vasy for suggesting this topic and his constant support. This project would have never been finished without his guidance and the many discussions we had. Furthermore, I am also grateful to Peter Hintz, Semyon Dyatlov and Jared Wunsch for preliminary comments on a draft; they were very helpful to improve the presentation of this manuscript. Finally, I would also like to thank Stanford's math department and the FRQ (`Fond de Recherche du Québec'; Doctoral Research Scholarship \#332587) for funding my PhD and allowing me to pursue my research, as well as the National Science Foundation for travel support under grant number DMS-2247004 (with PI Andras Vasy). 

\bibliographystyle{acm}
\bibliography{main}

\appendix

\section{Coisotropic Sobolev Spaces}
\label{AppendixA}

In this appendix we give a short introduction to Sobolev spaces associated with a coisotropic homogeneous manifold $L \subset T^{\ast}M$. Furthermore, we prove the propagation of singularity theorem for such spaces when the operator has an imaginary part. Overall, these spaces are useful when one expects solutions of a PDE to be more regular outside of $L$. Given a coisotropic homogeneous manifold with defining functions $\phi_j$ homogeneous of order 0 quantized to the pseudo-differential operators $\Phi_j$, we define the associated coisotropic Sobolev space of order $(s, k, \alpha)$ to be
\begin{equation}
    H^{s, k, \alpha}_L = \{ v \in H^s: \Phi_J v \in H^{s+r \alpha} \quad \forall J = (j_1, ..., j_r), r \leq k \}
\end{equation}
where $\Phi_J = \Phi_{j_1} ... \Phi_{j_l}$, equipped with the graph norm
\begin{equation*}
    \| v \|_{H^{s, k}_L} = \| v \|_{s} + \sum_{J = (j_1, ..., j_r), r \leq k} \| \Phi_J v \|_{s+r\alpha}
\end{equation*}
It will be important for our application to allow $\alpha$ to vary, hence we will let $\alpha$ be a smooth function on $S^{\ast}M$. Furthermore, we restrict to $0 \leq \alpha \leq 1$. Similary, one can define b-coisotropic spaces in the same fashion, replacing standard Sobolev spaces by b-Sobolev spaces. In words, these are functions for which we gain an $\alpha$ order of regularity each time we apply a PSDO to $v$ whose principal symbol vanish on $L$. In non-elliptic problems, these spaces appear a lot since the failure of elliptic regularity leads to solutions which are rougher inside the characteristic set. In our case, the spaces defined in section \ref{FredholmFunctionSpaces} are coisotropic spaces associated to the backward trapped set $\Gamma_u = p^{-1}(0) \cap \phi_u^{-1}(0)$:
\begin{equation}
    \mathcal{Y}^{s, \alpha} = H^{s, 1, \alpha}_{\Gamma_u} \quad \mathcal{X}^{s, \alpha} = P^{-1}(\mathcal{Y}^{s-m+1}) \subset H^{s, 1, \alpha}
\end{equation}
While strange at first, they satisfy the two most important properties we could wish for: they are acted upon by pseudo-differential operators and contain the smooth functions as a dense subspace. Thus, for most purposes, they behave like classical Sobolev spaces. This is the content of the next proposition, which is standard cf. \cite[Corollary 6.5]{melrose2011diffractionsingularitieswaveequation}.
\begin{proposition}
    Let $A \in \Psi^m(M)$ be a pseudo-differential operator of order $m$. Then 

    \begin{equation}
        A: H^{s, k, \alpha}_L \rightarrow H^{s-m, k, \alpha}_L
    \end{equation}
    continuously. Furthermore, the subspace $C^{\infty}(M)$ is dense in $H^{s, k, \alpha}_L$. The same result holds for $b$-coisotropic spaces.
\end{proposition}
\textbf{Proof:}
We first do the case $k = 1$. Note that $\Phi_j A v = A \Phi_j v + [\Phi_j, A] v$ with $[\Phi_j, A] \in \Psi^{m-1}(M)$ hence 

\begin{equation}
\begin{split}
    &\| \Phi_j A v \|_{s-m+\alpha} \leq \| A \Phi_j v \|_{s-m+\alpha} + \| [\Phi_j, A] v \|_{s-m+\alpha} \\ 
    & \qquad \leq \max(\|A\|_{H^{s+\alpha} \rightarrow H^{s-m+\alpha}}, \|[\Phi_j, A]\|_{H^s \rightarrow H^{s-m+\alpha}}) (\| \Phi_j v\|_{s+\alpha} + \| v \|_{s})
\end{split}
\end{equation}
which gives the first claim upon summing. For $k>1$, we use induction together with
\begin{equation*}
    \Phi_J Av = A \Phi_J v + \sum_{I = (i_1, ..., i_r), r < k} \Phi_I R_I
\end{equation*}
where the $R_I$ are of order $ord(A) - (k-r)$ and are given by commutators between $A$ and the factors $\Phi_j$ of $\Phi_J$. This gives the desired bound on $\Phi_J Av$ in terms of lower order coisotropic norms of $v$. As for the second claim, let $S_{\epsilon}$ be a family of smoothing operator tending strongly to the identity and uniformly bounded in Sobolev norms. Then $S_{\epsilon} v \rightarrow v$ in $H^s$ while 
\begin{equation}
    \Phi_J S_{\epsilon} v = S_{\epsilon} \Phi_J v + \sum_{I = (i_1, ..., i_r), r < k} \Phi_I R_I v
\end{equation}
and the commutators $R_I$ going to 0 in $H^{s+r\alpha}$ shows that $\Phi_J S_{\epsilon} v$ tends to $\Phi_J v$ in $H^{s+r\alpha}$ too. Thus smooth functions are dense.
$\Box$

\vspace{0.2in}

Note that asking for at most one order regularity improvement in $\Phi_j v$ is crucial to make these properties hold. This comes from the fact that commutators between pseudo-differential operators drops the order by only $1$ in general. Another important tool that is needed to work in these spaces is a propagation result for non-elliptic operators. For us, this is required in order to transport our estimate near the trapped set $\Gamma$ to the whole of its unstable manifold $\Gamma_u$. As we will need the case with a complex potential, we let $Q$ be a self-adjoint PSDO of order $m$ whose principal symbol $q$ is real, homogeneous of order $m$ and non-negative $q \geq 0$. In that case, a propagation result similar to the usual one holds, except that we lose a $1/2$ derivative in $\Phi_j v$ when propagating into $q > 0$. The case where $Q = 0$ and $\alpha = 1$ was proven in [Proposition 6.1 \cite{melrose2011diffractionsingularitieswaveequation}] in a different context. The main additional technicality is that one must be keep track of the commutators with $Q$ and $\rho^{\alpha}$ in order to handle losses of regularity.

\begin{proposition}
\label{PropagationCoisotropic}
    Let $P - i Q \in \Psi^{m}(M)$ be a PSDO of order $m$ for which the real part $p$ of its principal symbol is homogeneous of degree $m$ and non-degenerate, i.e. $dp \neq 0$ on $p^{-1}(0)$. Suppose that $H_p$ is tangent to the coisotropic submanifold $L$, that $q$ is non-negative $q \geq 0$ and that $\alpha$ is non-increasing along the flow of $H_p$: $H_p(\alpha) \leq 0$. Let $\chi, \psi, G$ be $0^{th}$ PSDO such that every points in $ell(\chi)$ reaches $ell(\psi)$ in finite time along the backward flow of $H_p$ all while remaining in $ell(G)$. Then the following inequality holds in a strong sense for all $v \in H^{-\infty}$, $N$ arbitrary:

\begin{equation}
\begin{split}
    &\sum_{J = (j_1, ..., j_r)} \| \Phi_J \chi v \|_{s+r\alpha} \lesssim \sum_{J = (j_1, ..., j_r)} \| \Phi_J G(P - iQ) v \|_{s+r\alpha-m+1} + \| \Phi_J \psi v \|_{s+r\alpha} \\
     & \qquad + \sum_{I = (i_1, ..., i_{l}), l < r} \| \Phi_{I} G_I(P - iQ) v \|_{s+r\alpha-m+1 - (r-l)/2} + \| \Phi_I \psi_I v \|_{s+r\alpha-(r-l)/2} + \|v\|_{-N}
\end{split}
\end{equation}
   where $G_I$ and $\psi_I$ are $0^{th}$ order with the same wave-front set as $G$ and $\psi$. If the propagation is done within $q^{-1}(0)$, that is $q = 0$ in $ell(G)$, then we can replace this by propagation of coisotropic regularity:
\begin{equation}
    \| \chi v \|_{H^{s, k, \alpha}_L} \lesssim \sum_{J = (j_1, ..., j_r), 0 \leq r \leq k} \| \Phi_J G_J(P - iQ) v \|_{s+r\alpha-m+1} + \| \Phi_J \psi_J v \|_{s+r\alpha}  + \|v\|_{-N}
\end{equation}
The same result holds in the b-category, that is if we replace the spaces above by b-coisotropic Sobolev spaces. Furthermore, if $q \leq 0$ and $H_p(\alpha) \geq 0$, the same result holds but now propagating backward along $H_p$ instead.
\end{proposition}
\textbf{Proof:} 
We present the proof in the case of $M$ closed, and then comments on how to modify it for a manifold with boundary. Let's start with the $k = 1$ case. For an operator $A$ with real principal symbol, one has the following commutator formula
\begin{equation}
\begin{split}
    &\left< i[P,A^{\dag}A] v, v \right> + \left< \frac{1}{i} A^{\dag} A (P - P^{\dag}) v, v \right> = \\
    & \qquad 2Im( \left< (P-iQ)v , A^{\dag} Av\right> ) + 2Re( \left< Qv, A^{\dag} Av \right> )
\end{split}
\end{equation}
The idea of the proof, as any other positive commutator estimate, is to arrange $A$ so that $[P, A^{\dag} A]$ is very large. Then using the equation above, one uses it to estimate $v$ in terms of $(P - iQ)v$. The second term on the right is the one that will impose a specific direction of propagation. It can be rewritten as 
\begin{equation}
\begin{split}
    Re( \left< Qv, A^{\dag} A v \right> ) &= Re( \left< QAv, Av \right> ) + Re( \left< [A, Q]v, Av \right> ) \\
    &= Re( \left< QAv, Av \right> ) + \frac{1}{2} \left< [A, [A, Q]]v, v \right> + \frac{1}{2} \left< \{ (A^{\dag} - A) [A, Q] - [A^{\dag} - A, Q] A \} v, v \right>
\end{split}
\end{equation}
 The first term can be estimated using the Sharp-G\r{a}rding inequality (this is where we use $q \geq 0$)
\begin{equation*}
    Re( \left< QAv, Av \right> ) \geq -c \| Av \|^2_{(m-1)/2}
\end{equation*}
where the $m-1$ denotes the $H^{m-1}$ norm, in other words we drop one order. The constant $c$ only depends on $Q$; this will be important when we wish to dominate this term later on. As for the double commutator terms and the ones involving $A - A^{\dag}$, note that it is one order lower than the commutator. Thus, we can recast the equation above as
\begin{equation}
\begin{split}
    & \left< -i [P,A^{\dagger}A] v, v \right> - \left< \frac{1}{i} A (P - P^{\dag}) v, A v \right> - c \| Av \|^2_{(m-1)/2}  + \left< Sv, v \right> \leq 2Im( \left< Av, A(P-iQ)v \right> )
\end{split}
\end{equation}
 with $S$ denoting the double commutator and the $A - A^{\dag}$ terms. Note that this imposes a sign on $H_p(a^2)$ as it must now be negative in order to obtain a bound. This corresponds to propagating the estimate \textit{forward}. By a partition of unity argument on the flowout of $ell(\chi)$, it suffices to restrict ourselves to small enough coordinate systems $(t, z)$ on $T^{\ast}M$ such that $\frac{1}{\rho^{m-1}} H_p = \partial_t$. Then, pick a smooth escape function $f(t)$ on $S^{\ast}M$ with support in $(-\epsilon, 1)$ that satisfy
\begin{equation}
    f'(t) = -b(t)^2 f(t) + e
\end{equation}
for some function $b$ bounded below and with the error $e$ supported in $(-\epsilon, \epsilon) \cup (1-\epsilon, 1)$, negative in $(1 - \epsilon, 1)$. One can take $b$ as large as we want, by considering functions of the form $e^{a^2/(t-1)}$ for large constants $a$. We will choose as commutant a quantization $A$ of
\begin{equation}
    a_j = \phi_j f(t) \eta(z) \theta \left( \frac{p^2}{\rho^{2m}} + \frac{q^2}{\rho^{2m}} \right) \rho^{(s+\alpha)-(m-1)/2} 
\end{equation}
where $\eta, \theta$ are cutoffs in $\delta$-sized ball, $\delta$ small, and respectively constrain the problem to a small neighborhood of the curve and near the characteristic set. Furthermore, we take $A = \Lambda \Phi_j \Op(\chi)$ with $\Phi_j$ a quantization of $\phi_j$, $\Op(\chi)$ of $f \eta \theta$ and $\Lambda$ of $\rho^{(s+\alpha) - (m-1)/2}$. The idea is that the derivative of $f$ will give the weight we need while the remaining terms of $H_p(a^2)$ will be errors we can control. The principal symbol of the left handside (ignoring $S$ for now) is then
\begin{equation}
\begin{split}
    &-H_p a_j^2 - 2p_1 a_j^2 - c a_j^2 \rho^{m-1} = (2b^2 + g - 2 ln(\rho) H_p(\alpha) - c) a_j^2 \rho^{m-1} - (p^2 + q^2) h - \Tilde{e} (\phi_j)^2
\end{split}
\end{equation}
where $g$ is some smooth factor (containing $-2p_1$, the coefficients of $H_p(\phi)$ and the derivative $H_p(\rho)$), the function $h$ is of order $2(s+\alpha)-2m$ and comes from the derivative term $\theta'$ (hence is supported away from $p = q = 0$) while the error $\Tilde{e}$ comes from $e$ in $f'(t)$. Since $g$ is bounded and $H_p(\alpha) \leq 0$ we can take $b$ large enough so that $b^2+g-ln(\rho)H_p(\alpha)-c = w$ is bounded below by a positive constant. Thus, the operator on the left-handside (ignoring $S$ again) becomes
\begin{equation}
\begin{split}
    &-i[P,A^{\dagger}A] - \frac{1}{i} A^{\dagger} A (P - P^{\dagger}) - c A^{\dagger} A \\
    & = -i[P, \Op(\chi) \Phi_j \Lambda^2 \Phi_j \Op(\chi)] -\frac{1}{i} \Op(\chi) \Phi_j \Lambda^2 \Phi_j \Op(\chi) (P - P^{\dagger}) - c \Op(\chi) \Phi_j \Lambda^2 \Phi_j \Op(\chi)  \\
    &= -i \{ [P, \Op(\chi)] \Phi_j \Lambda^2 \Phi_j \Op(\chi) + \Op(\chi) [P, \Phi_j] \Lambda^2 \Phi_j \Op(\chi) + \Op(\chi) \Phi_j ([P, \Lambda] \Lambda + \Lambda [P, \Lambda]) \Phi_j \Op(\chi) \\
    & \qquad + \Op(\chi) \Phi_j \Lambda^2 [P, \Phi_j] \Op(\chi)  + \Op(\chi) \Phi_j \Lambda^2 \Phi_j [P, \Op(\chi)] \} -\frac{1}{i} \Op(\chi) \Phi_j \Lambda^2 \Phi_j \Op(\chi) (P - P^{\dagger}) \\
    & \qquad \hspace{3em} - c \Op(\chi) \Phi_j \Lambda^2 \Phi_j \Op(\chi) \\
    &= \Op(\chi) \Phi_j \Lambda W \Lambda \Phi_j \Op(\chi) - \frac{1}{2} (\Tilde{E} \Phi_j \Lambda^2 \Phi_j \Op(\chi) - \Op(\chi) \Phi_j \Lambda^2 \Phi_j \Tilde{E}) \\
    & \qquad - \frac{1}{2}((P^{\dagger} + iQ)H \Phi_j \Lambda^2 \Phi_j \Op(\chi) (P-iQ) + (P^{\dagger} + iQ^{\dagger}) \Op(\chi) \Phi_j \Lambda^2 \Phi_j H (P-iQ))\\
    & \qquad \hspace{4em} - R_1 \Phi_j \Lambda^2 \Phi_j \Op(\chi) - \Op(\chi) R_2 \Lambda^2 \Phi_j \Op(\chi) - \chi \Phi_j \Lambda^2 R_3 \Op(\chi) - \Op(\chi) \Phi_j \Lambda^2 \Phi_j R_4 \\
    & \quad \hspace{8em} - \Op(\chi) \Phi_j (R_{\epsilon} \Lambda + \Lambda R_{\epsilon}) \Phi_j \Op(\chi)
\end{split}
\end{equation}
where the operators $W$, $\Tilde{E}$, $H$ quantize $w$, $\Tilde{e}$ and $h$ while the $R's$ denotes (different) reminder terms of order $m-2$ while $R_{\epsilon}$ denotes a reminder term of order $s+\alpha + (m-3)/2 + \epsilon$, $\epsilon > 0$ arbitrary, coming from commutators with $\Lambda$. The key thing is that the latter comes with two factors of $\Phi_j$ hence the $\epsilon$ loss of regularity doesn not pose problems. As for $S$, the double commutator $[A_j, [A_j, Q]]$ has the wrong sign, as its principal symbol is  
\begin{equation}
    \sigma([A_j, [A_j, Q]]) = -\sum_j H_{a_j}^2(q)
\end{equation}
In applications, $q$ will be a locally convex function hence this term is negative in $q > 0$ and crucially, will have support intersecting $\Gamma_u$. Thus, $S$ must be considered as a reminder error too, of the form $Sv = S_1 + S_2 \Phi_j \Op(\chi)$ where $S_1$ is of order $2(s + \alpha) - 1$ and $S_2$ of order $2(s + \alpha) - 1 + \epsilon$. In particular, the $\epsilon$ loss is also coming with a factor of $\Phi_j$ here, due to it only appearing when we commute $\Lambda$ with $Q$. The term $\left< Sv, v \right>$ is the reason that we get a weaker result when propagating into the region $q > 0$. We then quantize the operators above to obtain the following relation:
\begin{equation}
\begin{split}
    &\left<W \Lambda \Phi_j \Op(\chi) v, \Lambda \Phi_j \Op)(\chi) v \right> \\
     & \qquad \leq 2Im( \left< \Lambda \Phi_j \Op(\chi) v, \Lambda \Phi_j \Op(\chi) (P-iQ)v \right> \\
     & \qquad \hspace{1em} + Re(\left< \Lambda \Phi_j H (P - iQ) v, \Lambda \Phi_j \Op(\chi) (P-iQ) v \right>)
      + Re(\left< \Lambda \Phi_j \Tilde{E} v, \Lambda \Phi_j \Op(\chi) v \right>) \\
      & \qquad \hspace{2em} + |\left< Sv, v \right>| + \left< \Lambda \Phi_j R_1^{\dagger} v, \Lambda \Phi_j \Op(\chi) v \right> + \left< \Lambda R_2^{\dagger} \Op(\chi) v, \Lambda \Phi_j \Op(\chi) v \right> \\
      & \qquad \hspace{3em} + \left< \Lambda \Phi_j \Op(\chi) v, \Lambda R_3 \Op(\chi) v \right> + \left< \Lambda \Phi_j \Op(\chi) v, \Lambda \Phi_j R_4 v \right> \\
      & \qquad \hspace{4em} + \left< R_{\epsilon}^{\dagger} \Phi_j \Op(\chi) v, \Lambda \Phi_j \Op(\chi) v \right> + \left< \Lambda \Phi_j \Op(\chi) v, R_{\epsilon} \Phi_j \Op(\chi) v \right>
\end{split}
\end{equation}
Summing over $j$, using G\r{a}rding's inequality on the $W$ term, applying Cauchy-Schwartz on the pairings and pulling the $\Phi_j \Op(\chi) v$ terms on the left hand-side gives the following estimate
\begin{equation}
\begin{split}
\label{ErrorDoubleCommutator}
    \sum_j \| \Phi_j \Op(\chi) v \|_{s+\alpha}^2 &\lesssim
    \sum_j \| \Phi_j \Op(\psi) v \|_{s+\alpha}^2 + \| \Phi_j G (P-iQ) v \|_{s+\alpha-m+1}^2 \\
    &\qquad + \| \Tilde{G} (P-iQ) v \|_{s+\alpha-m}^2 + \| Ev \|_{s+\alpha-1} + \| Ev \|_{s+\alpha-1/2} + \| v \|_{-N}
\end{split}
\end{equation}
where $\psi$ quantize the positive part of $\Tilde{e}$, $G$ is a cutoff in a slightly larger neighborhood controlling $K$ and $\chi$, and $\Tilde{G}$ and $E$ controls various commutators/reminder errors but are (microlocally) supported in the same neighborhood. To improve the error terms $Ev$, we will apply a second commutator argument. This time it is just the standard propagation of singularity with complex potential. More precisely, we consider a commutant of the form
\begin{equation}
    a^2 = f^2(t) \eta(z)^2 \theta^2(\frac{p^2 + q^2}{\rho^{2m}}) \rho^{2(s+\alpha-1/2)-(m-1)} 
\end{equation}
and proceed in the exact same fashion (on a slightly enlarged neighborhood now) to obtain
\begin{equation}
    \| E v \|_{s+\alpha-1/2}^2 \lesssim \| \Op(\Tilde{\psi}) v \|_{s+\alpha-1/2}^2 + \| \Tilde{G} (P-iQ) v \|_{s+\alpha-1/2-m+1}^2 + \| v \|_{-N}^2
\end{equation}
again with the all localizers supported in a slightly enlarged neighborhood. Combining both results, one obtains
\begin{equation}
\begin{split}
    & \sum_j \| \Phi_j \Op(\chi) v \|_{s+\alpha}^2 \lesssim
    \sum_j \| \Phi_j \Op(\psi) v \|_{s+\alpha}^2 + \| \Phi_j G (P-iQ) v \|_{s+\alpha-m+1}^2 \\
     & \qquad + \| \Op(\Tilde{\psi}) v \|_{s+\alpha-1/2}^2 + \| \Tilde{G} (P-iQ) v \|_{s+\alpha-m+1/2}^2 + \| v \|_{-N}^2
\end{split}
\end{equation}
which is the desired result. Note that if we propagate only within $q^{-1}(0)$, that is $q = 0$ on the support of $\chi$, then the error term $\| Ev \|_{s+\alpha-1/2}$ due to the double commutator term drops out in equation \ref{ErrorDoubleCommutator} hence we instead only need to propagate the $s+\alpha - 1$ norm of $Ev$, giving a lossless result in coisotropic regularity, which corresponds to the second part of the proposition. To get the strong inequality, we must be able to justify the pairings and integration by parts done in the proof when $v \in H^{-N}$. This is done by adding a regulazing operator
\begin{equation*}
    S_{\epsilon} = \frac{1}{(1+ \epsilon \rho)^{2M}}
\end{equation*}
to the commutant and taking the $\epsilon \rightarrow 0$ limit. As this is rather standard, we refer to Hintz's notes on microlocal analysis [Section 8.4 \cite{PeterNote}] for the details, which are unaffected by the $\phi_j$ factor in $a_j$.
Now, for the case $k > 1$, we need to propagate higher regularity for $\phi_J v$, $|J| = r \leq k$. To do so we take a commutator of the form
\begin{equation}
    a_J^2 = (\phi_J)^2 f^{2}(t) \eta(z)^{2} \theta^{2} \left( \frac{p^2}{\rho^{2m}} + \frac{q^2}{\rho^{2r}} \right) \rho^{2(s+r\alpha)-(m-1)} 
\end{equation}
and proceed in the same fashion, expect that now the reminders we obtain are all of the form (up to reflecting around $\Lambda$)
\begin{equation}
\begin{split}
    &\Op(\chi) \Phi_{j_r} \Phi_{j_{r-1}} ... R ... \Phi_{j_2} \Phi_{j_1} \Lambda^2 \Phi_J \Op(\chi) \\
    & R (\Phi_{J})^{\dagger}\Lambda^2 \Phi_J \Op(\chi) \\
    & \Op(\chi) (\Phi_J)^{\dagger} R_{\epsilon} \Lambda \Phi_J \Op(\chi) 
\end{split}
\end{equation}
hence have only 1 less factors of $\Phi$ at most and thus can be estimated by lower coisotropic norms of $v$. More formally after commuting the reminder errors $R$ around, we end up with
\begin{equation}
\begin{split}
    &\sum_J \| \Phi_J \Op(\chi) v \|_{s+r\alpha}^2 \lesssim
    \sum_J \| \Phi_J \Op(\psi) v \|_{s+r\alpha}^2 + \| \Phi_J G (P-iQ) v \|_{s+r\alpha-m+1}^2 \\
     & \qquad + \sum_{I = (i_1, ..., i_{l}), l < r}  \| \Phi_I G_I (P-iQ) v \|_{s+r\alpha-m+1-(r-l)}^2 + \| E'_I \Phi_{I} v \|_{s+r\alpha - (r-l)}\\
     & \qquad \hspace{2em} + \| E_I \Phi_{I} v \|_{s+r\alpha - (r-l)/2} + \|v\|_{-N}
\end{split}
\end{equation}
 Where again the terms of order $s+r\alpha - (r-l)/2$ are due to the $S \sim [A, [A, Q]]$ term and drop out if $q = 0$ in our region of propagation. Using induction for the lower coisotropic norms then gives the result. Once again, if we stay away from $q \neq 0$, then the error of order $s+r\alpha-1/2$ drops out hence we get a lossless result. For $b$-manifolds, the proof is exactly the same except this time we have an extra factor $\tau^{-l}$ in the commutant, whose $H_p$ derivative $\sim -l H_p(\tau)$ can be controlled by increasing $|f'(t)|$.
$\Box$

\vspace{0.2in}

We conclude our short survey of coisotropic function spaces with one last useful thing to know about them, which is their dual. It is given by 
\begin{equation}
    (H^{s, k}_L)^{\ast}  = H^{-s} + \sum_{J = (j_1, ..., j_r), r \leq k} \Phi_J H^{-s-r\alpha}
\end{equation}
equipped with the quotient norm
\begin{equation}
    \| v \|_{(H^{s, k}_L)^{\ast}}^2 = inf_{v = v_0 + \sum_J \Phi_J v_J} \{ \|v_0\|_{-s}^2 + \sum_{J = (j_1, ..., j_r), r \leq k} \|v_J\|_{-s-r\alpha}^2 \}
\end{equation}
Note in particular that it is not symmetric in $L$, in other words we can't describe $H^{s, k}_L$ in terms of another coisotropic space. This is important for us since it means that the estimate we obtain for $P^{\dag}$ on coisotropic spaces associated to $\Gamma_s$ cannot be used to show a Fredholm estimate for it directly. To remedy this, we shall use a weaker estimate for it paired with propagation of singularity to show the existence of regular solutions to $P u = f$, hence showing that $P$ is Fredholm. 

\section{Results on Trapping Dynamics}
\label{AppendixB}
In this appendix, we describe how to check whether a given dynamical system has global normally hyperbolic trapping as described in definition \ref{GlobalNHT}. The approach is known and is nicely discussed by Dyatlov in the case of the Kerr-de Sitter metric in [Section 3.2 \cite{Dyatlov_2015}]. Here we describe the general case in order to justify that our scattering spacetime example indeed satisfies the assumptions globally. Consider a compact manifold with boundary $W$ with a vector field $H$ on it. We wish to show that either an integral curve of $H$ escapes at the boundary in finite time or it tends to some invariable submanifold $\Omega$. In pratice, $\Omega$ will be the union of the trapped sets and radial sets of the system:
\begin{equation}
    \Omega = \mathcal{R} \cup \mathcal{L} \cup \Gamma
\end{equation}
Local computations at $\Omega$ must be done to distinguish between $\mathcal{R}$, $\mathcal{L}$ and $\Gamma$. The boundary will be the a priori control region where we have estimates, say from boundary conditions or complex absorption. The simplest way to show this dynamical property of $\Omega$ is through an $\textit{escape function}$. More precisely, suppose we can find a function $f$ so that 
\begin{equation}
    H(f)(x) = 0 \implies H^2(f)(x) \geq 0 \text{ with $H(f) = H^2(f) = 0$ iff $x \in \Omega$}
\end{equation}
Then, from there we can construct an escape function $F$ defined by
\begin{equation}
    F = H(f) \exp(f) \quad H(F) = (H^2(f) + H(f)^2) \exp(f)
\end{equation}
In particular its derivative is bounded below away from $\Gamma$. Furthermore,  If $dH(f)$ are $dH^2(f)$ linearly independent, then
\begin{equation}
    \Omega = H(F)^{-1}(0) = H(f)^{-1}(0) \cap H^2(f)^{-1}(0)
\end{equation}
is a smooth codimension $2$ manifold (codimension $2k$ if $f$ is a $k$-vector). These properties is then sufficient to prove that $H$ escapes through the boundary or tends to $\Omega$: 
\begin{proposition}
    Let a function $f$ exists as above and $\gamma_t$ an integral curve of $H$. Then either $\gamma_t$ tends to $\Omega$ or reaches the boundary in finite time (in either time direction). 
\end{proposition}
\textbf{Proof:} We consider the case $t \rightarrow +\infty$ without loss of generality. Now suppose $\gamma_t$ does not leave the manifold in finite time. Then $\gamma_t$ exists for all time $t > 0$. Furthermore, if it does not tend to $\Omega$, then we can find a subsequence $t_j \rightarrow \infty$ such that $\gamma_{t_j}$ stays away from a neighborhood of $\Omega$. Furthermore, since $\gamma_{t_j}$ stays away from $H(F) = 0$ we have
\begin{equation}
    lim_{t_j \rightarrow \infty} \inf H(F)({\gamma_{t_j}}) > c > 0
\end{equation}
where $c$ is some constant. This is a contradiction, since by compactness $F$ is bounded on $M$, but has its $H$ derivative bounded below for infinite time. $\Box$

\vspace{0.2in}

This reduces the problem of checking trapping dynamics to finding $f$ and understanding the zero set $H(F)$. We now restrict our attention to the trapped set $\Gamma$, which here is the part of $\Omega$ that exhibits normally hyperbolic trapping. We can check for NHT by computing the eigenvalues of $H$ in the space of functions vanishing non-degenerately at $\Gamma$, denoted by $\mathcal{I}$. Indeed, $H$ acts on
\begin{equation*}
        H: \mathcal{I} / \mathcal{I}^2 \rightarrow \mathcal{I} / \mathcal{I}^2
\end{equation*}
by its tangency on $\Gamma$. This is finitely generated and the differential of the positive/negative eigenfunctions will generate the (dual of the) stable and unstable bundle $E_s / E_u$. Thus, we are looking for eigenvalues of opposite sign. Note that in terms of $H(f)$, $H^2(f)$ the action of $H$ is lower-triangular
\begin{equation}
    H \begin{pmatrix}
    H(f)\\
    H^2(f)
    \end{pmatrix} = 
    \begin{pmatrix}
    0 & 1 \\
    a & b &
    \end{pmatrix}
    \begin{pmatrix}
    H(f)\\
    H^2(f)
    \end{pmatrix} + O(H(f)^2, H^2(f)^2)
\end{equation}
hence if $a, b$ are constants over $\Gamma$ and $a < 0$, then NHT holds since the eigenvalues have opposite signs. As for whether $\Gamma$ is symplectic, this can be checked by verifying whether these eigenfunctions have non-trivial symplectic pairings; if they do than it is symplectic. It remain to discuss how to show that the unstable and stable manifolds exist and have defining functions. The easiest way to do so once we have expressed $\Gamma$ as a zero-set is to show that $r$-normal hyperbolicity holds for every $r$. This would give the existence of $\Gamma_{u/s}$, and we must then check the triviality of their normal bundles. To do so, we recall how $\Gamma_{u/s}$ are constructed from $r$-normal hyperbolicity. This is a simplification of [Theorem 4.1 \cite{InvariantManifolds}], which we also refer to for the proof. 
\begin{theorem} \label{UnstableStableManifoldTheorem} [Theorem 4.1 \cite{InvariantManifolds}] Suppose $\varphi_t$ is $r$-normally hyperbolic, that is
\begin{equation}
    r \beta < \nu 
\end{equation}
where $\beta$ is the maximal expansion rate along $\Gamma$ and $\nu = \min(\nu_s, \nu_u)$. Then there exists (unique) $C^r$ invariant submanifolds $\Gamma_{u/s}$ in a neighborhood of $\Gamma$ such that $T|_{\Gamma}\Gamma_{u/s} = T\Gamma \bigoplus E_{u/s}$. Furthermore, $\Gamma_{u/s}$ is a graph over
\begin{equation}
    \exp \{ v \in E_{u/s}; \| v \| < \epsilon \}
\end{equation}
where $\exp$ is the exponential map associated to the metric $g$, and it satisfies
\begin{equation}
    \varphi_{\mp t}(\Gamma_{u/s}) \subset \Gamma_{u/s}
\end{equation}
\end{theorem}
In particular, we see that if both $E_{u/s}$ are trivial over $\Gamma$ and $\Gamma$ has a trivial normal bundle then so does $\Gamma_{u/s}$ in a neighborhood of $\Gamma$. Indeed, by the last part of the theorem above,
 we see that 
 \begin{equation}
     \Gamma_u \times \mathbb{R}^{dim E_s} \simeq \Gamma \times \mathbb{R}^{dim E_u} \times \mathbb{R}^{dim E_s} \hookrightarrow W
 \end{equation}
where the last embedding is due to $\Gamma$ having a trivial normal bundle. Thus to check for the existence of local defining functions reduces to figuring out if $E_{u/s}$ are trivial, which can be easily done by exhibiting global sections for it. More precisely, the eigenfunctions of $H$ in the conormal bundle, associated to positive (u) and neative (s) eigenvalues span the duals $E^{\ast}_{u/s}$, hence their global existence over $\Gamma$ suffices to show the triviality of $E_{u/s}$. Within the framework above, this ammounts to diagonalizing $H$ in the space of functions vanishing non-degenerately at $\Gamma$ (modulo quadratic vanishing), which is an explicit computation if the coefficients are constants over $\Gamma$. 
We must then globalize $\Gamma_u$. This is done by flowing out along $H$ hence defining 
\begin{equation}
    \Gamma_u^T = \bigcup_{0 \leq t < T} \varphi_t(\Gamma_u \cap \mathcal{O}) \quad \Gamma_u = \overline{\bigcup_{0 \leq t < \infty} \varphi_t(\Gamma_u \cap \mathcal{O})}
\end{equation}
where $\mathcal{O}$ is the neighborhood of $\Gamma$ in which we did our construction and $T \geq 0$ is the propagation time. The main issue now is to show that $\Gamma_u^T$ is still embedded with trivial normal bundle. For an infinite time propagation, this would require using the monotone function $F$ to guarantee that $\Gamma_u$ does not 'wind back' upon itself, but since for our Fredholm theory we are mainly interested in $\Gamma_u^{T}$ we will not need that. Rather, we have the following result:
\begin{lemma}

  Let $N \subset \mathcal{O}$ be a submanifold of a neighborhood $\mathcal{O}$ of $\Gamma$. Define the flowout $N^T$ of $N$ by some time $T \in [0, \infty)$

    \begin{equation}
        N^T = \bigcup_{0 \leq t < T} \varphi_t(N)
    \end{equation}

    Suppose $N$ satisfies

    \begin{equation}
        \varphi_{\mp t}(N) \subset N
    \end{equation}

    for all $t > 0$. Then $N^T$ is embedded. Furthermore, if $N$ is a level set in $\mathcal{O}$, then so is $N^T$ for $T < \infty$ in $\varphi_T(\mathcal{O})$. 
    
\end{lemma}
\textbf{Proof:} The first part follows from $N^T = \varphi_{\pm T}(N)$ and $\varphi_T$ being a diffeomorphism, while for the second part one can transport any defining function of $N$ by $\varphi_T$ to a defining function for $N^T$. $\Box$
\vspace{0.2in}

Finally, we conclude by discussing how one obtains monotone function along $H$ on $\Gamma_u$ in this setup. This is required for our theory to work because the improvement order $\alpha$ of $\Phi_u v$ must be decreasing along $\Gamma_u$. To do so, we simply make use of the monotone function $F$. Indeed, we can define $\alpha$ globally on $M$ as
\begin{equation}
    \alpha = h(F)
\end{equation}
where $h$ is a decreasing cutoff function equal to $1$ in $(-\infty, a-\epsilon)$ and 0 in $[a, \infty)$. Then $H(\alpha) \leq 0$ and 
\begin{equation}
    \alpha(\varphi_{\tau}(x)) = 0 \quad \forall x \in \partial(\overline{\Gamma_u \cap \mathcal{O}}), \tau = a / \inf_{\overline{\Gamma_u} \cap \mathcal{O}^c}(H(F)) < \infty
\end{equation}
 Thus, after some finite time $\tau$, the improvement order in regularity of $\Phi_u v$ relative to $v$ becomes 0, which is exactly what we wanted. Note also that it is constant near $F = 0$, hence on a neighborhood around $\Omega$. Thus, it does not affect the trapping or radial estimates that we use in our Fredholm theory. We summarize the result of this appendix in the following proposition.
\begin{proposition}
\label{TrappingConclusion}
    Let $f$ be a smooth function such that
    \begin{equation}
    H(f)(x) = 0 \implies H^2(f)(x) \geq 0 \text{ with $H(f)(x) = H^2(f)(x) = 0$ iff $x \in \Omega$}
    \end{equation}
    Let $\Gamma \subset \Omega = H(f)^{-1}(0) \cap H^2(f)^{-1}(0)$ be a $r$-normally hyperbolic trapped set for every $r$ such that the unstable and stable bundles $E_{u/s}$ are trivial over $\Gamma$. Let $\mathcal{O}$ be an open neighborhood of $\Gamma$ in which the unstable and stable manifolds $\Gamma_{u/s}$ exists locally, with defining functions $\phi^{\mathcal{O}}_{u/s}$. Then for any finite $T > 0$ there exists embedded submanifolds $\Gamma_{u/s}^T$ obtained by flowing $\Gamma_{u/s} \cap \mathcal{O}$ along the flow $\varphi_T$ of $H$
    \begin{equation}
        \Gamma_{u/s}^T = \varphi_{\pm T}(\Gamma_{u/s} \cap \mathcal{O}) \quad 
    \end{equation}
    with defining functions $\phi^T_{u/s} = \phi^{\mathcal{O}}_{u/s} \circ \varphi_{\mp T}$ in $\varphi_{\pm T} (\mathcal{O})$. Furthermore, there exists a global function $\alpha$ equal to $1$ in $\mathcal{O}$ which is decreasing along the flow of $H$, i.e. $H(\alpha) \leq 0$, and satisfying
    \begin{equation}
    \begin{split}
        \alpha(\varphi_{\tau}(x)) &= 0 \quad \forall x \in \partial(\overline{\Gamma_u \cap \mathcal{O}}) \\
        \tau &= \frac{(\sup_{\mathcal{\Bar{O}}}(F) + \epsilon)}{\inf_{\overline{\Gamma_u} \cap \mathcal{O}^c}(H(F))} < \infty
    \end{split}
    \end{equation}
    where $F$ is given by $F = H(f) \exp(f)$ and $\epsilon > 0$ is arbitrary.
\end{proposition}
We can now apply the results of this appendix to the situations we find in this paper. Given a differential operator $P$ of degree $m$ on a compact manifold with boundary $M$, consider the Hamiltonian vector field $H_p$ of its principal symbol. It is homogeneous of degree $m-1$, hence if $\rho$ is homogeneous of degree $1$ and non-vanishing then
\begin{equation}
    \Tilde{H}_p = \frac{1}{\rho^{m-1}} H_p
\end{equation}
is a vector field on the co-sphere bundle $S^{\ast}M$, which is also a manifold with boundary. Then if $\Tilde{H}_p$ satisfies the assumptions of proposition \ref{TrappingConclusion} we can construct a pseudo-differential operator $\Phi_{u/s}$ of order 0 whose principal symbol is a defining function for $\Gamma_u^T$ in $\varphi_{T} (\mathcal{O})$ and non-vanishing outside $\Gamma_{u}$, where $T$ is large enough so that
\begin{equation}
    T > \tau = \frac{\sup_{\mathcal{\Bar{O}}}(F) + \epsilon}{\inf_{\overline{\Gamma_u} \cap \mathcal{O}^c}(H(F))}
\end{equation}
 For example, we can let $\Phi_{u}$ be a quantization of 
 \begin{equation}
     \Tilde{\phi}^{T}_{u} + i \psi
 \end{equation}
 where $\Tilde{\phi}^{T}_{u}$ is an extension of $\phi^T_u$ outside $\varphi_T(\mathcal{O})$ and $\psi$ vanishes on $\varphi_T(\mathcal{O})$ but non-vanishing outside. Furthermore, there is a global function $0 \leq \alpha \leq 1$ so that $H(\alpha) \leq 0$ and vanishing in a neighborhood $\mathcal{U}$ of $\Gamma_u \backslash \Gamma_u^T$:
\begin{equation}
    \alpha|_{\mathcal{U}} = 0 \quad \Gamma_u \backslash \Gamma_U^T \subset \mathcal{U}
\end{equation}
This latter property follows from $\alpha(\varphi_{\tau}(x)) = 0$ for $x \in \partial(\overline{\Gamma_u \cap \mathcal{O}})$. So the statement
\begin{equation}
    \Phi_u v \in H^{s+\alpha}
\end{equation}
 depends on $T$ and $\alpha$, but not on the extension of $\phi^T_u$ chosen outside $\varphi_T(\mathcal{O})$. This is because any such extension is non-vanishing outside $\mathcal{U}$ by assumption hence $\Phi_u$ is elliptic there. Thus, we have a reasonably invariant way of saying that $v$ is more regular outside $\Gamma_u^T$, in the sense that it only depends on a finite propagation time $T$ that we can take arbitrarely large and on a smooth improvement order $\alpha$ chosen. The operators $\Phi_u$ constructed above are the one we will be using in the Fredholm theory.

\end{document}